\documentclass[a4paper,10pt]{amsart}
\usepackage{amssymb,amsfonts,amscd,amsthm,epsfig}

\title{Noncommutative two-dimensional topological
  field theories and Hurwitz numbers for real algebraic curves}
\author{A.Alexeevski $^{1)}$, S.Natanzon $^{2)}$}
 \thanks{This research is partially supported by \\$^{1)}$ 
CRDF grant RM1-2088 \\  $^{2)}$ RFBR-01-01-00739, NSh-1972.2003.1 and 
INTAS-00-0259}

\begin{document}
\newcommand{\glue}{\mathop{\sf glue}\nolimits}
\newcommand{\cut}{\mathop{\sf cut}\nolimits}
\newcommand{\contr}{\mathop{\sf contr}\nolimits}
\newcommand{\Kl}{\mathop{\sf Kl}\nolimits}
\newcommand{\Mb}{\mathop{\sf Mb}\nolimits}
\newcommand{\Cyl}{\mathop{\sf Cyl}\nolimits}
\newcommand{\M}{\mathop{\rm M}\nolimits}
\newcommand{\conv}{\mathop{\sf Contr}\nolimits}
\newcommand{\cor}{\mathop{\rm Cor}\nolimits}
\newcommand{\Hom}{\mathop{\rm Hom}\nolimits}
\newcommand{\Ind}{\mathop{\rm Ind}\nolimits}
\newcommand{\aut}{\mathop{\rm Aut}\nolimits}
\newcommand{\tr}{\mathop{\rm tr}\nolimits}
\newcommand{\ad}{\mathop{\rm ad}\nolimits}
\newcommand{\id}{\mathop{\rm id}\nolimits}
\newcommand{\Ker}{\mathop{\rm Ker}\nolimits}
\newcommand{\Aut}{\mathop{\rm Aut}\nolimits}
\newcommand{\spin}{\mathop{\rm Spin}\nolimits}
\newcommand{\SL}{\mathop{\rm SL}\nolimits}
\newcommand{\GL}{\mathop{\rm GL}\nolimits}
\newcommand{\PSL}{\mathop{\rm PGL}\nolimits}
\newcommand{\PGL}{\mathop{\rm PGL}\nolimits}
\newcommand{\SO}{\mathop{\rm SO}\nolimits}
\newcommand{\Spin}{\mathop{\rm Spin}\nolimits}
\newcommand{\note}{ {\bf Remark.}\quad}
\newtheorem{theorem}{Theorem}[section]
\newtheorem*{definition}{Definition}
\newtheorem{lemma}{Lemma}[section]
\newtheorem{proposition}{Proposition}[section]
\newtheorem{corollary}{Corollary}[section]
\maketitle
\pagestyle{myheadings}
\markboth{Noncommutative extensions of 2D TFTs}
  {Noncommutative extensions of 2D TFTs}

{\bf Abstract}

It is well-known that classical two-dimensional topological field theories 
are in one-to-one correspondence  with commutative Frobenius algebras.
An important extension of  classical two-dimensional topological field theories is provided by open-closed two-dimensional topological field theories. In this paper we extend open-closed two-dimensional topological field theories to nonorientable 
surfaces. We call them Klein topological field theories (KTFT).

We prove that KTFTs bijectively  
correspond to algebras with certain additional structures, called structure algebras. 
Semisimple structure algebras are classified. 
Starting from an arbitrary  finite group, we construct a structure
algebra and prove that it is semisimple.

We define an analog of Hurwitz numbers for real algebraic curves and prove that they are correlators of a KTFT. The structure algebra of this KTFT is the structure algebra of the symmetric group.
\\{\it Key words:} topological field theory, Frobenius algebra, Hurwitz numbers.
\\MSC: 16W, 57M, 81T

\tableofcontents

%%%%%%%%%                

\section{Introduction}

The notion of a  (closed) topological field theory 
was introduced in \cite{ At} (see also \cite{Wit, Seg}).
According to  \cite {At}, a  topological field theory (TFT) in dimension $d$ is defined by the following data:

$(A)$ A finite dimensional vector space $Z(\Sigma)$ for each oriented closed smooth $d$-dimensional manifold $\Sigma$,

$(B)$ An element $Z(M)\in Z(\partial M)$, associated to each oriented smooth $(d+1)$-dimensional manifold (with boundary) $M$,

and by axioms, ensuring the topological invariance of these data.

TFTs  play important role both in pure mathematics and
in mathematical physics.
In particular, closed two-dimensional (2D) topological field theories 
 and their generalizations
are used in string theory and D-brane theory \cite{Laz}, G-flux and K-theories
\cite{Moo1}, WZW,  G/G and CS theories \cite{Gaw}.

In this work we consider not necessarily closed $2D$ TFTs.
For $2D$ closed TFT we have that:

$(A)$ $\Sigma$ is a closed contour and $A=Z(\Sigma)$ is a vector space,
the same for all contours,

$(B)$ $Z(\Omega)\in Z(\partial\Omega)$ is a tensor on $A$ associated to each oriented smooth surface with boundary and depending on topological type of $\Omega$ only.

Contracting  boundary contours of  $\Omega$ to points, we can transform surfaces  with boundary to surfaces without boundary, but with special points. Thus, we can consider a
closed  TFT as a functor from category of oriented surfaces without boundary, but with special poins into a category of tensors on $A$ . (A valency of the tensor is the number of the special poins). The axioms of the classical closed TFT are easily translated to this language.

Replasing in this definition oriented surfaces without boundary by arbitrary compact surfaces, we obtain some generalisation of the closed $2D$ TFTs, that we call a {\it Klein topological field theory} (KTFT). It is the main subject of this paper.

Our KTFT resembles to a lattice TFT, where partition functions for triangulated surfaces are defined instead of
a functor. An outline and a partial algebraic description of a lattice TFT was given in \cite{Kar-Mos}.
A restriction of KTFT to oriented surfaces (possibly with boundaries) is 
equivalent to an open-closed TFT \cite{Laz,Moo2}.

According to \cite{Dij1,Tu} (see also \cite{Ab,Dub}), closed TFTs correspond bijectiely to commutative Frobenius algebras with a nondegenerate 
associative bilinear form. In \cite{Laz} an algebraic object, corresponding to an open-closed TFT is constructed . In \cite{Moo2} this object is represented as an algebra with additional structures. (A certain algebraic interpretation of this algebra we call {\it Lazaroiu-Moore algebra}.) In \cite{Moo2} it
is claimed (without a proof) that the correspondence between open-closed TFTs and these algebras is a bijection.

In this paper we define certain algebras with additional structures ({\it structure algebras}) and prove that Klein topological field theories  bijectively correspond to structure algebras.

A structure algebra is an associative (typically noncommutative) algebra
endowed with an invariant scalar product and three features: 
a decomposition into direct sum of a commutative subalgebra and an ideal 
(reflecting two kinds of special points, interior and boundary ones),
an involutive antiautomorphism (reflecting the change of the 
local orientation at a special point), and an element $U$ 
(reflecting nonorientability of a surface). 
It is worth to note that KTFT gives an additional structure even for
the category of orientable surfaces; it is  an involutive antiautomorphism.
This involutive antiautomorphism is missing in closed and open-closed TFTs, 
since it is not needed when an orientation of the surface is fixed.

A structure algebra without the involutive antiautomorphism and  without the element $U$ is equivalent to a Lazaroiu-Moore algebra. It follow from our results that open-closed TFTs correspond
bijectively  to Lazaroiu-Moore algebras. Moreover, we prove that each  semisimple Lazaroiu-Moore algebra has an extension to a structure algebra, and the number of such extensions is finite. In physics, semisimple KTFTs correspond  to  massive systems.
Thus, each massive open-closed TFT is extended to a KTFT, and the number of such extensions is finite.

In the last part of the paper we apply the Klein topological field theories 
to study (generalized) Hurwitz numbers. 

More than 100 years ago Hurwitz \cite{Hur} raised the following 
problem. Let $\Omega$ be a complex algebraic curve of genus $g$, 
$p_1,\dots,p_m\in \Omega$ be pairwise distinct points and 
$\alpha_1,\dots,\alpha_m$ be partitions of $n$
(decompositions of the number $n$ into unordered sums of positive
integers). 
Denote by $S(\alpha_1,...,\alpha_m)$
the set of classes of birational equivalence
of morphisms $\pi:P\to \Omega$ of order $n$ from complex algebraic curves to $\Omega$, having critical values
$p_1,\dots,p_m$ such that local degrees of $\pi$ at points 
$\pi^{-1}(p_i)$ are $\alpha_i(i=1,...,m)$.
The problem is to compute the sum $\sum_{\pi\in S(\alpha_1,...,
\alpha_m)}\frac{1}{|\Aut(\pi)|}$, where $|\Aut(\pi)|$ is 
the order of the group of automorphisms of a covering
$\pi$. These sums are called (classical) {\it Hurwitz numbers}.

The generating function for Hurwitz numbers is the partition function for 
the $2D$ supersymmetric Yang-Mills theory \cite{Co-Mo-Ra}. These numbers are
connected with string theory \cite{Gr-Tay}, mirror symmetry \cite{Dij},
theory of singularities \cite{Arn}, matrix models \cite{Kos-Sta-Wyn}
and integrable systems \cite{Ok-Pan}.
There exist  combinatorial, differential and integral formulas
for the simplest of Hurwitz numbers. See, for example, \cite{Arn, Go-Jac2, Eke-Lan}.

According to \cite{Dij1}, the classical Hurwitz numbers are correlators for a closed TFT. The Frofenius algebra of this TFT is the Frobenius algebra of the symmetric group. Such approach gives an effective method for calculation of the classical Hurwitz numbers \cite{Dij}.

In this paper we extend the definition of Hurwitz numbers to real algebraic curves.
Moreover, our (generalized) Hurwitz numbers include classical Hurwitz numbers as a particular case.  Hurwitz numbers of real algebraic curves are important 
for string theory and supersymmetric Yang-Mills theory.\cite{Co-Mo-Ra}. Some special cases of them  were investigated in \cite{Med, Bar}.

We prove that the generalized Hurwitz numbers are correlators for a KTFT. Moreover, the structure algebra of this KTFT is the structure algebra of the symmetric group. This gives an effective method for the calculation of generalized Hurwitz numbers.

The paper is organized as follows. 

In section 2 we give an algebraic background for a Klein topological field theory, developed in section \ref{sec4} . We give  (subsection 2.1) an axiomatic definition of structure algebras and describe them by a set of structural constants and relations. 
(These relations will be interpreted geometricaly in sections \ref{sec3} and \ref{sec4}.)
In subsection \ref{2s2} we classify semisimple structure algebras. In subsection \ref{2s3} for any finite group we construct a semisimple structure 
algebra, which typically is not commutative.

Section 3 contains a geometrical background for the Klien topological field theory, developed in section \ref{sec4}. 
We consider  surfaces with  boundaries and
special points. The simplest classes of such surfaces (trivial, basic and simple
surfaces in the nomenclature of subsections 3.1 and 3.3 ) we exploit in section \ref{sec4} for constructing of sructure algebras.

We consider also (subsections 3.2 and 3.3) systems of nonintersecting generic cuts of surfaces, and in particular, complecte cut systems cuting surfaces into    basic surfaces. The simple surfaces generate elementary shifts of complecte cut systems. 

A central theorem of this section (theorem \ref{th4.1}, subsection 3.4) claims that any complete cut system of
a surface can be transformed to any other complete cut system of the same surface by 
elementary shifts. 
 This theorem playes an essential role in the proof of the main theorem of the paper (theorem \ref{MMain}) about the equivalence between structure algebras and Klein topological field theories.

In section 4 we  
use results of sections \ref{sec1} and \ref{sec3} in order to define a Klien topological field theory (subsection 4.1), reformulate it 
in terms of systems of correlators as it is usually done 
in physical literature (subsections 4.2), and prove the main theorem  \ref{MMain} , which states 
the correspondence between KTFT and structure algebras (subsection 4.3, 4.4).
As corollary we give an analog of this theorem for open-closed topological field
theoris and prove that any massive open-closed topological field
theory can be extended to a Klein topological field theory (subsection 4.5).

In section 5 the coverings over stratified surface are considered. It is shown (subsection 5.1) that 
singularities over boundary special points are classified by 'dihedral Yang diagrams', which in turn correspond to conjugacy classes of pairs of involutions in a symmetric
group $S_n$. 

Classical Hurwits numbers are generalized to coverings over 
stratified surfaces. A Klien topological field theory is associated (subsections 5.2, 5.3)
with these generalized Hurwitz numbers. This KTFT is called Hurwits topological field 
theory. 

It is proved that Hurwits topological field theory corresponds to the structure
algebra associated the symmetric group $S_n$. It allows us to obtain
the expressions (subsection 5.4) for generalized Hurwits numbers via structural constants of the structural
algebra.

 \textbf{Acknowledgment.} The final version of this paper was written during the second 
 author stay at Max-Planck-Institut f\"{u}r Mathematik in Bonn. He would like 
 to thank this institution for its support and hospitality.

%%%%%%%%%%%%%%%%%%%%%%%%%%%%%%%%%%%%%%%%%%

        	\section{Structure algebras}
\label{sec1}

This section gives an algebraic background for the Klein topological field theory
developed in section \ref{sec4} . We introduce (subsection 2.1) a purely algebraic object, which is an algebra 
endowed with additional structures ('structure algebra'). In 
section \ref{sec4}
we will prove that structure algebras bijectively correspond to 
Klien topological field theories. In addition to an axiomatic definition of the 
structure algebra we describe them by the set of structural constants and relations. 
These relations will be interpreted geometricaly in sections \ref{sec3} and \ref{sec4}.

Despite the complexity of the formal definition, semisimple structure algebras seems to be 
rather a handy object. In subsection \ref{2s2} we classify all of them. A semisimple 
structure algebra is the sum $H=A\oplus B$, where $B$ is the direct sum of full 
matrix algebras, and $A$ is a commutative subalgebra isomorphic to the direct sum 
of one dimensional algebras. Additional 
structures are: an isomorphism between the center of $B$ and a subalgebra of $A$, 
an invariant scalar product, an involutive antiauthomorphism and a certain 
element $U\in A$. 

In subsection \ref{2s3} for any finite group we construct a semisimple structure 
algebra, which typically is not commutative.
In section \ref{sec5} we prove that the structure algebras of symmetric groups
generate Hurwitz numbers.

	\subsection{Definition of a structure algebra}
\label{secDSA} All vector spaces (algebras) in the paper are vector spaces 
(algebras) over complex numbers.
Let $X$ be a finite dimensional  associative algebra with unit 
 endowed with a symmetric invariant nondegenerate scalar product 
$(.,.):X\otimes X\to\mathbb C$, i.e., $(x,y)=(y,x)$,  $(xy,z)=(x,yz)$
and $(x,a)=0$ for an arbitrary $x\in X$ implies $a=0$.

Fix a basis $E$ of $X$. We denote by $\alpha,\alpha',\dots\in E$
elements of the basis. We use the same letters as indices of tensors.
Thus, we denote by $F_{\alpha',\alpha''}$ the matrix
of the scalar product in basis
$E$ and we denote by $F^{\alpha',\alpha''}$ the tensor of the dual form.
Hence, 
$F_{\alpha',\alpha''}F^{\alpha'',\alpha'''}=\delta_{\alpha'}^{\alpha'''}$.
As usual,  repeating indices mean summation.

Denote by $\widehat{K}_X$ an element 
$F^{\alpha',\alpha''}\alpha'\otimes \alpha''$
of the tensor product $X\otimes X$ and by $K_X$
Casimir element. By definition, $K_X=F^{\alpha',\alpha''}\alpha'\alpha''\in X$.
By standard arguments, one can show that elements $\widehat{K}_X$  and $K_X$ 
do not depend on the choice of a basis.

Denote by $V_{K_X}$ the operator
$x\mapsto F^{\alpha',\alpha''}\alpha'x\alpha''$.

By  $x\mapsto x^*$ we denote an involutive antiautomorphism  of $X$,
i.e., $(x^*)^*=x$ and  $(xy)^*=y^*x^*$. 
Denote by $K_{X,*}$  twisted Casimir element
$F^{\alpha',\alpha''}\alpha'\alpha''{^*}$. Obviously,
it coincides with Casimir element of twisted scalar product 
$(x,y)_*=(x,y^*)$.

\begin{definition} A structure algebra 
$\mathcal H=\{H=A\oplus B, (.,.), x\mapsto x^*,  U\}$ 
is a finite dimensional associative algebra  $H$  
endowed with

a decomposition $H=A\oplus B$ of $H$ into the direct sum of two vector spaces; 

a symmetric invariant scalar product  $(.,.): H\otimes H\to\mathbb C$;

an involutive antiautomorphism $H\to H$, denoted by $x\mapsto x^*$;

an element $U \in A$, 
\\such that the following axioms hold:

$1^\circ$ $A$ is a subalgebra belonging to the center of algebra $H$; 
algebra $A$ has unit $1_A\in A$ and $1_A$ is also the unit of algebra $H$;

$2^\circ$ $B$ is a two-sided ideal of $H$ (typically noncommutative); algebra 
$B$ has a unit $1_B\in B$;

$3^\circ$ restrictions $(.,.)|_A$ and $(.,.)|_B$ are 
nondegenerate scalar products on algebras $A$ and $B$ resp. 

$4^\circ$  $(V_{K_B}(b_1),b_2)=(\widehat{K}_A,b_1\otimes b_2)$
for arbitrary $b_1,b_2\in B$ 
(this axiom reflects Cardy  relation \cite{Laz});

$5^\circ$ an involutive antiautomorphism preserves the decomposition 
$H=A\oplus B$ and the form $(.,.)$ on $H$, i.e.,
$A^*=A$, $B^*=B$, $(x^*,y^*)=(x,y)$;

$6^\circ$ $U^2=K_{A,*}$; 

$7^\circ$ $(U,b)=(K_{B,*},b)$ for any $b\in B$;

%$8^\circ$ $U^3=UK_A$; 

$8^\circ$ $(aU)^*=aU$ for any $a\in A$. 
\end{definition}        

\note Forms $(.,.)|_A$ and $(.,.)|_B$ are nondegenerate and 
we use them to raise or lower indices in tensors. 
Bilinear form $(.,.)$ on $H$ is not assumed to be 
nondegenerate. Linear subspaces $A$ and $B$ are not assumed to be 
orthogonal. Clearly, if an element $a\in A$ is orthogonal
to $B$ then $aB=0$. Therefore, if
$A$ is orthogonal to $B$ then $AB=0$.  

\medskip

\begin{definition}  Structure algerbras 
$\mathcal H=\{H=A\oplus B, (.,.), x\mapsto x^*,  U\}$ and 
$\mathcal H'=\{H'=A'\oplus B', (.,.)', x\mapsto x^{*'},  U'\}$ 
are called isomorphic if there existsts an isomorphism
$\varphi:H\to H'$ such that $\varphi(A)=A'$, $\varphi(B)=B'$, 
$(\varphi(x),\varphi(y))'=(x,y)$, $\varphi(x^*)=\varphi(x)^{*'}$, $\varphi(U)=U'$.    
\end{definition}

Fix a structure algebra  
$\mathcal H=\{H=A\oplus B, (.,.), x\mapsto x^*,  U\}$. 
Let $a$ be an element of $A$. Then
formula $(\phi(a),b)=(a,b)$ defines  an element $\phi(a)\in B$
because  bilinear form $(.,.)|_B$ is nondegenerate.

\begin{lemma} \label{2.1}
The mapping $\phi:A\to B$ is a homomorphism of algebra $A$
into the center of the algebra $B$  such that for any $a\in A$,
$b\in B$ the equality $ab=\phi(a)b$ holds.
\end{lemma}

\begin{proof} By invariance of the form $(.,.)$ we have
$(ab,b')=(a,bb')=(\phi(a),bb')=(\phi(a)b,b')$.
Therefore, $ab=\phi(a)b$, which means that 
$\phi$ is a homomorphism.
By definition, $A$ lies in the center of $H$. Hence
 $(ab,b')=(ba,b')=(b,ab')=(b,\phi(a)b')=(b\phi(a),b')$.
Therefore, $\phi(a)b=b\phi(a)$ and $\phi(a)$ belongs to
the center of $B$.
\end{proof}

We shall reformulate the definition of a structure algebra
$\mathcal H=\{H=A\oplus B, (.,.), x\mapsto x^*,  U\}$ in a
coordinate form. Fix basis $E_A$ of $A$ and $E_B$ of $B$. 
We denote by $\alpha,\alpha_1,\dots$ the elements of $E_A$ and
by $\beta,\beta_1,\dots$ the elements of $E_B$.
Indices of tensors are denoted below
by the same letters $\alpha,\alpha_1,\dots;\beta,\beta_1,\dots$.

Define tensors:

(1) $F_{\alpha_1,\alpha_2}=(\alpha_1,\alpha_2)$;

(2) $F_{\beta_1,\beta_2}=(\beta_1,\beta_2)$;

(3) $R_{\alpha,\beta}=(\alpha,\beta)$;

(4) $S_{\alpha_1,\alpha_2,\alpha_3}= (\alpha_1\alpha_2,\alpha_3)$; 

(5) $T_{\beta_1,\beta_2,\beta_3}=(\beta_1\beta_2,\beta_3)$;

(6) $R_{\alpha,\beta_1,\beta_2}=(\alpha\beta_1,\beta_2)$;

(7) $I_{\alpha_1,\alpha_2}=(\alpha_1^*,\alpha_2)$;

(8) $I_{\beta_1,\beta_2}=(\beta_1^*,\beta_2)$;

(9) $D_{\alpha}=(U,\alpha)$;

(10) $J_\alpha=(1_A,\alpha)$;

(11) $J_\beta=(1_B,\beta)$.

Tensors (1)-(11) completely define a structure algebra. Indeed, tensors
(1)-(3) define the bilinear form, tensors (4)-(6) define the multiplication
in algebra $H$, tensors (7)-(8) define the involutive antiautomorphism,
tensor (9) defines the element $U$ and tensors (10), (11) define units
$1_A\in A$ and $1_B\in B$. 

We call tensors (1)-(11) structure constants of  structure algebra
$\mathcal H$ in basis $E_A$, $E_B$.

Let us write down exact formulas.
Suppose  we are given by two linear spaces $A$ and $B$ with  basis 
$E_A$ of $A$, $E_B$ of $B$  and by {\it arbitrary} tensors
$F_{\alpha_1,\alpha_2}$, $\dots$, $J_\beta$ that have the same type 
as tensors from the left-hand sides of equalities (1) - (11) and are denoted
by the same letters.

Assume $F_{\alpha_1,\alpha_2}$ and $F_{\beta_1,\beta_2}$ are 
symmetric and nondegenerate. Using $F_{\alpha_1,\alpha_2}$ 
and $F_{\beta_1,\beta_2}$
and dual tensors we will raise or lower indices of tensors. In the case of 
asymmetric tensors we always raise the last index. 

Let $H=A\oplus B$. 
Define bilinear form on $H$  by formulas (1)-(3).
Define a multiplication in $H$ as follows:
$\alpha_1\alpha_2=S_{\alpha_1,\alpha_2}^\alpha \alpha$; 
$\beta_1\beta_2=T_{\beta_1,\beta_2}^\beta \beta$;
$\beta\alpha=\alpha\beta=R_{\alpha,\beta}^{\beta'}\beta'$.
By definition,  $A$ is a subalgebra and  $B$ is an ideal of $H$.
Define a linear map $x\mapsto x^*$ as follows:
$\alpha^*=I_{\alpha}^{\alpha'}\alpha'$;
$\beta^*=I_{\beta}^{\beta'}\beta'$.
Define elements $U\in A,1_A$, $1_B\in B$ by formulas $U=D^\alpha\alpha$,
$1_A=J^\alpha\alpha$, $1_B=J^\beta\beta$.
Denote by $S_{\alpha_1,\alpha_2,\alpha_3,\alpha_4}$
the contraction 
$S_{\alpha_1,\alpha_2}^{\alpha'}S_{\alpha',\alpha_3,\alpha_4}$
and by  $T_{\beta_1,\beta_2,\beta_3,\beta_4}$
the contraction
$T_{\beta_1,\beta_2}^{\beta'}T_{\beta',\beta_3,\beta_4}$.

\begin{lemma} \label{2.2} A set of data
$\mathcal H=\{ H=A \oplus B, (.,.), x\mapsto x^*,  U\}$ 
is a structure algebra if and only if
\\ (1) $F_{\alpha_1,\alpha_2}$, $F_{\beta_1,\beta_2}$ are 
symmetric nondegenerate tensors;
\\ (2) $S_{\alpha_1,\alpha_2,\alpha_3}$ and
$S_{\alpha_1,\alpha_2,\alpha_3,\alpha_4}$
are symmetric tensors;
\\ (3) tensors $T_{\beta_1,\beta_2,\beta_3}$ and 
$T_{\beta_1,\beta_2,\beta_3,\beta_4}$
are invariant under cyclic permutations; 
\\ (4) $R_{\alpha,\beta_1,\beta_2}=
R_\alpha^{\beta'}T_{\beta',\beta_1,\beta_2}$;
\\ (5) $R_{\beta}^{\alpha'}S_{\alpha',\alpha_1,\alpha_2}=
	R_{\alpha_1}^{\beta'}R_{\alpha_2}^{\beta''}
	T_{\beta',\beta'',\beta}$;
\\ (6) $R_{\alpha,\beta_1,\beta_2}=R_{\alpha,\beta_2,\beta_1}$; 
\\ (7) $R_{\alpha',\beta_1}F^{\alpha',\alpha''}R_{\alpha'',\beta_2}=
       T_{\beta_1,\beta'}^{\beta''}T_{\beta'',\beta_2}^{\beta'}$
\\ (8) $I_{\alpha_1}^{\alpha'}I_{\alpha',\alpha_2}=F_{\alpha_1,\alpha_2}$,
$I_{\beta_1}^{\beta'}I_{\beta',\beta_2}=F_{\beta_1,\beta_2}$
\\ (9) $I_{\alpha_1,\alpha_2}=I_{\alpha_2,\alpha_1}$,
$I_{\beta_1,\beta_2}=I_{\beta_2,\beta_1}$,          
$I_{\alpha}^{\alpha'}R_{\alpha',\beta}=I_{\beta}^{\beta'}R_{\alpha,\beta'}$
\\ (10) $I_{\alpha_1}^{\alpha'}I_{\alpha_2}^{\alpha''}I_{\alpha_3}^{\alpha'''}
S_{\alpha''',\alpha'',\alpha'}=S_{\alpha_1,\alpha_2,\alpha_3}$,              
$I_{\beta_1}^{\beta'}I_{\beta_2}^{\beta''}I_{\beta_3}^{\beta'''}
T_{\beta''',\beta'',\beta'}=T_{\beta_1,\beta_2,\beta_3}$,
\\ (11) $S_{\alpha,\alpha',\alpha''}D^{\alpha'}D^{\alpha''}=
       S_{\alpha,\alpha',\alpha''}I^{\alpha',\alpha''}$;
\\ (12) $D^{\alpha'} R_{\alpha',\beta}=I^{\beta',\beta''}T_{\beta',\beta'',\beta}$
%\\(13) $S_{\alpha,\alpha',\alpha'',\alpha'''}D^{\alpha'}D^{\alpha''}D^{\alpha'''}=
%S_{\alpha,\alpha',\alpha'',\alpha'''}D^{\alpha'}F^{\alpha'',\alpha'''}$
\\ (14) $D^{\alpha'}S_{\alpha',\alpha_1,\alpha_2}=
I_{\alpha_1}^{\alpha'}D^{\alpha''}S_{\alpha',\alpha'',\alpha_2}$
\\ (15) $J^{\alpha'}S_{\alpha',\alpha_1,\alpha_2}=F_{\alpha_1,\alpha_2}$,
$J^{\alpha'}R_{\alpha',\beta_1,\beta_2}=F_{\beta_1,\beta_2}$
\\ (16) $J^{\beta'}T_{\beta',\beta_1,\beta_2}=F_{\beta_1,\beta_2}$
\end{lemma}

\begin{proof} We shall list below the correspondence between relations
(1)-(12) and
axioms of a structure algebra. All proofs are by direct calculations.

Relation (1) coincides with axiom $3^\circ$.

Relation (2) is equivalent to the claim '$A$ is an associative
commutative algebra endowed with an invariant scalar product'.

Relation (3) is equivalent to the claim 
'$B$ is an associative algebra endowed with invariant scalar product'.

Define a map $\phi:A\to B$ by  $\phi(\alpha)=R_\alpha^\beta\beta$.
Clearly $(\alpha,\beta)=(\phi(\alpha),\beta)$.

Relation (4) is equivalent to the claim  $\alpha\beta=\phi(\alpha)\beta$.

Relation (5) is equivalent to the claim '$\phi$ is a homomorphism 
of algebras'.

Relation (6) is equivalent to the claim $\phi(\alpha)\beta=\beta\phi(\alpha)$.

It is easy to show that if $\phi$ is a homomorphism and its image lies in 
the center of algebra $B$ then $H$ is an associative algebra and the
form $(.,.)$ is invariant and vice versa.

Relation (7) is equivalent to axiom $4^\circ$.

Relation (8) is equivalent to the claim '$x\mapsto x^*$ is an involution'.

Relation (9) is equivalent to the claim 'the involution  
$x\mapsto x^*$ preserves the bilinear form  $(.,.)$'.

Relations (10) is equivalent to the claim
'the involution $x\mapsto x^*$ preserving the form $(.,.)$ is an 
antiautomorphism'. Indeed, two relations from (10) are the direct reformulation 
of this fact for subalgebras $A$ and $B$.            
Let us prove that
$(\alpha\beta)^*=\beta^*\alpha^*$. 
We have $(\phi(\alpha),\beta)=(\alpha,\beta)$.
Therefore, $(\phi(\alpha^*),\beta)=(\alpha^*,\beta)=(\alpha,\beta^*)$. 
Analogously, $(\phi(\alpha)^*,\beta)=(\phi(\alpha),\beta^*)=(\alpha,\beta^*)$. 
We obtain that $\phi(\alpha^*)=\phi(\alpha)^*$.
Hence, $(\alpha\beta)^*=(\phi(\alpha)\beta)^*=\beta^*\phi(\alpha)^*=
\beta^*\phi(\alpha^*)=\beta^*\alpha^*$.

Relation (11) is equivalent to axiom $6^\circ$.

Relation (12) is equivalent to axiom $7^\circ$.

%Relation (13) is equivalent to axiom $8^\circ$.

Relation (14) is equivalent to axiom $8^\circ$.

Relations (15), (16) are equivalent to the claims 
'$1_A$ is a unit of $H$' and '$1_B$ is a unit of $B$'.

\end{proof}

We call (1) - (16) {\it relations} for structure 
constants of a structure algebra.

\subsection{Semisimple structure algebras}
\label{2s2}
A structure algebra $\mathcal H$ is called {\it semisimple} if the
algebra $H$ is semisimple.

Let $\mathcal H$ be a finite dimensional complex semisimple structure algebra. 
Then both subalgebras $A$ and $B$ are semisimple algebras.
Indeed, $B$ is an ideal of $H$ and $A$ is isomorphic to the factor $H/B$.

By Wedderburn theorem, $B$ is isomorphic to the direct sum 
of matrix algebras, $B=\oplus_{i=1}^{k} \M_{n_i}$. 
$A$ is commutative semisimple algebra; therefore, $A$ is isomorphic to
the direct sum of one-dimensional algebras,  
$A=\overbrace{\mathbb C\oplus\dots\oplus\mathbb C}^m$.
Clearly, as an abstract algebra  $H$ is isomorphic to the direct sum 
$A\oplus B$. 
Despite this fact, the decomposition $H=A\oplus B$
that is included in the set of data $\mathcal H$ typically does  not coincide 
with a decomposition into the direct sum of two-sided ideals.

Let us write down formulas for  the multiplication $A\times B\to B$ and 
for the invariant bilinear form $(.,.):H\otimes H\to\mathbb C$.
There is a uniquely defined, up to permutations, complete system
$e_1,\dots e_m$ of orthogonal  idempotents in $A$ and 
$A=\mathbb Ce_1\oplus\dots\oplus\mathbb Ce_m$.
Orthogonal idempotents are orthogonal with respect to the invariant
scalar product. Therefore, we have $(e_i,e_j)=\lambda_i\delta_{i,j}$,
$\lambda_i\ne 0$.

Similarly, components $\M_{n_i}$ of algebra $B=\oplus_{i=1}^{k} \M_{n_i}$
are orthogonal with respect to any invariant scalar product on $B$.
All invariant scalar products on $\M_{n_i}$ are proportional 
to each other. Hence, if $X,Y\in \M_{n_i}$ then $(X,Y)=\mu_i\tr(XY)$,
$\mu_i\ne 0$. 

Denote  by $E_i$ the unit matrix of  $\M_{n_i}$. Clearly,
elements $\{E_i |i=1,\dots,k\}$ form a basis of the center of 
algebra $B$ and they are orthogonal idempotents. 
By lemma \ref{2.1}, morphism  $\phi:A\to B$ is a homomorphism into the
center of $B$. Therefore, the image 
of a complete system $\{e_i |i=1,\dots,m\}$ of orthogonal idempotents is 
a system of orthogonal idempotents in the center of $B$.  Hence,
$\phi(e_i)=\sum_{j\in N_i}E_j$, where $N_i\subset \{1,\dots,k\}$
and if $i\ne i'$ then $N_i\cap N_{i'}=\emptyset$. Thus,
$$e_ib=be_i=(\sum_{j\in N_i}E_j)b \  \text{for any}\ b\in B\ \ (2.1)$$
and this formula completely describes the multiplication between elements
of $A$ and $B$.

Denote by $E_{s,i,j}\in M_{n_s}$ the $n_s\times n_s$ matrix with all 
elements equal to zero except (i,j) element which is equal to $1$.
We get
$$(e_i,E_{s,i,j})=(\phi(e_i),E_j)=
	\left\{ \begin{array}{l}
			\mu_s\delta_{i,j}\mbox{ if $s\in N_i$}\\
			0 \mbox{ if $s\notin N_i$}\\
	  	\end{array} 
        \right.\ \ (2.2)$$

Conversely, let $A=\mathbb C^m$ and $B=\oplus_{i=1}^{k} \M_{n_i}$
be algebras with invariant scalar products defined by 
constants $\lambda_i$, $i=1,\dots, m$ and
$\mu_j$, $j=1,\dots,k$ resp.
Define a multiplication between $A$ and $B$ 
by formula 2.1 and the scalar product between $A$ and $B$ by formula 2.2. 
Then axioms $1^\circ$-$3^\circ$ of a structure
algebra are satisfied. 
Let us denote such an algebra endowed with an invariant bilinear form by
$\widetilde{H}=A\oplus B$.

\begin{lemma}  Cardy axiom $4^\circ$ is satisfied in algebra $\widetilde H$
if and only if 
\\ (1) for $i=1,\dots,m$ a subset $N_i$
of the set $\{1,\dots,k\}$ is either
empty set or contains one element;
\\ (2) $\{1,\dots,k\}=\sqcup_{i}N_i$;
\\ (3) if $N_i=\{j\}$ then $\lambda_i=\mu_j^2$. 
\end{lemma}

\begin{proof} Let us compute element
$\widehat{K}_A$ and transformation $V_{K_b}$.
Evidently, $\widehat{K}_A=\sum_{i=1}^m \frac{1}{\lambda_i}e_i\otimes e_i$
and  $K_A=\sum_{i=1}^m \frac{1}{\lambda_i}e_i$.

Evidently, all matrices $E_{s,i,j}$ form a basis of algebra $B$.
We obtain by direct calculations that
$K_B=\sum_{s=1}^{k}\frac{n_s}{\mu_s} E_s$ and
$V_{K_B}(E_{s,i,j})=\delta_{i,j}\frac{1}{\mu_s}E_s$.
Therefore, putting  $b_1=E_{s,i,j}$, $b_2=E_{s',i',j'}$ in axiom $4^\circ$ we
obtain that: 

a) left side $= \delta_{s,s'}\delta_{i,j}\delta_{i',j'}\frac{1}{m_s}\mu_s=
 \delta_{s,s'}\delta_{i,j}\delta_{i',j'}$;

b) right side $=\left\{
      \begin{array}{l}
		\frac{\mu_s\mu_{s'}}{\lambda_l}\delta_{i,j}\delta{i',j'}
					\mbox{ if $s,s'\in N_l$}\\
		0 \mbox{ if $s$ and $s'$ do not lie in the same set $N_l$}\\
		\end{array} \right.
	   $

The proof is  completed by comparing these formulas.

\end{proof}

Assume that axiom $4^\circ$ holds for an algebra $\widetilde H$. Then
we can and will reorder components of $A$
and $B$ in such a way that 

a) for $i=1,\dots,k$ we have $\phi(e_i)=E_i$, where $E_i$ is the unit
matrix of the component $\M_{n_i}$;

b) for $i=k+1,\dots,m$ we have $\phi(e_i)=0$.

Let us check axioms related  to an involutive antiauthomorphism
$x\mapsto x^*$ and an element $U$ for algebra $\widetilde H$. 
Suppose $x\mapsto x^*$ is an involutive transformation
satisfying axiom $5^\circ$.

Obviously, involutive antiautomorphism $x\mapsto x^*$ can either fix 
an idempotent $e_s$ or permute two idempotents $e_{s'}$ and $e_{s''}$. 
If the latter occurs then $\lambda_{s'}=\lambda_{s''}$ because this 
antiautomorphism preserves the scalar product. Thus, the involution 
of the set $\{1,\dots,m\}$ is induced. We denote it 
by the same sign $s\mapsto s^*$. 

Analogously, the involution $x\mapsto x^*$ either preserves a component 
$\M_{n_s}$ or permutes two components $\M_{n_s}$ and $\M_{n_s'}$ of 
algebra $B$. If the latter occurs then we have $n_s=n_{s'}$ and  
$\mu_s=\mu_{s'}$. Thus, the involution of the set $\{1,\dots,k\}$ 
is induced. This involution coincides with the restriction to 
$\{1,\dots,k\}$ of previously defined involution $s\mapsto s^*$ of 
the set $\{1,\dots,m\}$ because transformation $x\mapsto x^*$ preserves
bilinear form and $e_i$ is orthogonal to $M_{n_j}$ if and only if $i\ne j$.

Obviously, the restriction of the transformation $x\mapsto x^*$ to $A$
is completely defined by the involution $s\mapsto s^*$ of the set 
$\{1,\dots,m\}$. There is an ambiguity for the restriction of the 
involutive transformation to subalgebra $B$.
The possibilities are well known. We describe them below without proofs.

If $s\le k$ and $s^*\ne s$ then the restriction of $x\mapsto x^*$ to 
$\M_{n_s}\oplus\M_{n_{s^*}}$ is conjugated to the transformation 
$\{X,Y\}\to \{Y',X'\}$, where $X'$ denotes the transpose of $X$.Here
$\{X,Y\}$ is an element of $\M_{n_s}\oplus\M_{n_{s^*}}$, i.e., two matricies.
Thus, changing the basis of $B$ if necessary, we have $\{X,Y\}^*=\{Y',X'\}$.

If $s\le k$ and $s^*=s$ then an involutive antiautomorphism of $\M_{n_s}$ 
is conjugated to one of two canonical antiautomorphisms. 

The first of them is associated with a symmetric
bilinear form and coincides with the transpose of a matrix. 

The second of them is associated with a nondegenerate skew-symmetric
bilinear form. Hence, it exists only if $n_s=2r$.
Let us identify the set $\{1,\dots,n_s\}$ with the set 
$\mathbb Z_2\times\{1,\dots,r\}$, 
where $\mathbb Z_2$ is the group with two elements $\{0,1\}$. 
If $i=(\epsilon,i')$ then put $\epsilon(i)=\epsilon$ and 
$i^\tau=(\epsilon+1,i')$.
Define a transformation $\tau:\M_{n_s}\to \M_{n_s}$ by formula
$(E_{i,j})^\tau=(-1)^{\epsilon(i)+\epsilon(j)}E_{j^\tau,i^\tau}$.
One can check directly that $\tau$ is the involutive antiautomorphism preserving 
the invariant scalar product.

Therefore, if $s^*=s$ then in an appropriate basis we have either
$E_{s,i,j}^*=E_{s,j,i}$ or $E_{s,i,j}^*=(-1)^{\epsilon(i)
+\epsilon(j)}E_{s,j^\tau,i^\tau}$.
In order to distinguish these cases let us introduce an invariant
$\nu=\nu(s)\in{\pm 1}$ and put $\nu(s)=1$ in the former case, $\nu(s)=-1$ 
in the latter case.

Denote by $P$ a set of fixed points of involution $s\mapsto s^*$ and put 
$P_0=P\cap\{1,\dots,k\}$. 

As a result of the above considerations, we
obtain the following lemma.

\begin{lemma} Let $x\mapsto x^*$ be an involutive antiautomorphisms 
of algebra $\widetilde H$. Suppose,  $x\mapsto x^*$ satisfies axiom 
$5^\circ$. Then antiautomorphism $x\mapsto x^*$ induces 
\begin{itemize}
\item[$\circ$] involution $s\mapsto s^*$ of set $\{1,\dots,m\}$ 
such that 
$(\{1,\dots,k\})^*=\{1,\dots,k\}$ and $\mu_s=\mu_{s^*}$, $\lambda_s=\lambda_{s^*}$;
\item[$\circ$] numbers $\nu(s)\in\{\pm 1\}$ for  $s\in P_0$ where 
$P_0=\{s|s\le k, s^*=s\}$
\end{itemize} 
Antiautomorphisms $s\mapsto s^*$ up to inner automorphism of 
$\widetilde H$ are classified by pairs  $(s\mapsto s^*;\nu(s))$.
\end{lemma}

Let us fix additionally an involutive
antiautomorphism satisfying axiom $5^\circ$ and check the possibilities of
choosing an element $U\in A$ satisfying axioms $6^\circ-8^\circ$.

\begin{lemma} An element $U\in A$ satisfies axioms $6^\circ-8^\circ$
if and only if
$U=\sum_{i\in P_0}\frac{\nu_i}{\mu_i}e_i+\sum_{j\in P\setminus P_0}x^je_j$,
where $(x^j)^2=\frac{1}{\lambda_j}$. 
\end{lemma}

\begin{proof} Let $U=\sum_{i=1}^m x^ie_i$ be an element of $A$ satisfying 
axioms $6^\circ-8^\circ$. Then 
$U^2=\sum_{i=1}^m (x^i)^2e_i$. Twisted Casimir element
$K_{A,*}$ is equal to $\sum_i \frac{1}{\lambda_i}e_ie_{i^*}=
\sum_{i\in P}\frac{1}{\lambda_i}e_i$. Therefore,
 $U=\sum_{i\in P} x^ie_i$ and $(x^i)^2=\frac{1}{\lambda_i}$ (axiom $6^\circ$)
Evidently,
$(U,E_{s,i,j})=\left\{\begin{array}{l}
                     x_s\mu_s\delta_{i,j}\mbox{ if $s\in P_0$}\\
                                  0  \mbox{ if $s\notin P_0$}\\
				   \end{array}
				\right. $

Compute element $K_{B,*}$. By definition, 
$K_{B,*}=\sum_{s=1}^k \sum_{i,j}\frac{1}{\mu_s}E_{s,i,j}E_{s,j,i}^*$. 
If $s^*\ne s$ then $E_{s,i,j}E_{s,j,i}^*=0$, hence only summands 
corresponding to $s\in P_0$ are nonzero. 

If $\nu_s=1$ then 
$\sum_{i,j}\frac{1}{\mu_s}E_{s,i,j}E_{s,j,i}^*=
\sum_{i,j}\frac{1}{\mu_s}E_{s,i,j}E_{s,i,j}=
\frac{1}{\mu_s}\sum_{i}E_{s,i,i}=\frac{1}{\mu_s}E_s=\frac{1}{\mu_s}\nu_sE_s$. 

If $\nu_s=-1$ then 
$\sum_{i,j}\frac{1}{\mu_s}E_{s,i,j}E_{s,j,i}^*=
\sum_{i,j}\frac{1}{\mu_s}E_{s,i,j}
(-1)^{\epsilon(j)+\epsilon(i)}E_{s,i^*,j^*}=
\frac{1}{\mu_s}\sum_{i}E_{s,i,i^*}(-1)E_{s,i^*,i}=\frac{1}{\mu_s}\nu_sE_s$. 

Therefore, $K_{B,*}=\sum_{s\in P_0} \frac{1}{\mu_s}\nu_s E_s$ and 
$(K_{B,*}, E_{s,i,j})=\left\{\begin{array}{l}
				\nu_s\delta_{i,j}\mbox{ if $s\in P_0$}\\
				0  \mbox{ if $s\notin P_0$}\\
			   \end{array}
				\right. $

By axiom $7^\circ$ $x_s\mu_s=\nu_s$ for $s\in P_0$. Thus, we prove that 
$U=\sum_{i\in P_0}\frac{\nu_i}{\mu_i}e_i+\sum_{j\in P\setminus P_0}x^je_j$,
where $(x^j)^2=\frac{1}{\lambda_j}$. 

Conversely, let $U$ is given by the latter formula. Then reversing
our steps, we obtain that it satisfies axioms $6^\circ$ and $7^\circ$.
By direct calculations we obtain also that $U$ satisfies axioms
$8^\circ$.
\end{proof}

\note In the case of a semisimple structure algebra axiom 
$8^\circ$ follows from axioms $6^\circ$ and $7^\circ$. 
\medskip

Thus, we proved the following theorem.

\begin{theorem} Let 
$\mathcal H =\{H=A\oplus B, (.,.), x\mapsto x^*,  U\}$
be a semisimple structure algebra. Then
\\ (1) $A=\mathbb Ce_1\oplus\dots\oplus\mathbb Ce_m$, 
$B=\M_{n_1}\oplus\dots\oplus \M_{n_k}$, $k\le m$ and
$\phi(e_i)=E_i$ for $i\le k$;
\\ (2) $(X_i, \widetilde X_i)=\mu_i\tr(X_i\widetilde X_i)$, 
$\mu_i\ne 0$, for $X_i,\widetilde X_i\in M_{n_i}$,
$i=1,\dots,k$, and $(X_i, \widetilde X_j) =0$ for $i\ne j$; 
\\ (3) $(e_i,e_i)=\mu_i^2$ for $i=1,\dots, k$; $(e_i,e_i)=\lambda_i$
for $i>k$; $(e_i,e_j)=0$ for $i\ne j$; 
\\ (4) if $M_{n_i}^*=M_{n_j}$ then $e_i^*=e_j$; 
\\ a) if $i\ne j$ then $n_i=n_j$, $\mu_i=\mu_j$ and 
for $X\in M_{n_i}$ we have $X^*=X'\in M_{n_j}$;
\\b) if $i=j$ then either  $X^*=X'$  or $n_i$ is even and   
$X^*=X^\tau$; put $\nu_i=1$ in the first case and $\nu_i=-1$
in the second case;
\\ (5) for $i > k$ we have $e_i^*=e_j$ and if $i\ne j$ then $\lambda_i=
\lambda_j$;
\\ (6) $U=\sum_{i\in P_0}\frac{\nu_i}{\mu_i}e_i+
\sum_{j\in P\setminus P_0}x^je_j$,
where $P=\{x\in\{1,\dots,m\} | e_i^*=e_i\}$, $P_0=P\cap\{1,\dots,k\}$ and 
$(x^j)^2=\frac{1}{\lambda_j}$. 

Conversely, if conditions (1)-(6) hold then $\mathcal H$ is a semisimple 
structure algebra

\end{theorem}

\begin{corollary} \label{cor2.1} If a semisimple associative algebra
$H$, its decomposition $H=A\oplus B$ and bilinear form  $(.,.)$ satisfy
axioms $1^\circ-4^\circ$ then there exists at least one involutive
antiautomorphism satisfying axiom $5^\circ$. For each of them 
there exists $2^{p}>0$, where $p=|P\setminus P_0|$ elements $U\in A$ such
that $\mathcal H=\{H=A\oplus B, (.,.), x\mapsto x^*, U\}$ is a semisimple
structure algebra.
\end{corollary}

Proof is evident.

%%%%%%%%%%%%%%%%%%%%%%%%%%%%%%%%%%%%%%%%%%%%%

\subsection{Structure algebra of a finite group}
\label{2s3}

In this section we assign a structure algebra $\mathcal{H}=\mathcal{H}(G)$
to any finite group $G$. Denote by $\vert \Lambda\vert$
the cardinality of a set $\Lambda$.
Denote by $\mathbb C [G]$ the group algebra of  group $G$.
Denote by $A$ the center of algebra $\mathbb C [G]$.
Let us  assign an element $E_\alpha=\sum_{g\in\alpha}g\in \mathbb C[G]$
to a class $\alpha$ of conjugated elements. It is known, that 
elements  $E_\alpha$ form a basis of $A$.

Denote by $M(G)$ an algebra of all endomorphisms of vector 
space $\mathbb C [G]$. Elements of $G$ form a basis of $\mathbb C [G])$.
For an ordered pair $(g_1,g_2)$ of elements of $G$ 
denote by $E_{g_1,g_2}$ the matrix $((\delta_{g_1,g_2}))$.

Two ordered pairs $(g_1,g_2)$, $(g_1',g_2')$ of elements of $G$ are called 
conjugated if and only if there exists an element $g\in G$ such that 
$gg_1g^{-1}=g_1'$, $gg_2g^{-1}=g_2'$. 
Let $(s_1,s_2)$ be an ordered pair of involutive elements of $G$ and
$\beta$ be a class of ordered pairs conjugated to $(s_1,s_2)$.
Define an element $E_\beta\in M(G)$ by the formula
$E_\beta=\sum_{(s',s'')\in\beta}E_{s',s''}$. 
Denote by $B$ a linear subspace of $M(G)$  generated by all  $E_\beta$.
Obviously elements $E_\beta$ form a basis of $B$.
Action of $g\in G$ on group $G$ by conjugation $x\mapsto gxg^{-1}$ defines
an element of $M(G)$; we denote it by $V_g$.
Obviously mapping $g\mapsto V_g$ defines the
representation $V:\mathbb C[G]\to M(G)$.

\begin{lemma}. \label{V(g)}
\\ (1) $V(g)=\sum_{h\in G}E_{gh,hg}$.
\\ (2) $V(g)E_{g_1,g_2}=E_{gg_1g^{-1},g_2}$, 
      $E_{g_1,g_2}V(g)=E_{g_1,g^{-1}g_2g}$
\\ (3) $V(g)E_{g_1,g_2}V(g^{-1})=E_{gg_1g^{-1},gg_2g^{-1}}$
\end{lemma}
 
The proof is elementary.

\begin{lemma} 
\label{CharB}
Linear subspace $B$ coincides  with the set of 
matrices $X\in M(G)$ such that 

a) $X$ is a linear combination of elements
$E_{s_1,s_2}$, where $s_1,s_2$ are involutive elements of $G$; 

b) $X$ commutes with all $V_g$, $g\in G$.
\end{lemma}

The proof follows from lemma \ref{V(g)}. $\square$

Put $V_\alpha=V(E_\alpha)$.
Define the algebra structure on the direct sum  $H=A\oplus B$.
The multiplication 
of elements of $A$ and elements of $B$ follows from algebra structure on 
$A$ and $B$ respectively.
Define multiplication of elements $E_\alpha\in A$ and 
$E_\beta\in B$ by formulas: $E_\alpha E_\beta=V_\alpha E_\beta$, 
$E_\beta E_\alpha= E_\beta V_\alpha$.

\begin{lemma} $H=A\oplus B$ is a semisimple associative algebra,
$A$ is a central subalgebra, $B$ is an ideal. 
\end{lemma}

\begin{proof} Associativity of $H$ and properties of the decomposition
$H=A\oplus B$ can be easily checked by  direct calculations.
$A$ is a semisimple algebra as the center of a group algebra.

Let us proove that  $B$ is 
a semisimple algebra. Denote by 
$S$ the set of all involutions in $G$. Clearly,
the subspace ${\mathbb C}[S]\subset {\mathbb C}[G]$ is
invariant with respect to the representation $V:\mathbb C[G]\to M(G)$.
By lemma \ref{CharB}, $B$ coincides with the   
centralizer of $V(G)\subset M(S)$.
Therefore, $B$ is a semisimple algebra. Hence $H$ is a 
semisimple algebra.
\end{proof}

There are natural involutive antiautomorphism in both algebras $A$ and $B$.
Namely, denote by "$^*$" the linear extension of the involution 
$g\mapsto g^{-1}$, $g\in G$, to group algebra $\mathbb C[G]$. 
Obviously, $E_\alpha^*=E_{\alpha^*}$, where $\alpha^*$  
consists of inverse elements to elements of $\alpha$ and "$^*$" is 
an involutive antiautomorphism of $A$.

Denote by the same sign "$^*$" the transpose of matrices of $M(G)$, it is  
an involutive antiautomorphism of the algebra.
For a class $\beta=[(s_1,s_2)]$ put $\beta^*=[(s_2,s_1)]$.
Clearly, $E_\beta^*=E_{\beta^*}$. Therefore, $B^*=B$. So we have
defined the involutive transformation  of $H$.

\begin{lemma} Involutive transformation  $x\mapsto x^*$ is an antiautomorphism 
of algebra $H$.
\end{lemma}

\begin{proof} We should check that
$(E_\alpha E_\beta)^*=E_\beta^* E_\alpha^*$.
We have $(E_\alpha E_\beta)^*=(V_\alpha E_\beta)^*=E_\beta^* V_\alpha^*$. 
By lemma \ref{V(g)} $V_g=\sum_{h\in G}E_{gh,hg}$. Therefore,
\\$V_g^*=\sum_{h\in G}E_{hg,gh}=\sum_{h\in G}E_{g^{-1}h',h'g^{-1}}$, where 
$h'=ghg$. Hence $V_g^*=V_{g^{-1}}$ and $V_\alpha^*=V_{\alpha^*}$.
\end{proof}

It is well-known that any invariant symmetric bilinear form $(.,.)$ 
is uniquely defined by a linear form $f(x)$ such that $f(xy-yx)=0$ identically.
Bilinear form corresponding to $f(x)$ is defined as  $(x,y)=f(xy)$.
Define a linear form $f$ on $H$ by formulas: $f(E_\alpha)=
\frac{1}{\vert G\vert}\delta_{\alpha,1}$,
$f(E_\beta)=\frac{1}{\vert G\vert}\tr(E_\beta)$ (recall that $E_\beta$ 
is an element of matrix algebra $M(G)$).
Clearly, $f$ defines the invariant symmetric bilinear form
$(.,.)$ on $H$. 

The restriction $(.,.)|_A$ coincides with the restriction to the center
of the standard invariant form on $\mathbb C [G]$ given by formula
$(g_1,g_2)=\frac{1}{\vert G\vert}\delta_{g_1,g_2^{-1}}$ for $g_1,g_2\in G$.
Therefore, $(E_{\alpha_1}, E_{\alpha_2})=\frac{1}{\nu_{\alpha_1}}
\delta_{\alpha_1,\alpha_2^*}$, where
$\nu_{\alpha}=\frac{\vert G\vert}{\vert\alpha\vert}$.
Obviously $\nu_\alpha$ is equal to
the number of elements in the centralizer of any $g\in\alpha$.
Note that form  $(.,.)|_A$ is nondegenerate.

The restriction $(.,.)|_B$ coincides with the restriction to $B$ of the
standard invariant form on $M(G)$ given by formula 
$(X_1,X_2)=\frac{1}{\vert G\vert}\tr(X_1X_2)$. Therefore,
$(E_{\beta_1},E_{\beta_2})=\frac{1}{\nu_{\beta_1}} 
\delta_{\beta_1,\beta_2^*}$, where
$\nu_{\beta}=\frac{\vert G\vert}{\vert\beta\vert}$. Obviously 
$\nu_\beta$ is equal to the number of elements in the stabilizer 
of any element  $(s_1,s_2)\in\beta$. 
Clearly, the form  $(.,.)|_B$ is nondegenerate.
Note that scalar products $(E_\alpha,E_\beta)$
can be nonzero because by definition
$(E_\alpha,E_\beta)=f(E_\alpha E_\beta)=
\frac{1}{\vert G\vert}\tr(V_\alpha E_\beta)$.

Clearly, involution $x\mapsto x^*$ preserves this bilinear form.
According to notations of previous subsection put
$F_{\alpha_1,\alpha_2}=(E_{\alpha_1}, E_{\alpha_2})$,
$F_{\beta_1,\beta_2}=(E_{\beta_1},E_{\beta_2})$ and
$R_{\alpha,\beta}= (E_\alpha,E_\beta)$.

\begin{lemma} Cardy axiom $4^\circ$ holds for algebra $H$.
\end{lemma}

\begin{proof} First, let us compute  the left hand side
$L=(V_{K_B}(E_{\beta_1}),E_{\beta_2})$ of the identity in axiom $4^\circ$. 
By definition,
$V_{K_B}(E_{\beta_1})=F^{\beta',\beta''}E_{\beta'}E_{\beta_1}E_{\beta''}$.
We have $F^{\beta',\beta''}=\nu_{\beta'}\delta^{\beta',\beta''{^*}}$.
Hence $V_{K_B}(E_{\beta_1})=
\sum_{\beta}\nu_{\beta}E_{\beta}E_{\beta_1}E_{\beta^*}$ and
\\$L=\sum_{\beta}\nu_\beta\frac{1}{\vert G\vert}
\tr (E_{\beta}E_{\beta_1}E_{\beta^*}E_{\beta_2})=
\sum_{\beta}\frac{1}{\vert\beta\vert}\tr(E_{\beta}
E_{\beta_1}E_{\beta^*}E_{\beta_2})$.

Each summand in this formula is a product of matrices
$E_{x,y}E_{s_1,s_2}E_{u,v}E_{s_3,s_4}$, where 
$(x,y)\in \beta, (s_1,s_2)\in\beta_1,(u,v)\in \beta^*,(s_3,s_4)\in\beta_2$.
If the trace of the summand is nonzero then $y=s_1,u=s_2,v=s_3,x=s_4$.
Let this conditions be satisfied. Then the trace is equal to $1$ and
the pair $(s_4,s_1)$ is conjugated to the pair $(s_3,s_2)$. Therefore,
$$L=\sum_{(s_1,s_2,s_3,s_4,g|(s_1,s_2)\in\beta_1,
(s_3,s_4)\in\beta_2, gs_2g^{-1}=s_1,gs_3g^{-1}=s_4,g\in G}
\frac{1}{\vert[(s_4,s_1)]\vert}\frac{1}{\vert Z_G(s_4,s_1)\vert}=$$
$$=\frac{\vert\{(s_1,s_2,s_3,s_4,g|(s_1,s_2)\in\beta_1,
(s_3,s_4)\in\beta_2,g\in G, gs_2g^{-1}=s_1,gs_3g^{-1}=s_4\}\vert}
{\vert G\vert}$$

Let us compute the right hand side 
$R=(\widehat{K}_A,E_{\beta_1}\otimes E_{\beta_2})$  of the identity from 
axiom $4^\circ$. 
By definition,

$R=\sum_{\alpha',\alpha''} F^{\alpha',\alpha''}(E_{\alpha'},E_{\beta_1})
     (E_{\alpha''},E_{\beta_2})=\sum_{\alpha}{\nu_\alpha}
     (E_{\alpha},E_{\beta_1})\\ \cdot (E_{\alpha^*},E_{\beta_2})=
\sum_{\alpha}\nu_{\alpha} \frac{1}{\vert G\vert^2}
		\tr(V_{\alpha}E_{\beta_1})
                \tr(V_{\alpha^*}E_{\beta_2})$.

By lemma \ref{V(g)}, we have
$V_{\alpha}=\sum_{h\in G, g\in\alpha}E_{gh,hg}=
\sum_{x\in G, g\in\alpha}E_{x,g^{-1}xg}$. Analogously,
$V_{\alpha^*}=\sum_{y\in G, h\in\alpha}E_{y,hyh^{-1}}$.

Hence,

$R=\sum_{\alpha}\frac{1}{\vert G\vert\vert\alpha\vert}
     (\sum_{x\in G, g\in\alpha, (s_1,s_2)\in\beta_1}
    \tr(E_{x,g^{-1}xg}E_{s_1,s_2}))\cdot$
\\ $\cdot(\sum_{y\in G, h\in\alpha,(s_3,s_4)\in\beta_2}
                        \tr(E_{y,hyh^{-1}}E_{s_3,s_4}))$.

Nonzero summands arise from the cases 
$x=s_2, 
g^{-1}xg=s_1,
y=s_4, 
hyh^{-1}=s_3$; in this case a summand is equal to 
$\frac{1}{\alpha}$.

Therefore,

$R=\sum_{\alpha}\frac{1}{\vert G\vert\vert\alpha\vert}
\vert\{s_1,s_2,s_3,s_4,g,h |
(s_1,s_2)\in\beta_1, (s_3,s_4)\in\beta_2,g\in\alpha,h\in\alpha,
s_2=gs_1g^{-1},s_3=hs_4h^{-1}\}\vert$.

Element $h$ runs through the set $\alpha$ 
of elements conjugated to $g$. Hence, we can present it as 
$h=zgz^{-1}$, where $z$ runs through all elements of group $G$. Each 
element $h$ has $k$ presentations of this type, where 
$k=\vert Z_G(g)\vert$ is the cardinality of the centralizer of element $g$.
Therefore,
$R=\sum_{\alpha}\frac{1}{\vert G\vert\vert\alpha\vert}
		\frac{1}{\nu_\alpha} 
        \vert\{s_1,s_2,s_3,s_4,g,z |
		(s_1,s_2)\in\beta_1, (s_3,s_4)\in\beta_2,g\in\alpha, z\in G,
                s_2=gs_1g^{-1},s_3=(zgz^{-1})s_4(zg^{-1}z^{-1}\}\vert.$
Note that $\frac{1}{\vert G\vert\vert\alpha\vert}\frac{1}{\nu_\alpha}=
\frac{1}{\vert G\vert^2}$ and therefore,
this coefficient can be carried out of summing. Therefore,
$R=\frac{1}{\vert G\vert^2}
        \vert\{s_1,s_2,s_3,s_4,g,z |
		(s_1,s_2)\in\beta_1, (s_3,s_4)\in\beta_2,g\in G, z\in G,
                s_2=gs_1g^{-1},s_3=(zgz^{-1})s_4(zg^{-1}z^{-1}\}\vert.$

The equality $s_3=(zgz^{-1})s_4(zg^{-1}z^{-1})$ can be rewritten as 
$s_3'=gs_4'g^{-1}$, where $s_3'=z^{-1}s_3z$, $s_4'=z^{-1}s_4z$.  
For two fixed $z=z_1$ and $z=z_2$ the numbers of tuples 
$(s_1,s_2,s_3,s_4,g,z)$ satisfying conditions are equal since the
pairs $(s_3,s_4)$ run through all representatives of the class $\beta_2$.
Therefore,
\\$R=\frac{1}{\vert G\vert^2}\vert G\vert
        \vert\{s_1,s_2,s_3,s_4,g,z |
		(s_1,s_2)\in\beta_1, (s_3,s_4)\in\beta_2,g\in G,
                s_2=gs_1g^{-1},s_3=gs_4 g^{-1}\}\vert.$

We get that $L=R$.
\end{proof}

Thus, we have proved all axioms for a structure algebra except those
concerning an element $U\in A$. By corollary \ref{cor2.1},
there exists at least one
element $U\in A$ satisfying axioms $6^\circ-8^\circ$. Let us fix it.
Therefore, we prove theorem.

\begin{theorem} For any finite group $G$ the set of data 
$\mathcal H(G)=(H=A\oplus B, (.,.), x\mapsto x^*, U)$ is a 
structure algebra. 
\end{theorem}

\note It can be shown that in the case of the symmetric group $G=S_n$
the element $U$ is  uniquely determined.

%%%%%%%%%%%%%%%%%%%%%%%%%%%%%%%%%%%%%%

\section{Cuts of stratified surfaces}
\label{sec3}

This section contains a geometrical background for the Klien topological field theory developed in section \ref{sec4}. 
We consider  surfaces with  boundaries and
special points. The simplest classes of such surfaces (trivial, basic and simple
surfaces in the nomenclature of subsections 3.1 and 3.3 ) we exploit in section \ref{sec4} for constructing of sructure algebras.

We consider also (subsections 3.2 and 3.3) systems of nonintersecting generic cuts of the surfaces, and in particular, complecte cut systems cuting surfaces into basic surfaces. The simple surfaces generate elementary shifts of complecte cut systems. 

The central  theorem of the section (theorem \ref{th4.1},subsection 3.4) claims that any complete cut system of
a surface can be transformed to any other complete cut system of the same surface by 
 elementary shifts. 
 This theorem playes an essential role in the proof of the main theorem of the paper (theorem \ref{MMain}) about the equivalence between structure algebras and Klein topological field theories.

\subsection{Stratified surfaces}

Denote by $\partial X$ a boundary of a topological space $X$
and denote  by  $X^\circ$  its interior $\overline{X}\setminus\partial X$.
We deal below with compact topological manifolds possibly with boundary and
call them 'manifolds' for short. 

{\bf Examples} 
\begin{enumerate}
\item A connected one-dimensional manifold 
is homeomorphic either to a circle or to a segment.
\item A connected orientable two-dimensional manifold 
is homeomorphic to a sphere with 
$g$ handles and $s$ holes. We call it {\it a
surface of type} $(g,s,1)$.                        
\item A connected nonorientable two-dimensional manifold 
is homeomorphic either to a projective plane with 
$a$ handles and $s$ holes or to a Klein bottle 
with $a$ handles and $s$ holes. We call it 
{\it a surface of  type} $(g,s,0)$, where 
$g=a+1/2$ in the former case and 
$g=a+1$ in the latter case. 
\end{enumerate}
 
\begin{definition} Let $\Lambda$ be a finite set. 
The decomposition
$\Omega=\coprod_{\lambda\in\Lambda}\Omega_\lambda$ of a topological space 
$\Omega$ into the disjoint union of subspaces (strata)
$\Omega_\lambda\subset\Omega$ is called a stratification if and only if
\begin{enumerate}
\item each stratum $\Omega_\lambda$ is homeomorphic to a connected manifold
without a boundary;
\item each stratum is open in its closure;
\item the boundary $\partial\Omega_\lambda$ of a stratum coincides with
the union of strata of less dimension.
\end{enumerate}

We call topological space  $\Omega$ itself a base of the stratification.
Zero-dimensional strata are called {\it special points}.

\end{definition}

\note
\begin{enumerate}
\item A stratum is not necessarily homeomorphic to an open ball.
\item Joining of strata induces a partial ordering of  set
$\Lambda$.
\item There are several different definitions of
a stratified manifold (see \cite{Hug-Wei}). All these  definitions are
equivalent in dimensions $1$ and $2$. Fortunately, in this work we deal 
with dimensions $1$ and $2$ only. $\ \ \square$
\end{enumerate}

A homeomorphism $\phi:\Omega\to\Omega'$ of the bases of stratified manifolds
$\Omega=\coprod_{\lambda\in\Lambda}\Omega_\lambda$ and 
$\Omega'=\coprod_{\lambda\in\Lambda'}\Omega_\lambda'$
is called {\it an isomorphism of stratifications},
if and only if the restriction of $\phi$ to any stratum of $\Omega$ 
is the homeomorphism with an appropriate stratum of 
$\Omega'$. Clearly, an isomorphism of stratifications induces 
the isomorphism of partial ordered sets 
$\Lambda$ and $\Lambda'$.

\begin{definition} A stratification 
$\Omega=\coprod_{\lambda\in\Lambda}\Omega_\lambda$
of a manifold $\Omega$ is called  {\it special stratification} 
if and only if all strata of codimension one belong to the boundary 
of $\Omega$. 
\end{definition}

{\bf Examples} 
\begin{enumerate}
\item
 A connected specially stratified one-dimensional 
manifold is isomorphic either to a circle with the unique stratum or to 
a segment with natural stratification (strata are two end points and the
open interval). 
\item
Let $\Omega=\coprod_{\lambda\in\Lambda}\Omega_\lambda$ be a connected 
specially stratified two-dimensional manifold. 
We call the set of data
$G=(g, \varepsilon,m,m_1,\dots,m_s)$ {\it a type of} $\Omega$. Here
$(g,s,\varepsilon)$ is the type of surface $\Omega$, 
$m$ is the number of interior special points,  
$m_i$ is the number of special points on $i$-th boundary contour.
Two types that differ only in the order of $m_i$ are considered as equal.
It can be easily shown that an isomorphism class connected two-dimensional
specially stratified manifolds is uniquely determined by its type.
\item 
Let $\Omega$ be a stratified surface consisting of connected components 
$\Omega_i$. Denote by $G$ unordered set of types $G^i$ of surfaces $\Omega_i$.
We call $G$ a  type of $\Omega$. 
Clearly, up to isomorphism, $\Omega$ is uniquely determined by its type.
\end{enumerate}

\note A special stratification of a surface is uniquely defined by a set 
$\Omega_0$ of its special points.  

We call a two-dimensional specially stratified manifold $\Omega$ consisting of
finitely many connected components {\it a stratified surface} for short.
Let $\Omega$ consists of $c$ connected components.
Denote by
$G^i=(g^i,\varepsilon^i,m^i,m_1^i,\dots,m_s^i)$ ($i=1,\dots,c$) a type of $i$-th
connected component  $\Omega^i$. Define an invariant $\mu(\Omega)$ by formula

$$
\mu(\Omega)=
\sum_i 2g^i + \sum_i m^i + \sum_i s^i +\frac{1}{2}\sum_i\sum_j m_j^i - 2c
$$

From the definition follows that $\mu(\Omega)$ is half-integer.
Clearly, invariant $\mu$ is additive with respect to the  disjoint union of
stratified surfaces:                                        
$$\mu(\Omega_1\coprod\Omega_2)=\mu(\Omega_1)+\mu(\Omega_2)$$

We say that a connected stratified surface $\Omega$ is {\it a trivial surface} 
if and only if $\mu(\Omega)\le 0$. 
All trivial surfaces can be easily listed.

\begin{lemma} \label{3.1} Any trivial stratified surface is isomorphic to 
one of surfaces from the following list:
\\ (1)  sphere $S^2$ without special points ($\mu=-2$);
\\ (2)  projective plane $\mathbb{R}P^2$ without special points ($\mu=-1$); 
\\ (3)  disc  $D^2$ without special points ($\mu=-1$); 
\\ (4)  sphere $(S^2,p)$ with a unique interior special point  $p$ ($\mu=-1$);
\\ (5)  disc $(D^2,q)$ with a unique boundary special point $q$ and without interior 
special points ($\mu=-\frac{1}{2}$);
\\ (6)  sphere $(S^2,p_1,p_2)$ with two interior special points ($\mu=0$);
\\ (7)  projective plane $(\mathbb{R}P^2,p)$ with one interior special point
($\mu=0$);
\\ (8)  torus $T^2$ without special points ($\mu=0$); 
\\ (9)  Klein bottle $\Kl$ without special points  ($\mu=0$); 
\\ (10) disc  $(D^2,p)$ with a unique interior special point and 
      without boundary special points ($\mu=0$);
\\ (11) disc $(D^2,q_1,q_2)$ with two boundary special points and 
      without interior special points ($\mu=0$);
\\ (12) M\"obius band $\Mb$ without special points  ($\mu=0$); 
\\ (13) cylinder $\Cyl$ without special points ($\mu=0$). 
\end{lemma}

%%%%%%%%%%%%%%%%%%%%%%%%%%%%%%%%%%

\subsection{Cut systems}
\label{cutsys}
Let $\Omega$ be a stratified surface.
A generic not self-intersecting curve $\gamma\subset \Omega$
is called  {\it a simple cut}. This implies that
$\gamma$  does not meet any special point and either is (closed) contour 
consisting purely of interior points, or
is a segment, its end  points belong to the boundary of $\Omega$ and 
all interior points of the segment  are interior points of the surface.

A set $\Gamma$ of pairwise nonintersecting simple cuts $\gamma\subset\Omega$
is called  {\it a cut system of $\Omega$} (Fig.1a).

Let $\Omega$, $\Omega'$ be two stratified surfaces
endowed with cut systems $\Gamma$, $\Gamma'$. We say that
an isomorphism $\phi:\Omega\to\Omega'$ of stratified surfaces is
{\it an isomorphism of pairs}
$(\Omega,\Gamma)$, $(\Omega',\Gamma')$ if and only if 
$\phi(\Gamma)=\Gamma'$. 

\begin{definition} A triple  $(\Omega_*,\Gamma_*,\tau)$
consisting of 
\begin{itemize}
\item[$\circ$] a stratified surface $\Omega_*$;

\item[$\circ$] a subset $\Gamma_*\subset\partial\Omega_*$ such that 
each connected component of $\Gamma_*$ coincides with the closure 
of a one-dimensional stratum of $\Omega_*$;

\item[$\circ$] an involutive homeomorphism  $\tau:\Gamma_*\to\Gamma_*$ 
having no fixed points 
\end{itemize}

is called cut surface. In this case pair $(\Gamma_*,\tau)$  
is called a gluing system.

An isomorphism of cut surfaces $(\Omega_*,\Gamma_*,\tau)$,
$(\Omega_*',\Gamma_*',\tau')$ is an isomorphism 
$\phi:\Omega_*\to\Omega_*'$ of stratified surfaces such that 
$\phi(\Gamma_*)=\Gamma_*'$ and  $\tau'\circ\phi=\phi\circ\tau$.
\end{definition}

Let $(\Omega_*,\Gamma_*,\tau)$   be a cut surface. Gluing points
$x$ and $\tau(x)$ we obtain a surface $\Omega$ and 'gluing topological map'  
$\glue:\Omega_*\to\Omega$. 
Clearly, stratification of $\Omega_*$ induces the stratification 
of $\Omega$ and the image $\Gamma=\glue(\Gamma_*)$ 
is a cut system of $\Omega$ and 
$\glue:(\Omega_*,\Gamma_*,\tau)\to (\Omega,\Gamma)$ 
is a functor from the category of cut surfaces  to
the category of pairs $(\Omega,\Gamma)$ (in both categories 
morphisms are isomorphisms). 

Conversely, if $\Gamma$ is a cut system of $\Omega$ then one can
construct cut  surface $(\Omega_*,\Gamma_*,\tau)$ as follows. Points of 
$\Omega_*$ are points of $\Omega\setminus \Gamma$ and pairs 
$(x,c)$, where $x\in\Gamma$ and $c$ is a coorientation of $\Gamma$ in a 
neighborhood of point $x$. The stratification of $\Omega_*$ and gluing
system $(\Gamma_*,\tau)$ are defined evidently (Fig.1b). Clearly, we
constructed a functor $\cut:(\Omega,\Gamma)\to (\Omega_*,\Gamma_*,\tau)$.

\begin{lemma} Functors $\glue$ and $\cut$ establish
the equivalence of categories of cut surfaces 
(with morphisms defined as isomorphisms of 
cut surfaces) and  
pairs  $(\Omega,\Gamma)$ (with morphisms defined as isomorphisms of 
pairs).
\end{lemma}

\medskip
\epsfxsize=12truecm
\begin{figure}[tbhp]
\centerline{\epsffile{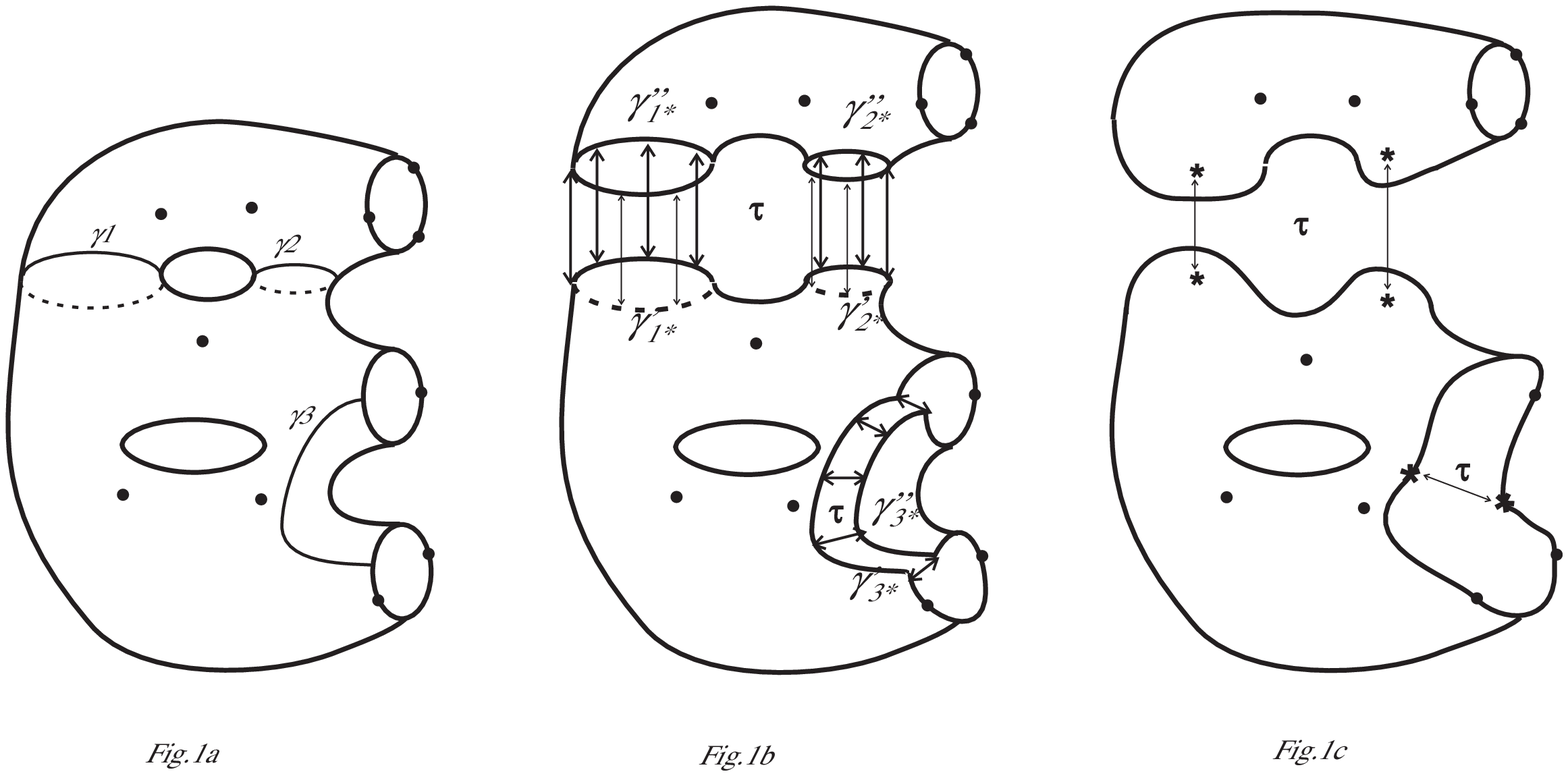}}
\caption{The construction of a contracted cut surface.\newline
Fig.1a. Cut system $\Gamma=\{\gamma_1,\gamma_2,\gamma_3\}$ of surface 
$\Omega$.\newline  
Fig.1b. Cut surface $(\Omega_*,\Gamma_*,\tau)$ obtained from cut system 
$\Gamma$. Here $\Gamma_*=\{\gamma_{1*}',\gamma_{1*}'',\gamma_{2*}',
\gamma_{2*}'', \gamma_{3*}',\gamma_{3*}''\}$.\newline
Fig.1c. Contracted cut surface $\Omega_\#=\Omega/\Gamma$. 
Special points obtained as contractions of the connected components 
of $\Gamma_*$ are marked by $*$.}
\end{figure}
\medskip

Proof is skipped.

Let $(\Omega_*,\Gamma_*,\tau)$ be a cut surface. Denote 
by $\Omega_{\#}$ a surface obtained by contracting each connected 
component of $\Gamma_*$ to a point. Clearly, the stratification
of $\Omega_*$ induces the stratification of  $\Omega_\#$. Special
points of this stratification are images of special points of 
$\Omega_*$ and points that are obtained as contracted connected 
components of $\Gamma_*$ (Fig.1c). Homeomorphism $\tau$ induces
the involution of the set of special points of $\Omega_\#$
coming from contracted connected components.
We denote this involution by the same letter $\tau$. Clearly,
we have constructed a functor $\contr:(\Omega_*,\Gamma_*,\tau)\to \Omega_\#$
from category of cut surfaces  to a category of stratified surfaces.

Therefore, we can assign
a stratified surface $\Omega_{\#}=\contr(\cut(\Omega,\Gamma))$ to any pair 
$(\Omega,\Gamma)$, where $\Gamma$ is a cut system of a stratified 
surface $\Omega$. We call $\Omega_{\#}$ {\it a contracted cut surface}
and denote it by $\Omega/\Gamma$.

\note The set $(\Omega_\#)_0$ of special 
points of a contracted cut surface $\Omega_\#=\Omega/\Gamma$ carries 
an additional structure. Namely, 1) $(\Omega_\#)_0$ contains the 
distingushed subset $\Omega_0$ consisting of special points
that are images of special points of  $\Omega$; 2)
there is fixed involution $\tau: (\Omega_\#)_0\setminus\Omega_0\to 
(\Omega_\#)_0\setminus\Omega_0$ ($\tau$ may have fixed points).

Let us classify simple cuts of a surface as follows. 
Denote by $\gamma$ an arbitrary simple cut of a 
connected stratified surface $\Omega$ of type 
$G=(g,\varepsilon,m,m_1,\dots,m_s)$. Two situations 
may occur. First, contracted cut surface $\Omega/\{\gamma\}$ is 
a connected surface. Denote its type by  
$G_\#=(g_\#,\varepsilon_\#,m_\#,m_{\#1},\dots,m_{\#s})$.
Second, $\Omega/\{\gamma\}$ consists of two connected components.
Then type $G_\#$ of $\Omega/\{\gamma\}$ is the set of types 
of its connected components. 
Denote types of components  by  
\\$G_\#^i=(g_\#^i,\varepsilon_\#^i,m_\#^i,m_{\#1}^i,\dots,m_{\#s^i}^i)$,
where $i=1,2$. In this case the image of a special 
point of $\Omega$ in $\Omega_\#$ belongs to one of components. 
Thus, $\gamma$ induces a division $\Omega_0=(\Omega_0)_1\sqcup(\Omega_0)_2$
of set $\Omega_0$  of special points of $\Omega$ into two subsets.

For a cut system ${\gamma}$  consisting of one simple cut we use notations
$(\Omega_*,\gamma_*,\tau)$ for cut surface obtained by cutting stratified 
surface $\Omega$ along ${\gamma}$.

\begin{lemma} Let $\Omega$ be a stratified surface.
Then
\begin{enumerate}
\item any simple cut  $\gamma$ of $\Omega$ belongs to just one 
of classes $1-9$ from the list below;
\item the identities between types  $G_{\#}$ and $G$ that 
are written in each item of the list holds.
\end{enumerate}
\end{lemma}
\medskip
\begin{center}
{\it List of classes of simple cuts}
\end{center}
\smallskip
A. {\it In classes 1-4 a simple cut  $\gamma$ is supposed to be
homeomorphic to a circle}

\smallskip
{\bf Class 1} consists of {\it separating contours} (Fig.2a).  
A simple cut $\gamma$ is called a separating contour if and only if 
$\gamma_*\subset\Omega_*$ consists of two contours and
surface  $\Omega_\#$ consists of two connected components,

$G=(g_\#^1+g_\#^2,\varepsilon_\#^1\varepsilon_\#^2,
   m_\#^1+m_\#^2-2, m_{\#1}^1,\dots,m_{\#s_1}^1,m_{\#1}^2,
   \dots,m_{\#s_2}^2)$

\smallskip
{\bf Class 2} consists of {\it cuts of a handle} (Fig.2b).   
A simple cut $\gamma$ is called 
a cut of a handle if and only if
$\gamma_*\subset\Omega_*$ consists of two contours,
  $\Omega_\#$ is connected surface and $\varepsilon_\#=\varepsilon$,
                              
$G=(g_\#+1,\varepsilon_\#,m_\#-2, m_{\#1},\dots,m_{\#s})$.

\smallskip
{\bf Class 3} consists of {\it cuts of a neck of Klein bottle} (Fig.2c).  
A simple cut $\gamma$ is called a cut of a neck of Klein bottle
if and only if  $\gamma_*\subset\Omega_*$ consists of two contours,
 $\Omega_\#$ is connected surface , $\varepsilon=0$
and $\varepsilon_\#=1$,

$G=(g_\#+1,0,m_\#-2, m_{\#1},\dots,m_{\#s})$.

\smallskip
{\bf Class 4} consists of  {\it M\"obius cuts} (Fig.2d).  
A simple cut $\gamma$ is called a M\"obius cut if and only if
$\gamma_*\subset\Omega_*$ is connected contour.
In this case $\varepsilon=0$,
$\Omega_\#$ is connected surface ,  $\varepsilon_\#$ is equal either to
$0$ or to $1$,

$G=(g_\#+\frac{1}{2},0,m_\#-1, m_{\#1},\dots,m_{\#s})$.                          
\bigskip          

B. {\it In classes 5-9 a simple cut $\gamma$  is supposed to be
homeomorphic to a segment.}
Hence,  $\gamma_*$ consists of two disjoint segments.

\smallskip
{\bf Class 5} consists of {\it cuts between two holes} (Fig.2e).  
A simple cut $\gamma$ is called a  cut between two holes if and only if 
end points of $\gamma$ belong to different boundary contours. 
In this case $\Omega_\#$ is a connected surface and 
in an appropriate numeration of boundary contours we obtain the identity
$G_\#=(g,\varepsilon,m,m_1+m_2+2,m_3,\dots,m_s)$.
\medskip          

C. {\it In classes 6-9 both boundary points of a segment $\gamma$ 
are supposed to belong to the same boundary contour $\omega$.}
Denote by  $\omega_*$ the preimage of $\omega$ in cut surface  $\Omega_*$ 
and by $\omega_\#$ the image of $\omega_*$ in contracted cut 
surface $\Omega_\#$.

\smallskip
{\bf Class 6} consists of {\it separating segments} (Fig.2f).  
A simple cut $\gamma$ is called a  separating segment if and only if 
$\Omega_\#$ consists of two connected components.
In this case $\omega_\#$ consists of two contours and 
in an appropriate numeration of boundary contours of $\Omega_{\#}$ 
we obtain the identity

$G=(g_\#^1+g_\#^2,\varepsilon_\#^1\varepsilon_\#^2,
	m_\#^1,+m_\#^2, m_{\#1}^1+m_{\#1}^2-2,m_{\#2}^1,\dots,m_{\#s_1}^1,
        m_{\#2}^2,\dots,m_{\#s_2}^2)$.

\smallskip
{\bf Class 7} consists of  {\it cuts of a handle through a hole} (Fig.2g). 
 
A simple cut $\gamma$ is called a cut of a handle through a hole,
if and only if
$\Omega_\#$ is connected surface, 
$\omega_\#$ consists of two contours and
$\varepsilon_\#=\varepsilon$.
In this case in an appropriate numeration of boundary contours of
$\Omega_{\#}$ we obtain the identity

 $G=(g_\#+1,\varepsilon_\#,m_\#, m_{\#1}+m_{\#2}-2,m_{\#3},\dots,m_{\#s})$.

\smallskip
{\bf Class 8} consists of  {\it cuts of a neck of a Klein bottle through a
hole} (Fig.2h).  
A simple cut $\gamma$ is called a cut of a neck of a Klein bottle through a
hole if and only if $\Omega_\#$ is connected surface, 
$\omega_\#$ consists of two contours,
$\varepsilon=0$ and $\varepsilon_\#=1$.
 
In this case
in an appropriate numeration of boundary contours of $\Omega_{\#}$ 
we obtain the identity 

 $G=(g_\#+1,0,m_\#, m_{\#1}+m_{\#2}-2,m_{\#3},\dots,m_{\#s})$.

\smallskip
{\bf Class 9} consists of  {\it cuts across M\"obius band} (Fig.2i).  
A simple cut $\gamma$ is called a cut across M\"obius band  if and only if 
boundary contour  $\omega_\#$ is connected contour.

In this case $\Omega_\#$ is a connected surface,
$\varepsilon=0$, $\varepsilon_\#$ equals either to $0$ or to $1$  
and 
in an appropriate numeration of boundary contours of $\Omega$ 
we obtain the identity

$G=(g_\#+\frac{1}{2},0,m_\#, m_{\#1}-2,m_{\#2},\dots,m_{\#s})$.

\medskip
\begin{figure}[tbhp]
\centerline{\epsffile{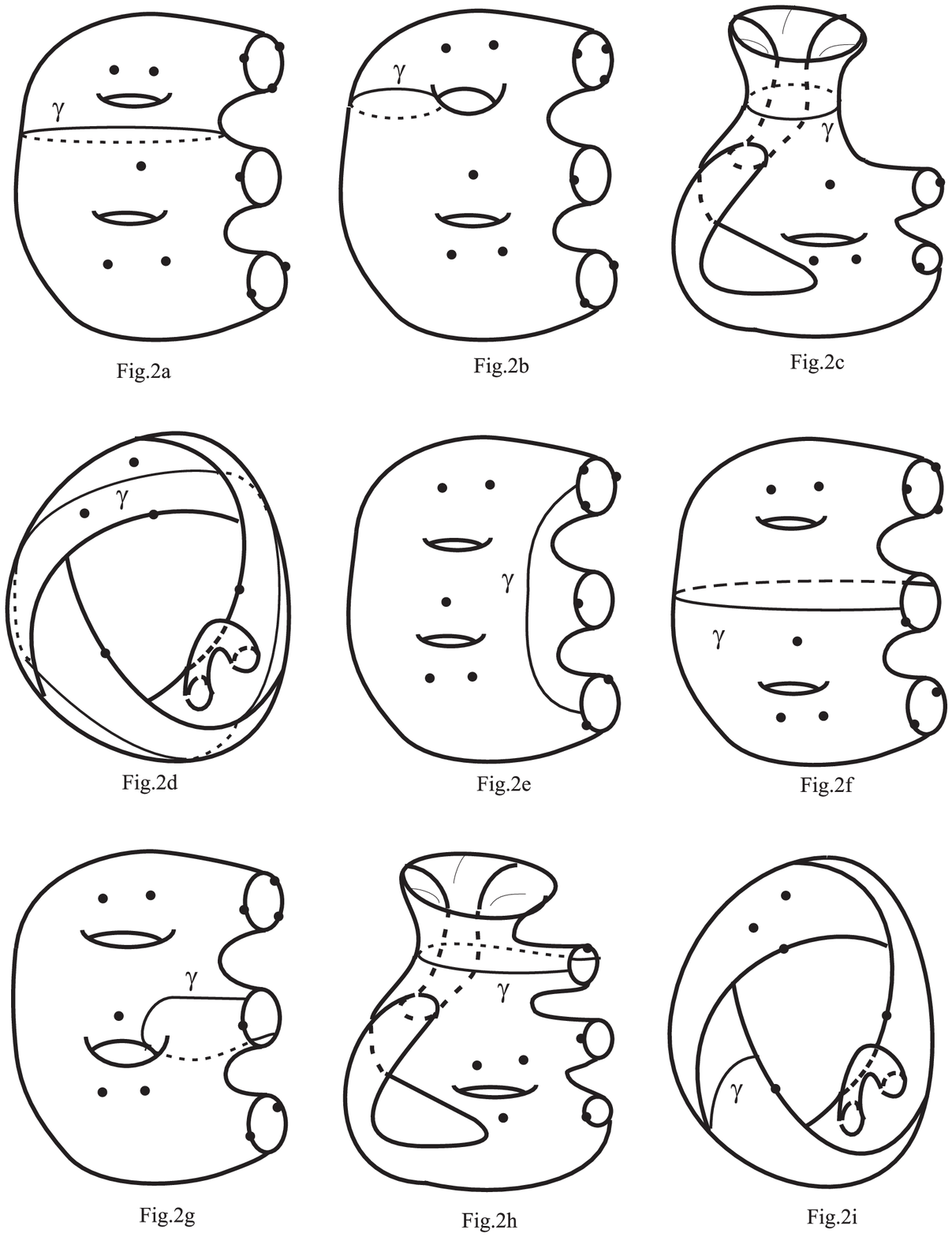}}
\caption{Examples of cuts of classes 1-9.
Each cut is marked by $\gamma$. 
Fig.2a. A separating contour. 
Fig.2b. A cut of a handle.
Fig.2c. A cut of a neck of Klein bottle.
Fig.2d. A M\"obius cut. 
Fig.2e. A cut between two holes. 
Fig.2f. A separating segment.
Fig.2g. A cut of a handle through a hole.
Fig.2h. A cut of a neck of a Klein bottle through a hole.
Fig.2i. A cut across M\"obius band}
\end{figure}
%\medskip
%\bigskip

\begin{lemma} Let $\Gamma$ be a cut system of a stratified surface 
$\Omega$. Then $\mu(\Omega)=\mu(\Omega_\#)$, where 
$\Omega_\#=\Omega/\Gamma$ is the contracted cut surface.
\end{lemma}

\begin{proof} The equality can be easily checked for 
every class of simple cuts.
\end{proof}

\begin{definition} Two cut systems $\Gamma'$ and $\Gamma''$ of $\Omega$ 
are called equivalent if there exists 
an isomorphism $\phi:\Omega\to\Omega$ of the stratified surfaces
conserving all special points, all one-dimensional strata and their 
orientations,  and such that $\phi(\Gamma')=\Gamma''$.
\end{definition}

\note In this  definition the requirement of conserving one-dimensional strata of 
a stratified surface $\Omega$ is essential only for boundary contours with 
zero, one or two special points.
If a boundary contour contains more that two special points then any isomorfism
of $\Omega$ conserving all special  points conservses also all one-dimensional strata
on the contour and their orientations.

\begin{lemma} \label{33} Two simple cuts $\gamma'$ and $\gamma''$
of a stratified surface $\Omega$ are equivalent if and only if
\begin{enumerate}
\item $\gamma'$ and $\gamma''$ belong to the same class;
\item types of contracted cut surfaces $\Omega_\#'=\Omega/\gamma'$ and
$\Omega_\#''=\Omega/\gamma''$ coincide;
\item if $\gamma'$ and $\gamma''$ are separating contours or separating 
segments the divisions of the set of special points of $\Omega$, 
induced by $\gamma'$ and $\gamma''$, coincide;
\item if $\gamma'$ and $\gamma''$ are segments then their end points
belong to the same one-dimensional strata.
\end{enumerate}
\end{lemma}
           
\begin{proof}
Obviously, if $\gamma'$ and $\gamma''$ are equivalent cuts of a stratified surface 
$\Omega$ then conditions (1)--(4) are satisfied. 

Conversely, let $\gamma'$ and $\gamma''$ be simple cuts satisfying conditions (1) -- (4).
Denote by $(\Omega'_*,\gamma'_*,\tau'_*)$ and
$(\Omega''_*,\gamma''_*,\tau''_*)$ cut surfaces obtained by cutting 
$\Omega$ along $\gamma'$ and $\gamma''$ resp. 
Conditions (1) -- (4) provide the equality of types of stratified 
surfaces  $\Omega'_*$ and $\Omega''_*$.
Therefore, there exists an 
isomorphism $\phi:\Omega'_*\to\Omega''_*$  of stratified surfaces.
Moreover, it can be easily shown for every class of simple cuts
that one can choose $\phi$ satisfying the properties: 
\begin{itemize}
\item[$\circ$] for 
any special point $r\in\Omega_0$ of $\Omega$  isomorphism $\phi$ 
brings the image $r'$ of $r$ in $\Omega'_*$ to the image
$r''$ of $r$ in $\Omega''_*$
\item[$\circ$] $\phi(\gamma'_*)=\gamma''_*$
\item[$\circ$] $\phi$ is the isomorphism of cut surfaces, i.e., 
$\phi\circ\tau'=\tau''\circ\phi$. 
\end{itemize}

Clearly, $\phi$ generates an homeomorphism $\overline{\phi}:\Omega\to\Omega$
such that $\overline{\phi}$ preserves all special points and 
$\overline{\phi}(\gamma')=\gamma''$.
\end{proof}

\begin{corollary} \label{34a}
Let  $\Gamma'$ and $\Gamma''$ be two cut systems of
a stratified surface $\Omega$. Suppose $\Gamma'$ contains a simple cut 
$\gamma'$ and $\Gamma''$ contains a simple cut $\gamma''$ such that
$\gamma'$ and $\gamma''$ are of the same class, types
of $\Omega/\gamma'$ and $\Omega/\gamma''$ coincide,  $\gamma'$
and $\gamma''$ induces the same division of the set $\Omega_0$ (if applicable) and 
end points of $\gamma'$ and $\gamma''$ belong to the same
one-dimensional stratum
(if applicable). Then there exists a cut system $\overline{\Gamma}$
such that it is equivalent to
$\Gamma''$ and contains  cut $\gamma'$.
\end{corollary}

Proof is evident.

\begin{lemma} \label{36a} Let $\Gamma'$ and $\Gamma$ be two equivalent 
cut systems of a stratified surface $\Omega$.
Fix local orientations in a small neighborhoods of all special points.
Then there exists
an isomorphism $\phi:\Omega\to\Omega$ of the stratified surface such that
$\phi(\Gamma')=\Gamma$, $\phi$ fixes all special points of $\Omega$ and 
preserves an orientation of a small neighborhood of any special point.
\end{lemma}  

Proof is skipped.

\begin{lemma} \label{35}
Let $\Gamma'$ and $\Gamma''$ be two cut systems of
a stratified surface $\Omega$. Suppose both 
$\Gamma'$ and $\Gamma''$ contain the same simple cut  $\gamma$.
If cut systems $\Gamma'\setminus\gamma$ and
$\Gamma''\setminus\gamma$ 
of contracted cut surface $\Omega/\gamma$ are equivalent then
$\Gamma'$ and $\Gamma''$ are equivalent cut systems of $\Omega$.
\end{lemma}

\begin{proof} 
Suppose, cut systems  $\widetilde\Gamma'=\Gamma'\setminus\gamma$ and
$\widetilde\Gamma''=\Gamma''\setminus\gamma$ of $\Omega/\gamma$
are equivalent.
Let $(\Omega_*,\gamma_*,\tau_*)$ be the cut surface, corresponding to
the pair $(\Omega,\gamma)$. Denote by $\omega_*$ the image 
of $\gamma_*$ in $\Omega/ \gamma$. Clearly, $\omega_*$ consists of either two or 
one special point of contracted cut surface $\Omega/ \gamma$.
Denote by $U$ the joint of neighborhoods of points from $\omega_*$
such that $U\cap (\widetilde\Gamma'\cup\widetilde\Gamma'')=\emptyset$.
By lemma \ref{36a},  there exists a homeomorphism $\widetilde\phi:\Omega/\gamma
\to\Omega/\gamma$, that is identical on $U$ and such
that $\widetilde\phi(\widetilde\Gamma')=
\widetilde\Gamma''$. Thus there exist a homeomorphism 
$\phi_*:\Omega_*\to \Omega_*$ such that
$\phi_*(\widetilde\Gamma')=\widetilde\Gamma''$ and $\tau_*\phi_*=
\phi_*\tau_*$. Gluing by $\tau_*$ gives a homeomorphism
$\phi:\Omega\to\Omega$ such that $\phi(\gamma)=\gamma$,
$\phi(\widetilde\Gamma')=\widetilde\Gamma''$
and $\widetilde\phi$ fixes all special points.
\end{proof}

\begin{corollary} \label{36}
Let  $\Gamma'$ and $\Gamma''$ be two cut systems of
a stratified surface $\Omega$. Suppose $\Gamma^\circ$ is a cut system 
such that $\Gamma^\circ\subset\Gamma'$ and $\Gamma^\circ\subset\Gamma''$. 
If cut systems $\Gamma'\setminus\Gamma^\circ$ and 
$\Gamma'\setminus\Gamma^\circ$ 
of $\Omega/\Gamma^\circ$ are equivalent
cut systems of $\Omega_\#=\Omega/\Gamma_0$
 then $\Gamma'$ and $\Gamma''$ are 
equivalent cut systems of $\Omega$.
\end{corollary}

%%%%%%%%%%%%%%%%%%%%%%%%%%%%%%%%%%%

\subsection{Basic and simple surfaces}

\begin{definition} A stratified surface $\Omega$ is called 
{\it a stable surface} if and only if any  connected component $\Omega_i$ of
$\Omega$ is a nontrivial surface (i.e., $\mu(\Omega_i)>0$)

A cut system  $\Gamma=\{\gamma_1,\dots,\gamma_k\}$ of a stable stratified
surface $\Omega$ is calledd {\it a stable cut system} if and only if
the contracted cut surface  $\Omega_{\#}=\Omega/\Gamma$ is also stable.

A connected stratified surface $\Omega$ is called {\it a basic surface}
if and only if there is no nonempty stable cut systems of $\Omega$.

A cut system $\Gamma$ of a stable stratified surface $\Omega$ is
called {\it a complete cut system} if and only if any connected component 
of the contracted cut surface  $\Omega_\#=\Omega/\Gamma$ is a basic surface. 

\end{definition}

Let us formulate several elementary statements.
A complete cut system is a stable cut system.
Let $\Gamma=\{\gamma_1,\dots,\gamma_k\}$ be a cut system of
a stratified surface $\Omega$. Choose a subset  
$\Gamma'=\{\gamma_{i_1},\dots,\gamma_{i_l}\}$ 
of $\Gamma$. Then  $\Gamma'$  is a cut system of $\Omega$ and 
the image of $\Gamma\setminus\Gamma'$ in the contracted cut 
surface $\Omega_\#=\Omega/\Gamma'$ is a cut system of $\Omega_\#$. 
If $\Gamma$ is a stable (complete) cut system  then $\Gamma\setminus\Gamma'$
is a stable (complete) cut system of $\Omega_\#$. 
Conversely, a cut system $\Gamma''$ on the contracted cut surface
$\Omega_\#=\Omega/\Gamma'$ can be lifted to a cut system
$\Gamma=\Gamma'\sqcup\overline{\Gamma}''$, where 
$\overline{\Gamma}''\subset\Omega$ is obtained from $\Gamma''$
in two steps: first, take the preimage of $\Gamma''$ in 
cut surface $\Omega_*$; second, take image of the preimage 
under glueing map $\glue:\Omega_*\to\Omega$. If  $\Gamma''$ 
is a stable (complete) cut system  then $\Gamma$
is a stable (complete) cut system of $\Omega_\#$. 
  
\begin{lemma} A basic stratified surface $\Omega$ is isomorphic to 
one of the three stratified surfaces:
\\ (1) sphere $(S^2,p_1,p_2,p_3)$ with three interior special points 
$p_1$, $p_2$, $p_3$ ($\mu=1$);
\\ (2) triangle, i.e., disc $(D^2,q_1,q_2,q_3)$ with three boundary special 
points $q_1$, $q_2$, $q_3$ ($\mu=\frac{1}{2}$);
\\ (3) disc $(D,p,q)$ with one interior special point $p$ and one  
boundary special point $q$ ($\mu=\frac{1}{2}$).
\end{lemma}

Proof is elementary 

\begin{lemma} For any stable stratified surface $\Omega$ there exists
a complete cut system of $\Omega$.
\end{lemma}

\begin{proof}
It is sufficient to prove lemma  for connected surfaces.
If genus $g$ of a stratified surface $\Omega$ is not zero, then 
there exist a cut of a handle or a M\"obius cut. Using one of 
this cuts we reduce the question to a surface with less genus and, by
inductive arguments, to a surface of genus zero. After it one can use 
inductive arguments with respect to the number of special points 
plus number of boundary contours. For example, if $\Omega$ contains at least
two interior special points and is not a basic surface then a contour $\gamma$
separating a disc with two of these points off the rest of surface is a stable cut and
contracted cut surface $\Omega/\gamma$ is the disjoint union of a sphere 
with three interior points and a stable surface with less 
parameter of induction. We skip elementary details.  
\end{proof}

Let us describe isomorphism  class of complete cut systems for 
several stratified surfaces.

\begin{lemma} \label{3.9}
Denote by $N$ be  the number of isomorphism classes of complete cut
systems of a stratified  surfaces. Then  
\\ (1) $N=3$ for a sphere with four interior special points (Fig.3a);
\\ (2) $N=2$ for a Klein bottle with one interior special point (Fig.3b);
\\ (4) $N=2$ for a disc with four boundary special points (Fig.3c);
\\ (5) $N=2$ for a disc with two boundary special points and one interior
special point (Fig.3d);
\\ (6) $N=3$ for a disc with one boundary special point and two interior
special points (Fig.3e);
\\ (7) $N=2$ for a M\"obius band with one boundary special point (Fig.3f);
\\ (8) $N=2$ for a cylinder having one boundary special point on each 
boundary contour, two boundary special points in total (Fig.3g).
\end{lemma}

\medskip
\epsfxsize=11truecm
\begin{figure}[tbph]
\centerline{\epsffile{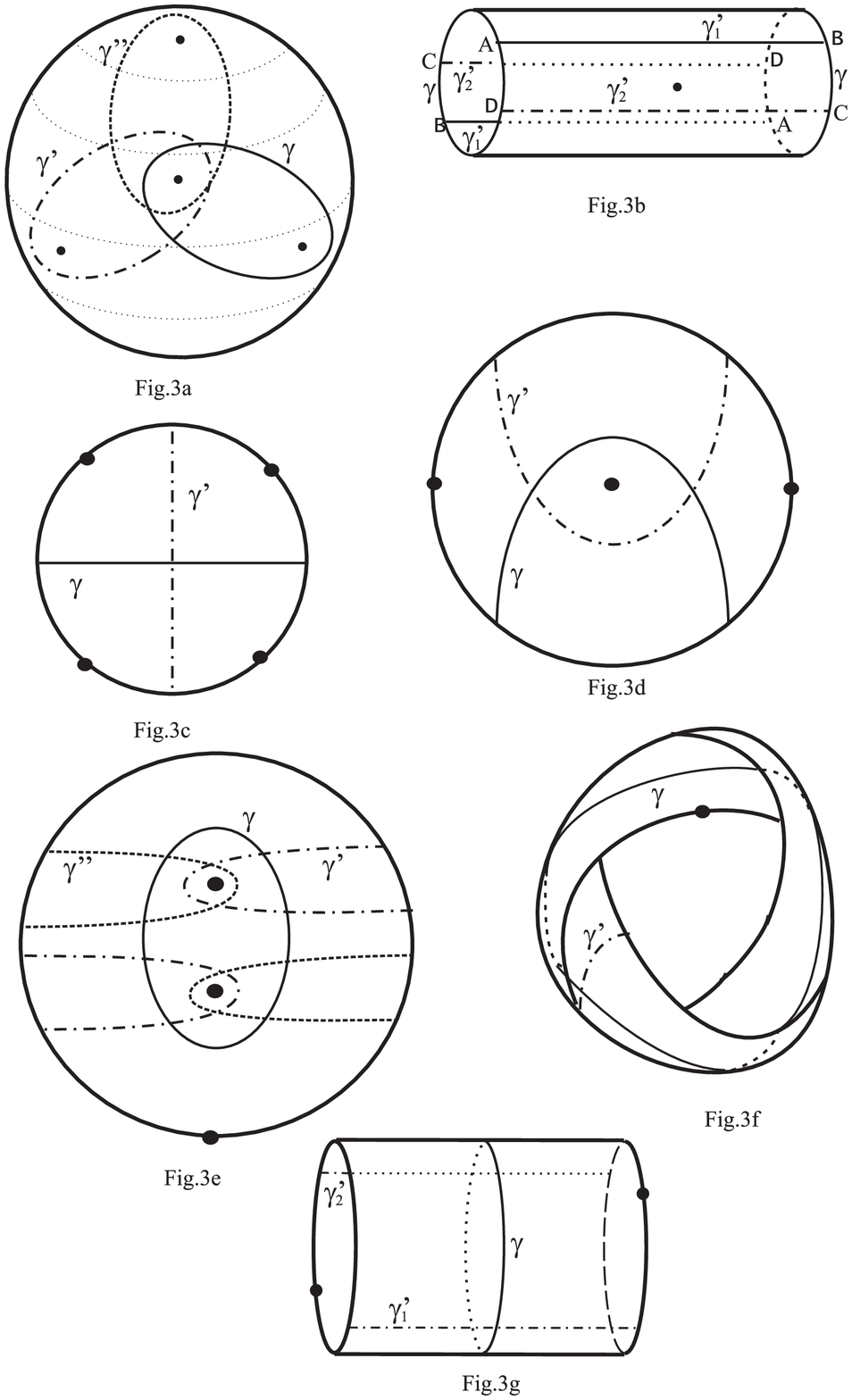}}
\caption{Non-equivalent complete cut systems on simple stratified surfaces.
\newline Different complete cut systems are drawn by diffent line types.
Cuts in a cut system are denoted by $\gamma$ with the same number
of accents and are distinguished by lower indeces.
Surfaces are numbered according to lemma \ref{3.9}.
Fig.3b presents Klein bottle cutted
by cut $\gamma$ forming one cut system.} 
\end{figure}
%\medskip

\begin{proof}

(1) Let $\Omega$ be a sphere with four interior special points.
Clearly, a complete cut system  $\Gamma$ of $\Omega$ consists of 
one separating contour that 
induces  the division of $\Omega_0$ into two pairs of special points. 
By lemma \ref{33} two contours inducing the same division are equivalent.
Therefore, there are three nonequivalent complete cut systems 
$\Gamma_d$ marked by a partition $d$ of the set $\{1,2,3,4\}$ into  
two-element subsets.

(2) Let $\Omega$ be a Klein bottle with one interior special point.
Any complete cut system of $\Omega$ contains a simple cut $\gamma$ 
such that the genus of $\Omega/\gamma$ is less than $1$. Clearly, only
two cases are possible: 

i) $\gamma=\gamma_{k}$, where $\gamma_k$ denotes  a cut of the 
neck of the Klein bottle. Clearly,
$\Gamma_1=\{\gamma_{k}\}$ is a complete cut system. 

ii) $\gamma=\gamma_{m}$, where $\gamma_m$ is a M\"obius cut. 
Clearly, there is a complete cut system 
$\Gamma_2=\{\gamma_{m},\gamma_{m}'\}$ consisting of  two M\"obius cuts 
and all complete cut systems consisting of two M\"obius cuts are equivalent.

(3) Let $\Omega$ be a surface of genus $g=\frac{3}{2}$
(i.e., a torus with a hole glued
by a M\"obius band) having no special points. 
There are the following  stable simple cuts (up 
to equivalence) of $\Omega$.
\begin{itemize}
\item[$\circ$] 
M\"obius cut $\gamma_{m1}$ such that $\Omega/\gamma_{m1}$ is 
a torus with one interior special point; 
\item[$\circ$] 
M\"obius cut $\gamma_{m2}$ such that $\Omega/\gamma_{m2}$ is 
a Klein bottle with one interior special point;
\item[$\circ$]
cut of a handle $\gamma_{h}$; in this case  $\Omega/\gamma_{h}$
is $\mathbb RP^2$ with two interior special points.
\end{itemize}             

Clearly, any complete cut system $\Gamma$  of $\Omega$ contains 
simple cut $\gamma$ reducing the genus of $\Omega$ by $\frac{1}{2}$ (i.e.,
genus of $\Omega/\gamma$ is equal to $1$). Hence, $\gamma$ is a M\"obius cut.

Suppose, $\gamma=\gamma_{m1}$. 
Then, evidently, $\Gamma=\Gamma_1$, 
where $\Gamma_1=\{\gamma_{m1},\gamma_{h}\}$ and 
all complete cut systems containing $\gamma_{m1}$ are equivalent.

Suppose, $\gamma=\gamma_{m2}$. Then $\Gamma\setminus\gamma$
is a complete cut system of a Klein bottle  $\Omega/\gamma$ having
one interior special point.  By (2) there are 
two isomorphism classes of complete cut systems of $\Omega/\gamma$. 
Therefore, either $\Gamma=\Gamma_2$, where 
$\Gamma_2=\{\gamma_{m2},\gamma_{m2}',\gamma_{m2}''\}$ or 
$\Gamma=\Gamma_3$, where 
$\Gamma_3=\{\gamma_{m2},\gamma_{h}\}$
(lifting of a cut of the neck of Klein bottle $\Omega/\gamma$ to $\Omega$ 
is a cut of a handle).
The latter is not equivalent to $\gamma_1$  because M\"obius cuts in 
them are not equivalent.

(3') 
Let $\Omega$ be a surface of genus ($g=\frac{3}{2}$) 
(i.e., a torus with a hole glued
by a M\"obius band) with one interior special point. 
There are the following  stable simple cuts
(up to equivalence) of $\Omega$.
\begin{itemize}
\item[$\circ$] 
M\"obius cut $\gamma_{m1}$ such that $\Omega/\gamma_{m1}$ is 
a torus with two interior special points; 
\item[$\circ$] 
M\"obius cut $\gamma_{m2}$ such that $\Omega/\gamma_{m2}$ is 
a Klein bottle with two interior special points;
\item[$\circ$]
separating contour $\gamma_{s1}$ such that $\Omega/\gamma_{s1}$ consists
of two connected components, namely,
a torus with one interior special point and $\mathbb RP^2$ with
two interior special points. 
\item[$\circ$] 
separating contour $\gamma_{s2}$ such that $\Omega/\gamma_{s2}$ consists
of two connected components, namely,
a Klein bottle with one interior special point and $\mathbb RP^2$ with
two interior special points;
\item[$\circ$]
cut of a handle $\gamma_{h}$; in this case  $\Omega/\gamma_{h}$
is $\mathbb RP^2$ with three interior special points.
\end{itemize}

Let $\Gamma$ be a complete cut system of $\Omega$. 

Suppose, $\Gamma$ contains a separating contour $\gamma_{s1}$. 
Then $\Omega/\gamma_{s1}$ is the disjoint union of a torus with
one interior special point and $\mathbb RP^2$ with two interior special points.
Clearly, there is only one equivalent class of complete cut systems of 
$\Omega/\gamma_{s1}$. The lifting $\Gamma_1$ of a complete cut system of
$\Omega/\gamma_{s1}$ to $\Omega$ consists of the following simple cuts:
$\Gamma_1=\{\gamma_{m1},\gamma_{s1},\gamma_{h}\}$.

Suppose, $\Gamma$ contains a separating contour $\gamma_{s2}$. 
Then $\Omega/\gamma_{s2}$ is the disjoint union of a 
a Klein bottle with one interior special point and $\mathbb RP^2$ with
two interior special points. By (2) there are two equivalency classes of 
complete cut systems of  $\Omega/\gamma_{s2}$. Their liftings to $\Omega$
are as follows:
$\Gamma_2=\{\gamma_{m2},\gamma_{s2},\gamma_{h}\}$ and 
$\Gamma_3=\{\gamma_{m2},\gamma_{m2}',\gamma_{m2}'',\gamma_{s2}\}$.

Suppose $\Gamma$ does not contain a separating contour and
contains a simple cut $\gamma_{m1}$.
Therefore, $\Omega/\gamma_{m1}$ is a torus with two interior special points.
It can be shown that in this case $\Gamma=\Gamma_4$, where 
$\Gamma_4=\{\gamma_{m1},\gamma_{h},\gamma_{h}'\}$.

Suppose, $\Gamma$ contains a simple cut $\gamma_{m2}$.
Therefore, $\Omega/\gamma_{m1}$ is a Klein bottle with two interior 
special points. It can be shown that in this case $\Gamma=\Gamma_5$, where 
$\Gamma_5=\{\gamma_{m2},\gamma_{h},\gamma_{h}'\}$.

(4) Let $\Omega$ be a  disc with four boundary special points.
Clearly, any complete cut system of $\Omega$ consists of a single 
simple cut $\gamma$, which is separating segment inducing the partition
of $\Omega_0$ into two pairs of consecutive boundary points. 
Mark boundary special points by $0,1,2,3$ going around boundary contour. 
There are two isomorphism classes of $\gamma$, namely,
$\Gamma_1=\{\gamma_{0,1}\}$, where $\gamma_{0,1}$ separates points $0,1$ of
$2,3$, and   $\Gamma_2=\{\gamma_{1,2}\}$, where $\gamma_{1,2}$ separates 
points $1,2$ of $3,0$.

(5) Let $\Omega$ be a
a disc with two boundary special points and one interior special point.
Clearly, any complete cut system of $\Omega$ consists of a single 
simple cut $\gamma$, which is separating segment, its end points belong 
to the same segment of the boundary and  $\gamma$ surrounds the interior 
special points. Thus, there are two isomorphism classes of complete cut 
systems of $\Omega$, namely,
$\Gamma_1=\{\gamma_l\}$ and $\Gamma_2=\{\gamma_r\}$, where 
end points of $\gamma_l$ (resp.  $\gamma_r$) belong to the 
'left' (resp. 'right') segment of the boundary contour.

(6) Let $\Omega$ be a
a disc with one boundary special point and two interior special points.
Clearly, there are three isomorphism classes of complete cut systems. 
First of them, denoted by  $\Gamma_1$, consists of a single separating
contour $\gamma_o$. Each of two other complete cut systems, 
$\Gamma_2$ and $\Gamma_3$ consist of two  
separating segments $\gamma_s$, $\gamma_s'$, such that  end points of 
any segment belong to the same segment of the boundary and  the segment  
surrounds one  interior special points. The only difference is 
the order of end points of two segments 
with respect to a fixed orientation of the boundary contour.

(7) Let $\Omega$ be a
a M\"obius band with one boundary special point. Clearly, there are two 
isomorphism classes of complete cut systems of $\Omega$. They are as follows.

$\Gamma_1=\{\gamma_m\}$, where $\gamma_m$ is a M\"obius cut of $\Omega$.

$\Gamma_2=\{\gamma_{mh}\}$, where $\gamma_{mh}$ is a cut across M\"obius band.

(8) Let $\Omega$ be a
a cylinder having one boundary special point on each boundary contour,
two boundary special points in total.
Clearly, there are two isomorphism
classes of complete cut systems of $\Omega$. They are as follows.

$\Gamma_1=\{\gamma_s\}$, where $\gamma_s$ is a separating contour.

$\Gamma_2=\{\gamma_{sh}\}$, where $\gamma_{sh}$ is a cut between 
two holes.

\end{proof}

\begin{definition}
Stratified  surfaces (1) -- (8)  from list in lemma \ref{3.9}
are called simple stratified surfaces.
\end{definition} 

%%%%%%%%%%%%%%%%%%%%%%%%%%%%%%%%%%%%

\subsection{Neighboring complete cut systems}

In this subsection we deal with complete cut systems of a fixed stratified surface.
It is supposed that any boundary contour of $\Omega$ contains at least one special
point. 

We call a pair of complete cut systems $\Gamma_1$ and
$\Gamma_2$ of $\Omega$  {\it an adjacent pair}
if the following condition holds.
Denote by $\Gamma_0$ the joint of all simple cuts that belongs
to $\Gamma_1$ and $\Gamma_2$ simultaneously and denote by 
$\Omega_{\#}$ the contracted cut surface $\Omega/\Gamma_0$.
Both $\Gamma_1\setminus\Gamma_0$ and $\Gamma_2\setminus\Gamma_0$ give 
rise to the complete cut systems of $\Omega_{\#}$. The condition is that 
cut systems $\Gamma_1\setminus\Gamma_0$ and $\Gamma_2\setminus\Gamma_0$
belong to the same connected component of $\Omega_{\#}$ and this component is 
a simple surface.

We call isomorphism classes $C_1$ and $C_2$ of complete 
cut systems  {\it an adjacent classes} if and only if
there exists representatives $\Gamma_1\in C_1$ and
$\Gamma_2\in C_2$  such that
$\Gamma_1$ and $\Gamma_2$ are adjacent complete cut systems.

We call two isomorphism classes $C'$, $C''$ of complete cut systems  
{\it a neighboring classes} if and only if there exists
a sequence of isomorphism class  $C_1,\dots,C_k$ such that 
$C_1=C'$, $C_k=C''$ and
$C_i$ is adjacent to $C_{i+1}$ for $i=1,\dots, k-1$.
We call two complete cut systems from neighboring isomorphism
classes  {\it neighbors}. 

We call a stable surface {\it an almost simple surface} if and only if 
it is isomorphic to one of the following stratified surfaces:
\\ (1) a disc with two or less interior special points
and an arbitrary number of boundary special points; 
\\ (2) a cylinder with at most one interior special points and an arbitrary
number of boundary special points; 
\\ (3) a M\"obius band with at most one interior special points and an arbitrary
number of boundary special points. 
\\ (4) a surface of genus ($g=\frac{3}{2}$) (i.e., a torus with
a hole glued by a M\"obius band) having no special points or  
with one interior special point. 

\begin{lemma} \label{4.1} All complete cut systems of an almost simple
stratified surface $\Omega$ are neighbors.
\end{lemma}

\begin{proof} If $\Omega$ is a disc without interior special points then
the statement of lemma can be reformulated in terms of triangulations
of a polygon and follows from the connectedness of an associhedron
\cite{Dev}.

Let $\Omega$ be a disc with one interior special point.
Then the statement of lemma in this case follows from the following 
claims. 

1) Any complete  cut system of $\Omega$ contains a unique simple cut 
$\gamma$ such that $\Omega/\gamma$ consists of two connected components 
and one of components is a disc with one interior and one boundary
special points. Let us call a cut of this type a type H cut. Obviously
end points of $\gamma$ belong to the same boundary segment. 

2) All complete cut systems of $\Omega$
containing fixed cut $\gamma$ of type H are neighbours. 

3) A complete 
cut system $\Gamma$ containing $\gamma$ of type H is adjacent to 
a complete cut system
$\Gamma'$ containing a cut $\gamma'$ of type H such that its end points 
lie on the neighbor segment to a segment containing end points of $\gamma$.

Let $\Omega$ be a disc with two interior special point.
In this case the following claims leads to the proof of lemma.

1) A complete cut system $\Gamma$ contains either two cuts of type H or 
a contour separating a disc with two interior special points from
the rest of surface. We call the last cut {\it a cut of type} P. 

2) A complete cut system containing a cut of type 
H is a neighbor of  a complete cut system containing a cut of type P.

3) All complete cut systems containing a fixed cut
$\gamma$ of type H are neighbors. 

4) A complete 
cut system $\Gamma$ containing $\gamma$ of type H is adjacent to a complete 
cut system $\Gamma'$ containing a cut $\gamma'$ of type H such that its 
end points lie on the neighbor segment to a segment containing end 
points of $\gamma$.

5)Let $\Omega$ be a surface of genus $g=\frac{3}{2}$
(i.e., a torus with a hole glued by a M\"obius band) having no special points. 
We will prove that there are just 3 nonequivalent complete cut systems  of 
$\Omega$ and all of them are neighbours.

There are the following  stable simple cuts (up 
to equivalence) of $\Omega$.
\begin{itemize}
\item[$\circ$] 
M\"obius cut $\gamma_{m1}$ such that $\Omega/\gamma_{m1}$ is 
a torus with one interior special point; 
\item[$\circ$] 
M\"obius cut $\gamma_{m2}$ such that $\Omega/\gamma_{m2}$ is 
a Klein bottle with one interior special point;
\item[$\circ$]
cut of a handle $\gamma_{h}$; in this case  $\Omega/\gamma_{h}$
is $\mathbb RP^2$ with two interior special points.
\end{itemize}             

Clearly, any complete cut system $\Gamma$  of $\Omega$ contains 
simple cut $\gamma$ reducing the genus of $\Omega$ by $\frac{1}{2}$ (i.e.,
genus of $\Omega/\gamma$ is equal to $1$). Hence, $\gamma$ is a M\"obius cut.

Suppose, $\gamma=\gamma_{m1}$. 
Then, evidently, $\Gamma=\Gamma_1$, 
where $\Gamma_1=\{\gamma_{m1},\gamma_{h}\}$ and 
all complete cut systems containing $\gamma_{m1}$ are equivalent.

Suppose, $\gamma=\gamma_{m2}$. Then $\Gamma\setminus\gamma$
is a complete cut system of a Klein bottle  $\Omega/\gamma$ having
one interior special point.  By (2) there are 
two isomorphism classes of complete cut systems of $\Omega/\gamma$. 
Therefore, either $\Gamma=\Gamma_2$, where 
$\Gamma_2=\{\gamma_{m2},\gamma_{m2}',\gamma_{m2}''\}$ or 
$\Gamma=\Gamma_3$, where 
$\Gamma_3=\{\gamma_{m2},\gamma_{h}\}$
(lifting of a cut of the neck of Klein bottle $\Omega/\gamma$ to $\Omega$ 
is a cut of a handle).
The latter is not equivalent to $\gamma_1$  because M\"obius cuts in 
them are not equivalent.

Finally, it can be easely shown, that $\Gamma_2$ is a neighbor of $\Gamma_3$ and 
$\Gamma_1$ is a neighbor of $\Gamma_3$. Indeed, both cut systems in a pair have 
equivalent simple cut.
 
(5') 
Let $\Omega$ be a surface of genus ($g=\frac{3}{2}$) 
(i.e., a torus with a hole glued by a M\"obius band) with one interior special point. 

We will prove that there are just 5 nonequivalent complete cut systems  of 
$\Omega$ and all of them are neighbours.

There are the following  stable simple cuts
(up to equivalence) of $\Omega$.
\begin{itemize}
\item[$\circ$] 
M\"obius cut $\gamma_{m1}$ such that $\Omega/\gamma_{m1}$ is 
a torus with two interior special points; 
\item[$\circ$] 
M\"obius cut $\gamma_{m2}$ such that $\Omega/\gamma_{m2}$ is 
a Klein bottle with two interior special points;
\item[$\circ$]
separating contour $\gamma_{s1}$ such that $\Omega/\gamma_{s1}$ consists
of two connected components, namely,
a torus with one interior special point and $\mathbb RP^2$ with
two interior special points. 
\item[$\circ$] 
separating contour $\gamma_{s2}$ such that $\Omega/\gamma_{s2}$ consists
of two connected components, namely,
a Klein bottle with one interior special point and $\mathbb RP^2$ with
two interior special points;
\item[$\circ$]
cut of a handle $\gamma_{h}$; in this case  $\Omega/\gamma_{h}$
is $\mathbb RP^2$ with three interior special points.
\end{itemize}

Let $\Gamma$ be a complete cut system of $\Omega$. 

Suppose, $\Gamma$ contains a separating contour $\gamma_{s1}$. 
Then $\Omega/\gamma_{s1}$ is the disjoint union of a torus with
one interior special point and $\mathbb RP^2$ with two interior special points.
Clearly, there is only one equivalent class of complete cut systems of 
$\Omega/\gamma_{s1}$. The lifting $\Gamma_1$ of a complete cut system of
$\Omega/\gamma_{s1}$ to $\Omega$ consists of the following simple cuts:
$\Gamma_1=\{\gamma_{m1},\gamma_{s1},\gamma_{h}\}$.

Suppose, $\Gamma$ contains a separating contour $\gamma_{s2}$. 
Then $\Omega/\gamma_{s2}$ is the disjoint union of a 
a Klein bottle with one interior special point and $\mathbb RP^2$ with
two interior special points. By (2) there are two equivalency classes of 
complete cut systems of  $\Omega/\gamma_{s2}$. Their liftings to $\Omega$
are as follows:
$\Gamma_2=\{\gamma_{m2},\gamma_{s2},\gamma_{h}\}$ and 
$\Gamma_3=\{\gamma_{m2},\gamma_{m2}',\gamma_{m2}'',\gamma_{s2}\}$.

Suppose $\Gamma$ does not contain a separating contour and
contains a simple cut $\gamma_{m1}$.
Therefore, $\Omega/\gamma_{m1}$ is a torus with two interior special points.
It can be shown that in this case $\Gamma=\Gamma_4$, where 
$\Gamma_4=\{\gamma_{m1},\gamma_{h},\gamma_{h}'\}$.

Suppose, $\Gamma$ contains a simple cut $\gamma_{m2}$.
Therefore, $\Omega/\gamma_{m1}$ is a Klein bottle with two interior 
special points. It can be shown that in this case $\Gamma=\Gamma_5$, where 
$\Gamma_5=\{\gamma_{m2},\gamma_{h},\gamma_{h}'\}$.

Like in the case (5), it can be easely shown, that $\Gamma_1$ is a neighbor of $\Gamma_2$,
$\Gamma_1$ is a neighbor of $\Gamma_4$ and 
$\Gamma_2$, $\Gamma_3$, $\Gamma_5$ are neighbors.
Indeed, both cut systems in each pair have equivalent simple cut.

6) For the rest of almost simple surfaces the proof is by similar arguments,
we skip them here.

\end{proof}

A pair of complete cut systems is said to be parallel if they contain
at least one common simple cut. A pair of isomorphism classes of complete
cut systems is said to be parallel if there exists parallel representatives
of the classes.

\begin{lemma} \label{3.11}
All complete cut systems of an orientable surface $\Omega$ without a boundary are
neighbors. 
\end{lemma}

\begin{proof} We use inductive arguments with respect to genus $g$ of
the surface $\Omega$. The statement of lemma is known for the 
sphere \cite{Kee}. If $g>0$ then any complete  cut system $\Gamma$
contains a cut of a handle. By lemmas \ref{33} and corollary \ref{34a},
any pair of isomorphism classes of complete cut 
systems is parallel. The proof follows from lemma \ref{35} 
and  the inductive arguments.
\end{proof}

A simple cut  $\gamma$ of $\Omega$ is called {\it projective} if it is a
M\"obius cut and $\Omega/\gamma$ is a nonorientable surface.

\begin{lemma} \label{3.12}
All complete cut systems of a nonorientable surface $\Omega$ without a boundary are
neighbors. 
\end{lemma}

\begin{proof} We  use inductive arguments with respect to genus
$g$ of the surface $\Omega$. Let $\Gamma$ be a complete cut system 
of $\Omega$.

If $g=\frac{1}{2}$ then any complete cut system of $\Omega$ contains
a M\"obius cut. By lemma \ref{3.11} and \ref{35}, we obtain that all 
complete cut systems are neighbors. Similar arguments are valid for 
$g=1$. Using a slightly more complicated reasoning one can prove 
lemma for $g=\frac{3}{2}$.

Let $g>\frac{3}{2}$ and lemma is proven for smaller genera. 
It can be shown, that  any  complete cut system $\Gamma$ 
contains a simple cut  $\gamma$ such that genus $g_\#$ of 
$\Omega/\gamma$ is less than the genus of $\Omega$ and 
$\Omega/\gamma$ is nonorientable. By inductive arguments,
$\Gamma\setminus\gamma$ is equivalent to a complete cut system 
containing a projective cut. By lemma \ref{35}, $\Gamma$ is equivalent
to a complete cut system containing a projective cut. Therefore,  all 
isomorphism classes of complete cut systems of $\Omega$ are parallel. 
The proof of lemma is completed by inductive arguments using lemma \ref{35}.
\end{proof}

A stable cut $\gamma$ of a connected stratified surface $\Omega$
is called  {\it normal} if it is either a contour 
or a separating segment such that one of the connected components
of a contracted cut surface is homeomorphic to a disk without
interior special points (no restrictions for boundary special points). 
If $\Omega$ is disconnected stratified surface then $\gamma$ is
said normal if it is a normal cut of a connected component of $\Omega$.

Clearly, any simple cut of a stratified surface without boundary is normal. 

We call a simple cut, that is not normal cut,  a {\it special} cut.

\begin{definition} A complete cut system $\Gamma$ is called normal 
cut system if $\Gamma$ consists of  normal simple cuts  only.
\end{definition}

For example, any complete cut system of a stratified surface
without boundary is normal.

Clearly, if $\Gamma$ is a normal cut system of $\Omega$ and $\Gamma_0\subset\Gamma$
is a subsystem of cuts then $\Gamma\setminus\Gamma_0$ is normal 
cut system of $\Omega/\Gamma_0$.

\begin{lemma} \label{4.6}
Any stable surface admits a normal complete cut system. 
\end{lemma}

\begin{proof} 
The proof is by inductive arguments with respect to the genus and the 
number of boundary contours. The step of induction is from 
$\Omega$ to $\Omega/\gamma$ where $\gamma$ is a  normal cut.
\end{proof}

\begin{lemma} \label{4.7}
All normal cut systems of a stratified surface $\Omega$ are neighbors. 
\end{lemma}

\begin{proof} 
If $\Omega$ is a surfaces without a boundary then the lemma follows from lemmas
\ref{3.11} and \ref{3.12}. 

If $\Omega$ is a disk without interior special points
then the lemma follows from lemma \ref{4.1}. 

In order to prove the lemma for an arbitrary $\Omega$ we 
use inductive arguments with respect to the number of boundary
contours. Let $\Gamma'$ and $\Gamma''$
be normal cut systems of $\Omega$.
Let $\omega$ be a boundary contour of $\Omega$. 
It  can be  easily shown that 
any  normal cut system  contains a separating oval  $\gamma$ 
with the following properties. The connected component of 
 $\Omega/\gamma$ that contains $\omega$ is 
homeomorphic to a disc with one interior special point (coming from $\gamma$).
The homotopy class of $\gamma$  does not depend on the choice of a normal 
cut system. Therefore, all isomorphism classes of normal cut systems of $\Omega$ 
are parallel. The proof of the lemma is completed by inductive arguments.
\end{proof}

\begin{lemma} \label{4.8}
Let  $\Omega$ be a stable surface and $\gamma$ be 
a stable special cut of $\Omega$. 
Suppose, $\Omega$ is not a simple or almost simple surface.
Then  there exists a normal cut $\gamma'$ homeomorphic 
to an oval such that $\{\gamma,\gamma'\}$  
is a stable cut system of $\Omega$. 
\end{lemma}

\begin{proof} 

If contracted cut surface 
$\Omega_{\#}=\Omega/\gamma$ admits a stable simple cut homeomorphic
to an oval then the claim of lemma is obviously true.
Suppose,  $\Omega_{\#}$ does not admit a simple cut with this property. 
Then any connected component of $\Omega_{\#}$ is either
a basic surface or a disk with no more than one interior
special point. It can be proven by exhaustion that 
$\Omega$ is either a simple or  almost  simple surface.
\end{proof} 

Let $\Gamma$ be a complete cut system of $\Omega$. Denote 
by $\Gamma_n$ the subset of all normal cuts from $\Gamma$.
$\Gamma$ is called {\it an almost normal} if 
one connected component of $\Omega/\Gamma_n$ is either a
simple or an almost simple surface  and all other connected components 
are basic surfaces.

By lemmas \ref{4.1} and \ref{4.6}, we obtain the following statement.

\begin{lemma} \label{4.9} Any almost normal cut system
has a neighbor that is a  normal cut system.
\end{lemma}

\begin{lemma} \label{4.10}
Any special stable cut $\gamma$ can be included in  an 
almost normal complete cut system. 
\end{lemma}

\begin{proof} 
Use inductive arguments with respect to genus, number
of boundary contours and number of interior special points. 

If $\Omega$ is either a simple or an almost simple surfaces
then the claim of lemma is trivial.

Otherwise choose a normal simple cut $\gamma'$ 
provided by lemma \ref{4.8}. By inductive hypothesis,
there exists an almost normal cut system $\Gamma_\#$ of 
the contracted cut surface $\Omega/\gamma'$ such that $\gamma\subset\Omega_\#$.
Clearly, $\Gamma_\#\sqcup\gamma'$ is an almost simple 
cut system of $\Omega$.
\end{proof}

\begin{theorem} \label{th4.1}
Let $\Omega$ be an arbitrary stable stratified surface. Then
all complete cut systems of $\Omega$ are  neighbors. 
\end{theorem}

\begin{proof}
By lemma \ref{4.7}, it is sufficient to prove that any complete cut system 
$\Gamma$ is the neighbor of a normal cut system. 
Use inductive arguments with respect to the number of  special
cuts. Suppose, $\Gamma$ contains 
$n>0$ special cuts. Let $\gamma\in\Gamma$ be 
a special cut. By lemma \ref{4.10}, there exists an almost normal cut system
$\Gamma'$ such that $\gamma\subset\Omega'$.

By inductive hypothesis, complete cut systems 
$\Gamma\setminus\gamma$ and $\Gamma'\setminus\gamma$ 
of the contracted cut surface $\Omega/\gamma$ are neighbors.
Therefore, $\Gamma$ and  $\Gamma'$ are neighbors. 
By lemma \ref{4.9}, $\Gamma'$ is the neighbor of a normal cut system.
Hence,  $\Omega$ is the neighbor of a normal cut system.
\end{proof}

%%%%%%%%%%%%%%%%%%%%%%%%%%%%%%%%%%%%%%%%%%%%%

\section{Two-dimensional topological field theory}
\label{sec4}

In this section we
use results of sections \ref{sec1} and \ref{sec3} in order to define a Klien topological field theory (subsection 4.1), reformulate it 
in terms of systems of correlators as it is usually done 
in physical literature (subsections 4.2) and prove the main theorem  \ref{MMain} , which states 
the correspondence between KTFT and structure algebras (subsection 4.3, 4.4).
As a corollary we get an analog of this theorem for open-closed topological field
theoris and prove that any massive open-closed topological field
theory can be extended to a Klein topological field theory (subsection 4.5).

\subsection{Definition of Klein topological field theory}

First, we shall fix a tensor category of surfaces (a 'basic' category $2\mathcal D$) 
and a functor from it to the tensor category of vector spaces.

\smallskip
{\it A set of local orientations.}
For any special point $r\in Q$  fix an orientation $o_r$ 
of its small neighborhood. Denote the set of these local orientations 
by $\mathcal O$. Let $\Omega$ be a connected stratified surface.
A set $\mathcal O$ of local orientations is called {\it admissible} if and only if
either $\Omega$ is orientable surface and all local orientations
are induced by an orientation of $\Omega$ (Fig.4a) or $\Omega$
is a nonorientable surface and local orientations at all special
points from any boundary contour $\omega_i$ are compatible
with one of orientations of $\omega_i$ (Fig.4b).
Thus, there are two admissible sets of local
orientations in the former case and  there are $2^s$, where 
$s$ is a number of connected components of $\partial\Omega$, 
admissible sets of local orientations in the latter case.

A set of local orientations for disconnected $\Omega$ is
said admissible if and only if it is admissible for each
connected component of $\Omega$.

\medskip
\begin{figure}[tbph]
\centerline{\epsffile{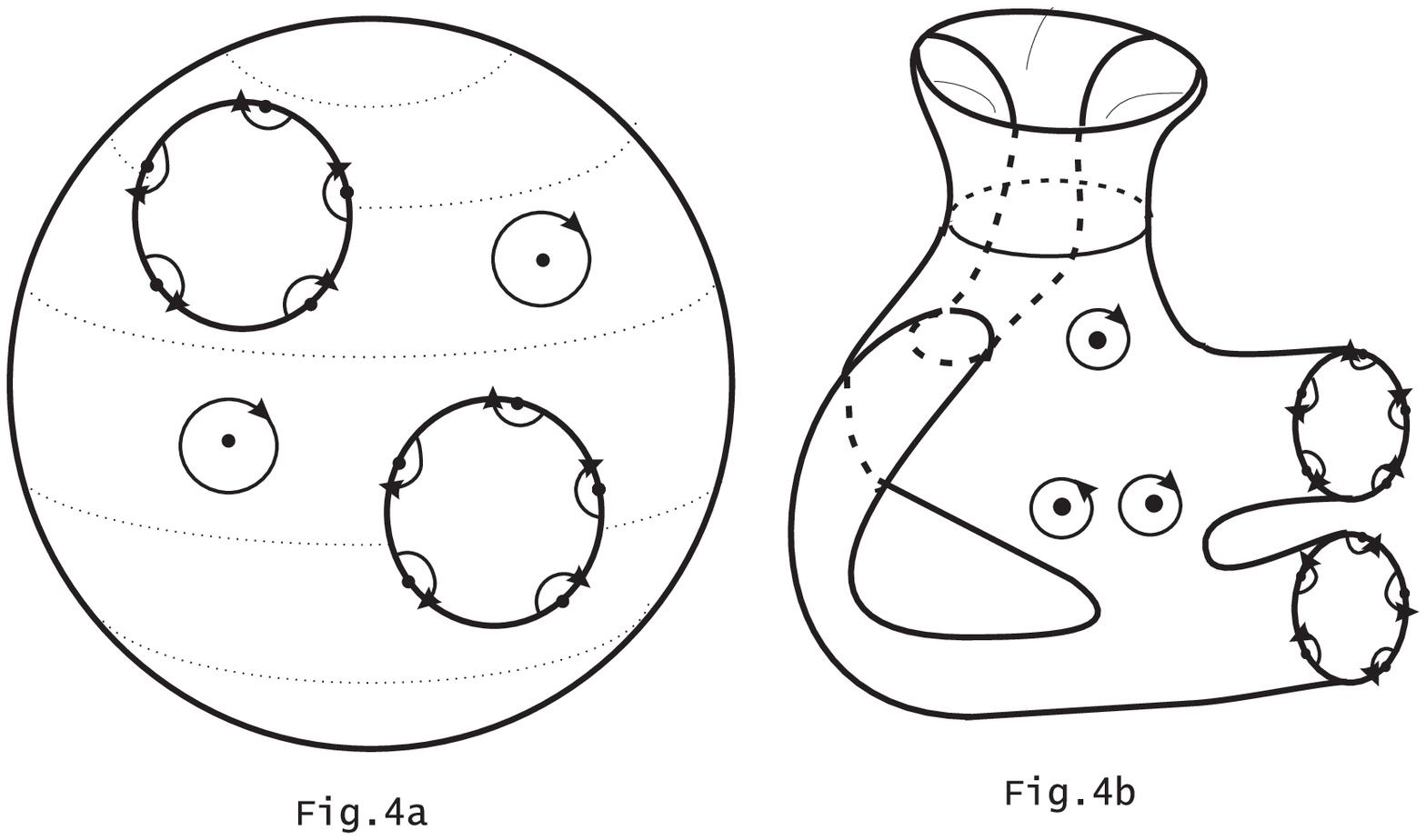}}
\caption{Admissible set of local orientations of (4a) oriented  
surface; (4b)  non-oriented surface.}
\end{figure}
\medskip

\begin{lemma} Let $(\Omega',\mathcal O')$ and $(\Omega'',\mathcal O'')$
be two pairs consisting of stratified surfaces and admissible sets of local 
orientations at special points. If stratified surfaces $\Omega'$ and $\Omega''$
are isomorphic then there exists an isomorphism $\phi:\Omega'\to\Omega''$ 
such that $\phi(\mathcal O')=\mathcal O''$.
\end{lemma}

Proof follows from standard properties of surfaces.

Let $\Omega$ be a stratified surface. Denote by
$\Omega_a$ the set of all interior special points 
and by $\Omega_b$ the set of all boundary special points.
Put also $\Omega_0=\Omega_a\sqcup\Omega_b$.
Fix a local orientation $o_r$ in a small neighborhood of any
special point $r\in \Omega_0$ and denote by $\mathcal O$
the set of all these local orientations.
Pairs $(\Omega,\mathcal O)$ are {\it objects of the basic category $2\mathcal D$}.
{\it Morphisms} are any combinations of morphisms of types 1)-4).
\smallskip

1) {\it Isomorphism} $\phi:(\Omega,\mathcal O)\to(\Omega',\mathcal O')$.
By definition, $\phi$ is an isomorphism $\phi:\Omega\to\Omega'$ of 
stratified surfaces compatible with local orientations at special points.    
\smallskip

2) {\it Changing of local orientations} $\psi:(\Omega,\mathcal O)\to(\Omega,\mathcal O')$.
Thus,  there is one such morphism for 
any pair $(\mathcal O,\mathcal O')$ of sets of local orientations on a stratified surface 
$\Omega$.
\smallskip

3) {\it Adding a special point} $\xi:(\Omega,\mathcal O)\to(\Omega',\mathcal O')$.
Morphism $\xi$ depends on a point $r\in\Omega\setminus \Omega_0$
endowed with a local orientation $o_r$. By definition, topological surface 
$\Omega'$ coincides with surface $\Omega$.
The stratification of $\Omega'$ is defined as  a refinement 
of the stratification of $\Omega$ by additional point $r$. 
Therefore, $\Omega_0'=\Omega_0\sqcup \{r\}$. Local orientation $o_r$
completes the set of local orientations, i.e., $\mathcal O'=\mathcal O\sqcup \{o_r\}$. 
\smallskip

4) {\it Cutting} $\eta:(\Omega,\mathcal O)\to(\Omega_\#,\mathcal O_\#)$.
Morphism $\eta$ depends on a cut system $\Gamma$ endowed with 
orientations of all simple cuts $\gamma\in\Gamma$.  $\Omega_\#$ is defined as 
contracted cut surface $\Omega/\Gamma$. Stratified surface $\Omega_\#$
inherits special points of $\Omega$ and local orientations at any of them.
A special point $r\in (\Omega_\#)_0\setminus\Omega_0$ appears as contracted 
connected component of a cut contour $\gamma_*$ in cut manifold $\Omega_*$.  
In obvious way the orientation of simple cut $\gamma$ induces the local orientation at 
$r$.
The set of all these local orientations at all points of
$(\Omega_\#)_0$ form the set $\mathcal O_\#$.

\smallskip

Disjoint union provides monoidal structure on the category $2\mathcal D$.

%Define {\it  'tensor product'}
%$\theta: (\Omega',\mathcal O')\times(\Omega'',\mathcal O'')\to
%(\Omega,\mathcal O)$ of two pairs $(\Omega',\mathcal O')$ and
%$(\Omega'',\mathcal O'')$ as their disjoint union
%$(\Omega,\mathcal O)=(\Omega'\sqcup \Omega'',\mathcal O'\sqcup \mathcal O'')$.

We shall define below a functor $(\Omega,\mathcal O) \to V(\Omega,\mathcal O)$ 
from the basic category of surfaces to the category of vector spaces.

Let $\{X_m | m\in M\}$ be a finite set of $n=|M|$  vector spaces $X_m$.
The action of symmetric group $S_n$ on the set $\{1,\dots,n\}$ induces
its action on  the linear space  
$\left(\oplus_{\sigma:\{1,\dots,n\}\leftrightarrows M}
X_{\sigma(1)}\otimes\dots\otimes X_{\sigma(n)}\right)$, an element 
$s\in S_n$ takes $X_{\sigma(1)}\otimes\dots\otimes X_{\sigma(n)}$
to $X_{\sigma(s(1))}\otimes\dots\otimes X_{\sigma(s(n))}$.
Denote by $\otimes_{m\in M} X_m$ the subspace of all invariants
of this action.

Vector space  $\otimes_{m\in M} X_m$ is canonically isomorphic to a
tensor product of all $X_m$ in any fixed order, the isomorphism
is a projection of vector space $\otimes_{m\in M} X_m$ to the summand
that is equal to the tensor product of $X_m$ in a given order.
Assume that all $X_m$ are equal to a fixed vector space $X$.
Then any bijection $M\leftrightarrow M'$ of sets induces the isomorphism 
$\otimes_{m\in M} X_m\leftrightarrow \otimes_{m'\in M'} X_{m'}$.

Let $A$ and $B$  be finite dimensional vector spaces over complex numbers
endowed with involutive transformations $A\to A$ and $B\to B$, which we denote by
$x\mapsto x^*$ ($x\in A$) and $y\mapsto y*$ ($y\in B$) resp.

Let $(\Omega,\mathcal O)$ be a pair consisting of a stratified surface $\Omega$ 
and a set $\mathcal O$ of local orientations at its special points.
Assign a copy $A_p$ of a vector 
space $A$ to any point $p\in\Omega_a$ and 
a copy $B_q$ of a vector space $B$ to any point $q\in\Omega_b$.
Elements of $A_p$ (resp. $B_q$) are called interior (resp. boundary) fields.
Put $V_{\Omega}=(\otimes_{p\in \Omega_a} A_p)\otimes
(\otimes_{q\in\Omega_b} B_q)$.

For any morphism of pairs $(\Omega,\mathcal O)\to (\Omega',\mathcal O')$
define a morphism of vector spaces $V_{\Omega}\to V_{\Omega'}$ as follows.

1) An isomorphism $\phi:(\Omega,\mathcal O)\to(\Omega',\mathcal O')$
induces the isomorphism 
$\phi_*:V_{\Omega}\to V_{\Omega'}$ because
$\phi$ generates the bijections $\Omega_a\leftrightarrows\Omega_a'$
and $\Omega_b\leftrightarrows\Omega_b'$ of sets of special points.

2) For a changing of local orientations 
$\psi:(\Omega,\mathcal O)\to(\Omega,\mathcal O')$
define a linear map $\psi_*:V_{\Omega}\to V_{\Omega}$
as $(\otimes_{r\in \Omega_0} \psi_r)$, 
where for any $r\in \Omega_0$

$$\psi_r(x)=\left\{\begin{array}{ll} x & \mbox{ if $o_r=o_r'$}\\
                               x^* & \mbox{ if $o_r=-o_r'$}\\
		\end{array}
	\right.$$

3) In order to define 
a morphism $\xi_*:V_{\Omega}\to V_{\Omega'}$
for adding a special point $\xi:(\Omega,\mathcal O)\to(\Omega',\mathcal O')$
we need to fix elements $1_A\in A$ and $1_B\in B$. These notations of 
elements are motivated in the frame of a Klein topological 
field theory, where they are 'trivial fields'.

In this case there exists a canonical isomorphism
$V_{\Omega'}=V_{\Omega}\otimes X$, where $X$ is equal either 
to $A$ (if adding special point $r$ belongs to the interior of $\Omega$) or
to $B$ (if adding special point $r$ belongs to the boundary of $\Omega$).
For $x\in V_{\Omega}$ put $\xi_*(x)=x\otimes 1_X$, where $1_X$ is
either $1_A$ or $1_B$ resp.

4) In order to define 
a morphism $\eta_*:V_{\Omega}\to V_{\Omega_\#}$
for any cutting morphism $\eta:(\Omega,\mathcal O)\to(\Omega_\#,\mathcal O_\#)$
we need to fix elements $\widehat{K}_{A,*}\in A\otimes A$, 
$\widehat{K}_{B,*}\in B\otimes B$ and $U\in A$. (Notations of them will be clear from 
the sequel.) 

Evidently, it is sufficient to define $\eta_*$ for an arbitrary  
oriented simple cut $\gamma\subset \Omega$. In this case  we have 
a canonical isomorphism
$V_{\Omega_\#}=V_{\Omega}\otimes X$, where 
$$X=\left\{\begin{array}{ll}
A\otimes A & \mbox{ if $\gamma$ is an contour and is not M\"obius cut}\\
B\otimes B &  \mbox{ if $\gamma$ is a segment}\\
A & \mbox{ if $\gamma$ is a M\"obius cut}\\
\end{array}\right.$$

For $x\in V_{\Omega}$ put $\eta_*(x)=x\otimes z$, where $z$ is
either $\widehat{K}_{A,*}$, or $\widehat{K}_{B,*}$, or $U$ resp.
\smallskip

Finally, for the tensor product  
$(\Omega',\mathcal O')\otimes (\Omega'',\mathcal O'')$ of objects of the category
$2\mathcal D$
there is evident canonical linear map 
$\theta_*: V_{\Omega'}\otimes V_{\Omega''}\to  
V_{\Omega'\sqcup\Omega''}$

\bigskip
                                     
We say that a set of data 
$\mathcal T=(A,x\mapsto x^*,B,y\mapsto y^*,\{\Phi_{(\Omega,\mathcal O)}\})$,
where $A$ and $B$ are linear spaces, $x\mapsto x^*$ and $y\mapsto y^*$ are 
involutive  linear transformations of $A$ and $B$ resp.,
$\{\Phi_{(\Omega,\mathcal O)}\}$ is a system of linear forms
$\Phi_{(\Omega,\mathcal O)}: V_{\Omega}\to \mathbb C$, is
a {\it Klein topological field theory} 
if the following axioms are satisfied.
\smallskip

%%%%%%%%%%%%%%%%%%%%%%%%%%%%%%%%%%%%%%%

$1^\circ$ {\it Topological invariance.}  

For any isomorphism of pairs
$\phi:(\Omega,\mathcal O) \to (\Omega', \mathcal O')$ the following
identity holds

$$\Phi_{(\Omega,\mathcal O)}(x)=\Phi_{(\Omega',\mathcal O')}(\phi_*(x).)$$

\smallskip
$2^\circ$ {\it Invariance of a change of  local orientations.} 

For any change of local orientations
$\psi:(\Omega,\mathcal O)\to(\Omega,\mathcal O')$
the following identity holds 

$$\Phi_{(\Omega,\mathcal O)}(x)=
\Phi_{(\Omega,\mathcal O')}(\psi_*(x)).$$
\smallskip

$3^\circ$ {\it  Nondegeneracy.}

Define first a bilinear form $(x,x')_A$ on vector space $A$. Namely, 
let $(\Omega,\mathcal O)$ be a pair, where $\Omega$ 
is a sphere  with 2 interior special points $p,p'$ and the set 
    $\mathcal O=\{o_p,o_{p'}\}$ is such that 
local orientations $o_p$, $o_{p'}$ induce the same global orientations of 
the sphere. Put $(x,x')_A=\Phi_{(\Omega,\mathcal O)}(x_p\otimes x'_{p'})$,
where $x_p$ and $x'_{p'}$ are images of $x\in A$ and $x'\in A$ in $A_p$
and $A_{p'}$ resp. The correctness of this definition follows from axioms
$1^\circ$ and $2^\circ$. Evidently, $(x,x')_A$  is a symmetric bilinear form.

Similarly, define a bilinear form $(y,y')_B$ on vector space $B$ using 
a disc  with 2 boundary special points $q,q'$ instead of a sphere with two 
interior special points $p,p'$. As in the previous case local orientations 
$o_q$, $o_{q'}$ must induce the same global orientations of the disc.
Evidently, $(y,y')_B$  is a symmetric bilinear form.

Axiom is that both forms $(x,x')_A$ and $(y,y')_B$ are nondegenerate.

\smallskip
$4^\circ$ {\it Invariance of adding unit fields.}

Axioms $1^\circ - 3^\circ$ allows us to choose elements 
$1_A\in A$, $1_B\in B$. 
Indeed, let $\Omega$ be a sphere with one special point $p$ and $o_p$ 
be a local orientation at $p$. Then linear form $\Phi_{(\Omega,\{o_p\})}$
is an element of the vector space dual to $A$. Identify dual vector
space with $A$ by  means of nondegenerate bilinear form $(.,.)_A$.
Thus, we obtain an element of $A$, which we denote by $1_A$.
Analogously, starting from a disc with one boundary special point we obtain
an element of $B$, which we denote by $1_B$. An element $1_A$ (resp. $1_B$)
is called trivial interior (resp. boundary) field. We shall use just
these elements in the definition of morphisms $\xi_*$ of type 3).

Axiom is that for any adding of special point
$\xi:(\Omega,\mathcal O)\to(\Omega',\mathcal O')$ 
the following identity holds

$$\Phi_{(\Omega,\mathcal O)}(x)=
\Phi_{(\Omega',\mathcal O')}(\xi_*(x)).$$

\smallskip
$5^\circ$ {\it Cut invariance.}

Axioms $1^\circ-3^\circ$ allows us to choose elements 
$\widehat{K}_{A,*}\in A\otimes A$, $\widehat{K}_{B,*}\in B\otimes B$ and
$U\in A$. Indeed, any nondegenerate bilinear form 
on a vector space $X$ canonically defines an element
$\widehat{K}_X\in X\otimes X$.
Taking forms $(x,x')_{A,{*}}=(x,{x'}^{*})_A$ and
$(y,y')_{B,*}=(y,{y'}^{*})_B$
we obtain elements $\widehat{K}_{A,*}$ and $\widehat{K}_{B,*}$.

Linear form $\Phi_{(\Omega,\mathcal O)}$ for a projective plane $\Omega$ 
with one interior special point is an element of the vector space
dual to $A$. We denote by $U$ the image of this element
in $A$ under the isomorphism induced by nondegenerate bilinear form
$(x,x')_A$.

We shall use just these elements for morphisms $\eta_*$ of type 4).

Axiom is that for any cut system $\Gamma$ endowed with orientations 
of all simple cuts  the following identity is required 

$$\Phi_{(\Omega,\mathcal O)}(x)=
\Phi_{(\Omega_\#,\mathcal O_\#)}(\eta_*(x))$$

(here $\eta_*$ depends on $\Gamma$, see 4) above)

\smallskip
$6^\circ$  {\it Multiplicativity.}

Axiom is that the product $(\Omega',\mathcal O')\otimes
(\Omega'',\mathcal O'')$ of 
 any two pairs $(\Omega',\mathcal O')$ and $(\Omega'',\mathcal O'')$
the following identity holds 

$$\Phi_{(\Omega',\mathcal O')}(x_1)\Phi_{(\Omega'',\mathcal O'')}(x_2)=
\Phi_{(\Omega'\sqcup\Omega'',\mathcal O' \sqcup \mathcal O'')}
(\theta_*(x_1\otimes x_2))$$

\bigskip

\begin{definition} Klein topological field theories 
$\mathcal T=(A,x\mapsto x^*,B,y\mapsto y^*,\{\Phi_{(\Omega,\mathcal O)}\})$ and
$\mathcal T'=(A',x\mapsto x^{*'},B',y\mapsto y^{*'},\{\Phi'_{(\Omega,\mathcal O)}\})$
are called isomorphic if there exist isomorphisms $\varphi_A:A\to A'$, 
$\varphi_B:B\to B'$ of vector spaces such that 
$\varphi_A(x^*)=\varphi_A(x)^{*'}$, $\varphi_B(y^*)=\varphi_B(y)^{*'}$ and 
$\Phi'_{(\Omega,\mathcal O)}(\varphi_{\Omega}(z))=
\Phi_{(\Omega,\mathcal O)}(z)$  where $z\in V_{\Omega}$ and 
$\varphi_{\Omega}:V_{\Omega}\to V'_{\Omega}$ is the isomorphism induced by 
$\varphi_A$ and $\varphi_B$
\end{definition}

The example of a Klein topological field theory is given in section 5.

%%%%%%%%%%%%%%%%%%%%%%%%%%%%%%%%%%%%%%%%%%%%%%%%%%%%%%%%%%%%%

\subsection{Correlators of Klein topological field theory}

Axioms of a topological field theory 
can be reformulated in terms of so called 'systems
of correlators'. 
In terms of correlators the dependence of a linear form 
$\Phi_{(\Omega,\mathcal O)}$
on an individual stratified surface $\Omega$ endowed 
with a set of local orientations 
$\mathcal O$  is reduced to the dependence on a topological 
type $G=(g,\varepsilon,m_1,\dots,m_s)$
of a surface. 

In this section we present the reformulation of a Klein 
topological field theory in terms of correlators.
Let $A$ and $B$ be two linear spaces, $x\mapsto x^*$  and $y\mapsto y^*$ be
fixed involutive transformations of $A$ and $B$ resp.

Denote by $G=(g,\varepsilon,m,m_1,\dots,m_s)$ a data that are 
the topological type of a connected stratified surface. 

We assume throughout this  section that 
\begin{itemize}
\item[$\circ$] $m_i>0$ for $i=1,\dots,s$ ($m$ may be zero; $s$ may be zero) 
\item[$\circ$] $m+\sum m_i>0$
\item[$\circ$] numbers $m_i$ are given in ascending  ordered   $m_1\le m_2\le\dots\le m_s$. 
\end{itemize}

Denote by $V_G$  tensor product 
$A^{\otimes m}\otimes B^{\otimes (m_1+ \dots + m_s)}$ of vector spaces.

For any $G$ fix one linear function $f_G:V_G\to\mathbb C$. According
to the style of notations used in field theories we denote
corresponding polylinear function as follows: 
$$\langle x_1,\dots,x_m, (y_1^1,\dots, y_{m_1}^1),\dots, (y_1^s,\dots,
y_{m_s}^s)\rangle_{(g,\varepsilon)}:=f_G(z),$$
where $z=x_1\otimes\dots\otimes x_m\otimes y_1^1\otimes \dots\otimes
y_{m_1}^1 \otimes \dots\otimes  y_1^s\otimes \dots\otimes y_{m_s}^s\in V_G$.

Note that vector spaces $V_G$ and $V_{G'}$ coincide for some $G\ne G'$ but 
corresponding polylinear functions  'remember' all data 
from $G=(g,\varepsilon,m,m_1,\dots,m_s)$ due to bracketing of arguments. 
We call a bracket  $(y_1^j,\dots, y_{m_j}^j)$ {\it a block} of arguments.

We use capital letters $X, X_1$ etc. to denote  subsets consisting of 
several $x_i$ and several blocks of arguments. 
For example, $\langle x_1, X \rangle_{g,\varepsilon}$,
where 
\\$X=(x_2,\dots,x_m, (y_1^1,\dots,y_{m_1}^1),\dots, (y_1^s,\dots, y_{m_s}^s))$ 
denotes a function 
\\$\langle x_1,\dots,x_m, (y_1^1,\dots, y_{m_1}^1), \dots, (y_1^s, \dots, y_{m_s}^s) \rangle_{g,\varepsilon}$. 

We write down $\langle X \rangle_G$
instead of $\langle X \rangle_{g,\varepsilon}$
in order to emphasize explicitly the dependence  on a type $G$.

We denote by $X_a$ the set of all  $x_i$ in a correlator  and
we denote by $X_b$ the set of all blocks of arguments. Thus, we write down

$\langle X\rangle_G=\langle X_a,X_b \rangle_G=
\langle x_1,\dots,x_m, (y_1^1,\dots, y_{m_1}^1), \dots, (y_1^s, \dots,
y_{m_s}^s) \rangle_{g,\varepsilon}$

We use capital letters $Y, Y_1,$ etc. to denote 
sequential sets of arguments in a block. For example,  
$\langle X, (y_1^s,Y) \rangle_{g,\varepsilon}$, where 

$X=(x_1,\dots,x_m, (y_1^1,\dots,
y_{m_1}^1),\dots, (y_1^{s-1},\dots,y_{m_{s-1}}^{s-1}))$ and 
$Y=(y_2^s,\dots, y_{m_s}^s)$, denotes  polylinear function 

$\langle x_1,\dots,x_m, (y_1^1,\dots, y_{m_1}^1),\dots, (y_1^s,\dots,
y_{m_s}^s)\rangle_{g,\varepsilon}$. 

For an arbitrary $Y=(y_1,\dots,y_k)$ put $Y^*=(y_k^*,\dots,y_1^*)$.

Choose a basis $E_A$ of $A$ and $E_B$ of $B$ and denote basic elements by
$\alpha, \alpha', \dots \in E_A$  and $\beta,\beta',\dots\in E_B$. Then 
we can present $f_G$ as a tensor 
$$\langle \alpha_1,\dots,\alpha_m, (\beta_1^1,\dots, \beta_{m_1}^1),\dots, 
(\beta_1^s,\dots,\beta_{m_s}^s)\rangle_{(g,\varepsilon)}$$ 
\medskip

Let $\mathcal C=(A,x\mapsto x^*,B,y\mapsto y^*, 
\{\langle X \rangle_G\})$ be a set of data consisting of  
\begin{itemize}
\item[$\circ$] linear spaces $A$ and $B$, 
\item[$\circ$] involutive linear transformations $x\mapsto x^*$
($x\in A$) and $y\mapsto y^*$ ($y\in B$) of
$A$ and $B$ resp.,      
\item[$\circ$] polylinear function  $\langle X \rangle_G$
for any topological type  $G=(g,\varepsilon,m,m_1,\dots,m_s)$
of a stratified surface .
\end{itemize}

\begin{definition} \label{defcor}
A set of data 
$\mathcal C=(A,x\mapsto x^*,B,y\mapsto y^*, \{\langle X \rangle_G\})$
is called a system of correlators and any function $\langle X \rangle_G$
is called a correlator if and only if the following axioms hold.
\smallskip
\\ (i) Any correlator 
$\langle x_1,\dots,x_m, (y_1^1,\dots, y_{m_1}^1),\dots, (y_1^s,\dots,
y_{m_s}^s)\rangle_{g,\varepsilon}$ 
is invariant  with respect to a permutation of $x_i$, a permutation of blocks 
$(y_1^i,\dots, y_{m_i}^i)$ of equal size, a cyclic permutations of
arguments inside each block $(y_1^i,\dots, y_{m_i}^i)$.
\smallskip
\\ (ii) If $\varepsilon=1$ then a correlator is invariant with respect to 
replacing of  $x_i$ by $x_i^*$ and block $Y^j=(y_1^j,\dots, y_{m_j}^j)$  
by $(Y^j)^*$
for all $i=1,\dots,m$, $j=1,\dots,s$ simultaneously.

If $\varepsilon=0$ then a  correlator is invariant with respect to 
replacing of  $x_i$ by $x_i^*$ for any fixed $i$ and with respect to
replacing a block $Y^j=(y_1^j,\dots, y_{m_j}^j)$ by
$(Y^j)^*$ for any fixed $j$.
\\ (iii) Bilinear maps $\langle x_1,x_2 \rangle_{0,1}$ and 
$\langle (y_1,y_2)\rangle_{0,1}$ are nondegenerate bilinear forms 
on $A$ and $B$ resp.; below we denote them by $(x_1,x_2)_A$ and $(y_1,y_2)_B$
resp.
\smallskip
\\ (iv) Fix a basis $E_A=\{\alpha\}$ of $A$ and a basis $E_B=
\{\beta\}$ of $B$. Denote by $F^{\alpha',\alpha''}$ the tensor dual to the tensor 
$F_{\alpha',\alpha''}$ of bilinear form $(\alpha',\alpha'')_A$, 
by $F^{\beta',\beta''}$ the tensor dual to the tensor
$F_{\beta',\beta''}:=(\beta',\beta'')_B$ and by $D^\alpha$ the tensor obtained by lifting 
the index of  tensor $D_\alpha=\langle \alpha \rangle_{\frac{1}{2},0}$. 
Axiom is that 
for arbitrary sets $X$, $X_1$, $X_2$, $Y_1$, $Y_2$  of arguments denoted 
according to 
the conventions described above the  relations $(1)-(9)$ hold.
                                                                               
$$(1)\ \ 
\langle X_1,X_2 \rangle_{g_1+g_2,\varepsilon_1\varepsilon_2}=\sum_{\alpha',\alpha''\in {E_A}} 
	\langle X_1,\alpha' \rangle_{g_1,\varepsilon_1}F^{\alpha',\alpha''}
\langle \alpha'',X_2 \rangle_{g_2,\varepsilon_2}
$$

$$(2)\ \ 
\langle X\rangle_{g+1,\varepsilon}=
\sum_{\alpha',\alpha''\in{E_A}} \langle X,\alpha', \alpha''
\rangle_{g,\varepsilon}F^{\alpha',\alpha''}
$$

$$(3)\ \ 
\langle X\rangle_{g+1,0}=
\sum_{\alpha',\alpha''\in{E_A}} \langle X,\alpha',\alpha''{^*}
\rangle_{g,1}F^{\alpha',\alpha''}
$$

$$(4)\ \ 
\langle X\rangle_{g+\frac{1}{2};0}=
\sum_{\alpha'\in{E_A}} \langle X,\alpha' \rangle_{g,\varepsilon}
D^{\alpha'}
$$

$$(5)\ \ 
\langle X,(Y_1), (Y_2) \rangle_{g,\varepsilon}=
\sum_{\beta',\beta''\in{E_B}} 
\langle X,(Y_1,\beta',Y_2,\beta'') \rangle_{g,\varepsilon} 
F^{\beta',\beta''}                             
$$      

$$(6)  
\langle X_1,X_2, (Y_1,Y_2)\rangle_{g_1+g_2,\varepsilon_1\varepsilon_2}=
\sum_{\beta',\beta''\in{E_B}} 
\langle X_1,(Y_1,\beta') \rangle_{g_1,\varepsilon_1} 
F^{\beta',\beta''}
\langle X_2,(Y_2,\beta'') \rangle_{g_2,\varepsilon_2}                  
$$

$$(7)\ \ 
\langle X, (Y_1, Y_2) \rangle_{g+1,\varepsilon}=
\sum_{\beta',\beta''\in{E_B}} 
\langle X,(Y_1,\beta'),(Y_2,\beta'') \rangle_{g,\varepsilon} 
F^{\beta',\beta''}                              
$$          

$$(8)\ \ 
\langle X, (Y_1, Y_2) \rangle_{g+1,0}=
\sum_{\beta',\beta''\in{E_B}} 
\langle X,(Y_1,\beta'),(Y_2^*,\beta''{^*}) \rangle_{g,1} 
F^{\beta',\beta''}
$$          

$$(9)\ \ 
\langle X, (Y_1, Y_2) \rangle_{g+\frac{1}{2},0}=
\sum_{\beta',\beta''\in{E_B}} 
\langle X,(Y_1,\beta',Y_2^*,\beta''{^*}) \rangle_{g,\varepsilon} 
F^{\beta',\beta''}
$$          
         
\end{definition}
\bigskip

\note  We pay no attention to the order in which  blocks of
arguments are written in a correlator because due to  (i)
of the definition of systems of correlators all
orderings of blocks in ascending order 
of their size give the same function.

\note We shall define a correlator for a type  $\overline{G}$ of
disconnected surfaces.
First, fix any ordering of all types $G$ of 
connected stratified surfaces (use for example
the lexicographic ordering). 
A type  $\overline{G}$ of disconnected surfaces is a set
$\overline{G}=(G_1,\dots,G_k)$ of types 
$G_i=(g^{(i)},\varepsilon^{(i)},m^{(i)},m_1^{(i)},\dots,m_{s_i}^{(i)})$ 
of connected stratified surface.
We require types $G_i$ to be numbered according the ordering.

Put 
$V_{\overline{G}}=A^{\otimes(\sum m^{(i)})}\otimes B^{\otimes (\sum m_j^{(i)})}$.
Define polylinear function $\langle X\rangle_{\overline{G}}$ on $V_{\overline{G}}$
as a product of correlators 
$$
\langle X\rangle_{\overline{G}}=
\langle X^{(1)} \rangle_{G_1}\dots\langle X^{(k)} \rangle_{G_k}
$$

Here  $X=(X_a^{(1)},\dots,X_a^{(k)},X_b^{(1)},\dots,X_b^{(k)})$
and $X^{(i)}=(X_a^{(i)},X_b^{(i)})$.

A type of an arbitrary stratified surface possibly disconnected will
be denoted below by $G$.

\medskip

Fix a basis $E_A$ of $A$ and a basis $E_B$ of $B$ and use nondegenerate 
tensors $F_{\alpha',\alpha''}=(\alpha',\alpha'')_A$,
$F_{\beta',\beta''}=(\beta',\beta'')_B$  and  dual to  lower and raise indices.

Denote by $1_A$ the element 
$\langle \alpha'\rangle_{0,1}F^{\alpha',\alpha''}\alpha''\in A$. 
Denote by $1_B$ the element 
$\langle \beta'\rangle_{0,1}F^{\beta',\beta''}\beta''\in B$. 
Denote by $U$ the element 
$\langle \alpha'\rangle_{\frac{1}{2},0}F^{\alpha',\alpha''}\alpha''\in A$. 

\begin{lemma} \label{4.2}
The following identities hold for any acceptable sets
of arguments $X$, $Y$.
\\ (1) $\langle 1_A, X \rangle_{g,\varepsilon}=
\langle X \rangle_{g,\varepsilon}$
\\ (2) $ \langle X,(1_B,Y) \rangle_{g,\varepsilon}=
\langle X, (Y) \rangle_{g,\varepsilon}$
\\ (3) $\langle U_A, X \rangle_{g,\varepsilon}= 
\langle X \rangle_{g+\frac{1}{2},0}$
\end{lemma}

\begin{proof} (1) Use identity (1) from the definition \ref{defcor} and put 
in it $X_1=\emptyset,g_1=0,\varepsilon_1=1, X_2=X, g_2=g,\varepsilon_2=\varepsilon$.  

(2) Use identity (6) from the definition \ref{defcor} and put 
in it $X_1=\emptyset ,Y_1 = \emptyset ,g_1=0,\varepsilon_1=1, 
X_2=X, Y_2=Y, g_2=g,\varepsilon_2=\varepsilon$.  

(3) Use identity (1) from the definition and put 
in it $X_1=\emptyset,g_1=\frac{1}{2},\varepsilon_1=0, X_2=X, g_2=g,\varepsilon_2=\varepsilon$.  
\end{proof}

\note Let $\gamma$ be an unstable simple cut of a stratified surface 
$\Omega$. Then one of the connected components of 
contracted cut surface $\Omega_\#$ is a trivial surface. This component 
is one of the surfaces (4), (5), (6), (7), (11) from the list
of trivial surfaces given in lemma \ref{3.1} due to the assumption that
any component of boundary of $\Omega$ contains at least one special point. 
In cases (4), (5), (7) simple cut $\gamma$ induces one of the relations 
from lemma \ref{4.2}. In cases (6), (11) simple cut $\gamma$ induces trivial
relation, following from the definition of bilinear forms $(x',x'')_A$
and $(y',y'')_B$. 

\medskip

We shall assign a system of correlators to a Klein topological field
theory.

Let $\mathcal T=(A,x\mapsto x^*,B,y\mapsto y^*,\{\Phi_{(\Omega,\mathcal O)}\})$ 
be a Klein topological field theory. For any type $G=(g,\varepsilon,m,m_1,\dots,m_s)$
choose a connected stratified surface $\Omega$  of type $G$ and 
an admissible set of local orientations $\mathcal O$. Note that
for any boundary contour $\omega_i$ admissible set of local orientations
generates a cyclic ordering of special points from $\omega_i$.

Denote by  $M$  the set of multipliers of the tensor product 
$V_G=A^{\otimes m}\otimes B^{\otimes (m_1+ \dots + m_s)}$.
A bijection $\mathcal N:M \leftrightarrow \Omega_0$ between set  $M$
and set $\Omega_0$ of special points on a stratified surface $\Omega$
endowed with an admissible set of local orientations $\mathcal O$
is called {\it an admissible bijection} if  
\begin{itemize}
\item[$\circ$] multipliers $A$ correspond to interior special points;

\item[$\circ$] first $m_1$ multipliers $B$ correspond to special points from
boundary contour $\omega_1$ such that $|\omega_1|=m_1$, next $m_2$ multipliers
$B$ correspond to  special points from
boundary contour $\omega_2$ such that $|\omega_2|=m_2$ etc.;

\item[$\circ$] the linear order of special points on any $\omega_i$ induced by  the order of multipliers 
$B$ in $V_G$ is compatible with the cyclic order induced by $\mathcal O$.
\end{itemize}

Fix  an admissible bijection $\mathcal N:M \leftrightarrow \Omega_0$. 
Denote by $\phi_{\mathcal N}:V_G\to V_{(\Omega,\mathcal O)}$ the isomorphism 
of vector spaces induced by this bijection.

For any  element 
$z=x_1\otimes\dots\otimes x_m\otimes y_1^1\otimes \dots\otimes  y_{m_1}^1
\otimes \dots\otimes  y_1^s\otimes \dots\otimes y_{m_s}^s\in V_G$
put                                                              
$$\langle x_1,\dots,x_m (y_1^1,\dots,y_{m_1}^1),\dots,(y_1^s,\dots,y_{m_s}^s) \rangle_G 
=\Phi_{(\Omega,\mathcal O)} (\phi_{\mathcal N}(z)). $$
                                              
Thus, we obtain a linear form $\langle X \rangle_G$ for any 
topological type $G$ of a stratified surface.

\begin{lemma} Linear form $\langle X \rangle_{g,\varepsilon}$
does not depend on the choice of a stratified surface $\Omega$, 
an admissible set of local orientations $\mathcal O$ and 
an admissible bijection $\mathcal N$.
\end{lemma}                  
                           
\begin{proof} Lemma
follows immediately from topological invariance of Klein field theory.
\end{proof}

Denote by $\mathcal C(\mathcal T)$ the set of data 
$(A,x\mapsto x^*,B,y\mapsto y^*, \{\langle X \rangle_G\})$.

\begin{lemma} The set of data $\mathcal C(\mathcal T)$ is a system of correlators.
\end{lemma}

\begin{proof} Identities (i) - (iii) follows from appropriate axioms of a Klein
topological field theory. 

We will pay more attention to identity (iv).
Let $\gamma$ be a simple cut of a stratified surface $\Omega$ of type 
$G=(g,\varepsilon,m,m_1,\dots,m_s)$ with $m_i>0$. 
Then  contracted cut 
surface $\Omega_\#=\Omega/\gamma$ is either a connected surface 
or consists of two connected components 
$(\Omega_\#)_1$, $(\Omega_\#)_2$. Denote by $G_\#$ the type of $\Omega_\#$.
Thus, $G_\#=((G_\#)_1$, $(G_\#)_2)$ in the latter case.

Cut invariance of a Klein topological field theory leads to 
the relation between correlators $\langle X \rangle_G$ and 
$\langle X \rangle_{G_\#}$.In order to
establish this relation we need to fix 

admissible sets of local orientations $\mathcal O$ on $\Omega$ and 
$\mathcal O_\#'$ on $\Omega_\#$;

admissible bijections $\mathcal N:M\leftrightarrow \Omega_0$ 
and $\mathcal N_\#:M_\#\leftrightarrow (\Omega_\#)_0$.

(Note that  a cut morphism $\eta:(\Omega,\mathcal O)\to (\Omega_\#,\mathcal O_\#)$
may bring an admissible set of local orientations $\mathcal O$ to 
nonadmissible set of local orientations $\mathcal O_\#$.)  

We have the chain of morphisms                 

$V_G \stackrel{\phi_{\mathcal N}}{\longrightarrow} 
V_{(\Omega,\mathcal O)} \stackrel{\eta_*}{\longrightarrow}  V_{\Omega_\#}
\stackrel{\phi_{\mathcal N'}^{-1}}{\longrightarrow} \to V_{\Omega_\#'}\to V_{G_\#}$.

Using cut invariance of linear forms $\Phi_{(\Omega,\mathcal O)}$  and the definition
of $\langle X \rangle_{g,\varepsilon}$ we obtain the identity, namely,
tensor $\langle X \rangle_G$ is equal to a contraction of 
tensor $\langle X_\# \rangle_{G_\#}$. 

We did an appropriate choice of admissible sets of local orientations
$\mathcal O$ and $\mathcal O'$ separately for  all 9 classes of simple cuts.
The resulting relations for correlators just form a list of relations in 
the definition of a sysem of correlators. 

\end{proof}

Conversely, assume we are given a system of correlators
$\mathcal C=(A,x\mapsto x^*,B,y\mapsto y^*, \langle X \rangle_G\})$.

First, define a linear form  $\Phi_{(\Omega,\mathcal O)}$ for 
stratified surfaces of type 
\\ $G=(g,\varepsilon,m,m_1,\dots,m_s)$ such that $m_i>0$ for $i=1,\dots,s$ and 
an admissible set of local orientations $\mathcal O$. For this purpose 
use an admissible bijection $\mathcal N:M\to\Omega_0$ and invert 
the above consideration. The function
$\Phi_{(\Omega,\mathcal O)}$ is well defined because the correlators
are symmetric.

Second, using axioms of a Klein topological field theory extend 
the definition of linear forms $\Phi_{(\Omega,\mathcal O)}$ to 
all pairs $(\Omega,\mathcal O)$, i.e., to arbitrary sets of local orientations 
$\mathcal O$ (using axiom $2^\circ$) and to surfaces with boundaries
having no special points (using axiom $4^\circ$).
It can be easily checked that constructed set of data 
$\mathcal T(\mathcal C)=(A,x\mapsto x^*,B,y\mapsto y^*,\{\Phi_{(\Omega,\mathcal O)}\})$
satisfies all axioms of a Klein topological field theory. 
From definitions it follows that the correspondence $\mathcal C 
\rightsquigarrow \mathcal T(\mathcal C)$ is inverse to the correspondence
$\mathcal T \rightsquigarrow \mathcal C(\mathcal T)$. Thus
    
\begin{theorem} \label{th5.1} Correspondence   
$\mathcal T\rightsquigarrow \mathcal C(\mathcal T)$ is a bijection  between 
Klein topological field theories and systems of correlators.
\end{theorem}

%%%%%%%%%%%%%%%%%%%%%%%%%%%%%%%%%%%%%%

\subsection{From system of correlators to structure algebras}

Fix a system of correlators 
$\mathcal C=(A,x\mapsto x^*,B,y\mapsto y^*,\{\langle X\rangle_G\})$.
In this subsection we consider correlators as tensors in  
basis $E_A$ and $E_B$ of vector spaces $A$ and $B$ resp.
Thus, in notation $\langle X\rangle_G$ of a correlator we assume that
$X$ is a set of 
basic elements
$\alpha_i\in E_A$, $\beta_i^j\in E_B$.
An index with star like $\alpha^*$ means the contraction with the matrix 
$I_{\alpha'}^{\alpha''}$ of involutive transformation $x\mapsto x^*$. For example,
$\langle \alpha^*,X \rangle_{(g,\varepsilon)}=
I_{\alpha}^{\alpha'}\langle \alpha',X \rangle_{(g,\varepsilon)}$

We deal here with correlators corresponding to a type $G$ of
not necessarily connected stratified surface. 
It was shown above that any simple cut $\gamma$ 
of a stratified surface $\Omega$ of type
$G$  generates an identity with 
tensor $\langle X \rangle_G$ in the left hand side 
and a contraction of tensor $\langle X_\# \rangle_{G_\#}$ in the 
right hand side. Here $G_\#$ is the type of the contracted cut 
surface $\Omega_\#=\Omega/\gamma$. The explicit formulas for these identities
corresponding to simple cuts of all classes are those that are listed 
in the definition of systems of correlators. 

Identities for correlators that differ only by
symmetries (i)-(ii), we consider as equal. Accepting this agreement we may claim, 
that the identity corresponding to a simple cut 
$\gamma\subset \Omega$ does not depend on the choice of $\Omega$ and 
an admissible set of local orientations $\mathcal O$. 
The identities corresponding to isomorphic pairs $(\Omega,\gamma)$
and $(\Omega',\gamma')$ are equal.

We denote this identity as follows: 
$$\langle X \rangle_G=\conv_{(\Omega,\gamma)}(\langle X_\# \rangle_{G_\#}).$$
Here $X_\#$ is tensor of a type encoded by the type  $G_\#$ of
contracted cut surface $\Omega_\#=\Omega/\gamma$, namely,
if $G_\#=(g_\#,\varepsilon_\#,m_\#,m_{\#1},\dots,m_{\#s_{\#}})$ then
$X_\#$ has $m_\#$ first indices from $A$ and $\sum m_{\#i}$
remaining indices from $B$. Right hand side expression is a contraction of
this tensor by means of tensors $F^{\alpha',\alpha''}$, 
$F^{\beta',\beta''}$, $D^\alpha$ and possibly by their contractions 
with matrices of star involutions $I_{\alpha'}^{\alpha}$,
$I_{\beta'}^{\beta}$. The explicit form of the contraction is 
defined by a pair $(\Omega,\gamma)$ and is encoded in the notation
of operator $\conv_{(\Omega,\gamma)}(\dots)$.  Indices of $X_\#$ that are 
not involved into contraction are just indices of tensor $X$. The type of $X$
is encoded by $G$.

A cut system $\Gamma=\{\gamma_1,\dots,\gamma_k\}$  
generates a chain of identities and finally we obtain the identity as follows:
$$\langle X \rangle_G=\conv_{(\Omega,\gamma_1)}(
 \conv_{(\Omega/\{\gamma_1\},\gamma_2)}(\dots
\conv_{(\Omega/\{\gamma_1,\dots,\gamma_{k-1}\},\gamma_k)}
(\langle X_\# \rangle_{G_\#})\dots)),$$
where $G_\#$ is topological type of contracted cut surface 
$\Omega_\#=\Omega/\Gamma$.                                       
Clearly, the result does not depend on the ordering of simple cuts.
We denote this identity as 
$$\langle X \rangle_G=\conv_{(\Omega,\Gamma)}(\langle X_\# \rangle_{G_\#}).$$

Isomorphic pairs  $(\Omega,\Gamma)$ and $(\Omega',\Gamma')$ 
generate  equal identities.
Note that operator $\conv_{(\Omega,\Gamma)}(\dots)$ 
can be applied to any tensor $\overline{X}_\#$ of type corresponding
to the type  $G_\#$
of the contracted cut surface $\Omega/\Gamma$.

The following identities are evident.

\begin{lemma} \label{5.4} (1) If $\Gamma^\circ\subset\Gamma$
is a subsystem of a cut system $\Gamma$
then $\conv_{(\Omega,\Gamma)}(\overline{X}_\#)=
\conv_{(\Omega,\Gamma^\circ)}
(
\conv_{(\Omega/\Gamma^\circ,\Gamma\setminus\Gamma^\circ)}(\overline{X}_\#) 
).                                             
$
\\ (2) Let $\Gamma'$ be a cut system on $\Omega'$, $\Gamma_1$ and $\Gamma_2$ be
two cut systems on $\Omega^\circ$ and $\overline{X}^\circ_\#$ (resp. $\overline{X}'_\#$)
be a tensor of type corresponding
to the type  $G^\circ_\#$ (resp. $G'_\#$) of the contracted cut surface 
$\Omega'/\Gamma'$ (resp.$\Omega^\circ/\Gamma^\circ$). If
$\conv_{(\Omega^\circ,\Gamma_1)}(\overline{X}_\#)=
\conv_{(\Omega^\circ,\Gamma_2)}(\overline{X}_\#)$ 
then
$\conv_{(\Omega^\circ\sqcup \Omega',\Gamma^\circ\sqcup\Gamma_1)}
(\overline{X}^\circ_\#\cdot \overline{X}'_\#)=
\conv_{(\Omega^\circ\sqcup \Omega',\Gamma^\circ\sqcup\Gamma_2)}
(\overline{X}^\circ_\#\cdot \overline{X}'_\#)$.
\end{lemma}

Correlators corresponding to types $G$ of trivial  surfaces are called 
{\it trivial correlators}. There are five of them; they were used above 
for definitions of bilinear forms $(\alpha',\alpha'')_A$,
$(\beta',\beta'')_B$ and elements $1_A\in A$, $1_B\in B$, $U\in A$.

Correlators corresponding to types $G$ of basic surfaces are called 
{\it basic correlators}. Therefore, there are three basic correlators
$$ \langle x_1,x_2,x_3\rangle_{0,1}, \ 
 \langle (y_1,y_2,y_3)\rangle_{0,1} ,\ 
 \langle x, (y)\rangle_{0,1}.$$

\begin{lemma} \label{5.5} The correlator $\langle X\rangle_G$ for the type $G$ of any 
connected stable stratified surface $\Omega$   is equal to a contraction
of a product of basic correlator.
\end{lemma}

\begin{proof} Indeed, take a stratified surface $\Omega$ of type $G$ and 
a complete cut system $\Gamma$ of $\Omega$. Then we get the identity 
$\langle X \rangle_G=\conv_{(\Omega,\Gamma)}(\langle X_\# \rangle_{G_\#})$.
By definition, connected components of $\Omega_\#$ are basic surfaces.
Therefore,  correlator $\langle X_\# \rangle_{G_\#}$ is a product of basic correlators.
\end{proof}

Let $\Gamma'$ and $\Gamma''$ be two complete cut systems of a stable connected 
stratified surface $\Omega$. Then we have two identities 
$\langle X \rangle_G=\conv_{(\Omega,\Gamma')}(\langle X_\#' \rangle_{G_\#'})$
and $\langle X \rangle_G=\conv_{(\Omega,\Gamma'')}(\langle X_\#'' \rangle_{G_\#''})$.
Hence 
$$
\conv_{(\Omega,\Gamma')}(\langle X_\#' \rangle_{G_\#'})=
\conv_{(\Omega,\Gamma'')}(\langle X_\#'' \rangle_{G_\#''}),
$$
and we obtain the relation between basic correlators.
Obviously, replacing $\Gamma'$ or $\Gamma''$ by an equivalent 
complete cut system brings this relation to the equal one.
Moreover, if two triples $(\Omega_1,\Gamma_1',\Gamma_1'')$ and 
$(\Omega_2,\Gamma_2',\Gamma_2'')$ are equivalent, i.e., there exists 
an isomorphism $\phi:\Omega_1\to\Omega_2$ such that 
$\phi(\Gamma_1')=\Gamma_2'$, $\phi(\Gamma_1'')=\Gamma_2''$ (or
$\phi(\Gamma_1')=\Gamma_2''$, $\phi(\Gamma_1'')=\Gamma_2'$), then
corresponding relations are equal.

We call the relations corresponding to nonequivalent pairs of complete
cut systems on simple surfaces {\it defining relations}. The 
relations corresponding to unstable cuts
of simple and trivial surfaces we  call {\it defining relations} too.

Using the list of simple surfaces and the classification of complete
cut systems on them (see lemma \ref{3.9}) we obtain the following list
of defining relations between
basic correlators. In the list we omit relations that are consequence of 
given ones. The enumeration of relations corresponds to the enumeration 
of simple surfaces. In order to simplify formulas we use nonbasic 
correlators $\langle \alpha,\alpha',\alpha'',\alpha'''\rangle_{0,1}$ and 
$\langle (\beta,\beta',\beta'',\beta''')\rangle_{0,1}$ because  they have clear 
expression via basic correlators.                           

Defining relations corresponding to complete cut systems of
simple surfaces
\medskip
\\ $(1)\ \ 
\langle \alpha_1,\alpha_2,\alpha'\rangle_{0,1}F^{\alpha',\alpha''}
\langle \alpha'',\alpha_3,\alpha_4\rangle_{0,1}=$

$=
\langle \alpha_1,\alpha_3,\alpha'\rangle_{0,1}F^{\alpha',\alpha''}
\langle \alpha'',\alpha_2,\alpha_4\rangle_{0,1}
$
\medskip
\\ $(2)\ \
\langle \alpha,\alpha',\alpha''{^*}\rangle_{0,1}F^{\alpha',\alpha''}=
\langle \alpha,\alpha',\alpha''\rangle_{0,1}D^{\alpha'}D^{\alpha''}
$
\medskip
\medskip
\\ $(4)\ \
\langle (\beta_1,\beta_2,\beta')\rangle_{0,1}F^{\beta',\beta''}
\langle (\beta'',\beta_3,\beta_4)\rangle_{0,1}=$
\\ \hspace{2cm} $=\langle (\beta_2,\beta_3,\beta')\rangle_{0,1}F^{\beta',\beta''}
\langle (\beta'',\beta_4,\beta_1)\rangle_{0,1}
$
\medskip
\\ $(5)\ \
\langle (\beta_1,\beta_2,\beta')\rangle_{0,1}F^{\beta',\beta''}
\langle \alpha,(\beta'')\rangle_{0,1}=
\langle (\beta_2,\beta_1,\beta')\rangle_{0,1}F^{\beta',\beta''}
\langle \alpha,(\beta'')\rangle_{0,1}
$       
\medskip
\\ $(6)\ \
\langle \alpha_1,\alpha_2,\alpha'\rangle_{0,1}F^{\alpha',\alpha''}
\langle \alpha'',(\beta)\rangle_{0,1}=$
\\ \qquad $=\langle (\beta,\beta',\beta'')\rangle_{0,1}
F^{\beta',\beta'''}\langle \alpha_1,(\beta''')\rangle_{0,1}
F^{\beta'',\beta^{(iv)}}\langle \alpha_2,(\beta^{(iv)})\rangle_{0,1}
$
\medskip
\\ $(7)\ \
\langle \alpha',(\beta)\rangle_{0,1}D^{\alpha'}=
\langle (\beta,\beta',\beta''{^*})\rangle_{0,1}F^{\beta',\beta''}
$
\medskip
\\ $(8)\ \ 
\langle \alpha',(\beta_1)\rangle_{0,1}F^{\alpha',\alpha''}
\langle \alpha'',(\beta_2)\rangle_{0,1}=
\langle (\beta_1,\beta',\beta_2,\beta'')\rangle_{0,1}F^{\beta',\beta''}
$
\medskip

Defining relations completely describe properties of 
systems of correlators. Precise statements are formulated and proved 
in sequel subsection using structure algebras.

Denote trivial and basic correlators as follows.

(1)  $F_{\alpha_1,\alpha_2}=\langle \alpha_1,\alpha_2\rangle_{0,1}$

(2)  $F_{\beta_1,\beta_2}=\langle(\beta_1,\beta_2)\rangle_{0,1}$;

(3)  $R_{\alpha,\beta}=\langle \alpha,(\beta)\rangle_{0,1}$;

(4)  $S_{\alpha_1,\alpha_2,\alpha_3}=\langle\alpha_1,\alpha_2,\alpha_3\rangle_{0,1}$; 

(5)  $T_{\beta_1,\beta_2,\beta_3}=\langle(\beta_1,\beta_2,\beta_3)\rangle_{0,1}$;

(6)  $R_{\alpha,\beta_1,\beta_2}=\langle \alpha, (\beta_1,\beta_2)\rangle_{0,1}$;

(7)  $I_{\alpha_1,\alpha_2}=\langle \alpha_1^*,\alpha_2 \rangle_{0,1}$;

(8)  $I_{\beta_1,\beta_2}=\langle(\beta_1^*,\beta_2)\rangle_{0,1}$;

(9)  $D_{\alpha}=\langle\alpha \rangle_{\frac{1}{2},0}$;

(10) $J_\alpha=\langle\alpha \rangle_{0,1}$;

(11) $J_\beta=\langle (\beta)\rangle_{0,1}$.

\begin{theorem} \label{th5.2} Tensors (1)-(11) are structure constants of 
a structure algebra, which we  denote by $\mathcal H(\mathcal C)$. 
\end{theorem}

\begin{proof} By lemma \ref{2.2} it is sufficient to verify relations
(1)--(16) for these tensors.
These relations either follow from the symmetries of correlators or 
coincide with defining relations for basic correlators (1)-(8).

\end{proof}

%%%%%%%%%%%%%%%

\subsection{From structure algebra to Klein topological field theory}

Let $\mathcal H=(H=A\oplus B, (.,.), x\mapsto x^*, U\in A)$ 
be a structure algebra. We shall construct a Klein topological
field theory 
$\mathcal T(\mathcal H)=(A,x\mapsto x^*,B,y\mapsto y^*,\{\Phi_{(\Omega,\mathcal O)}\})$
with the same vector spaces $A$ and $B$ and involutive transformations
$x\mapsto x^*$ and $y\mapsto y^*$ which are restrictions of involutive 
antiautomorphism of algebra $H$ to its summands.

First, define linear forms
$\Phi_{(\Omega,\mathcal O)}:V_{\Omega}\to\mathbb C$
for trivial and  basic surfaces $\Omega$. Choose basis $E_A$ of
vector space $A$ and $E_B$ of vector space $B$. It is sufficient to
define a linear form on basic elements 
$(\otimes_{p\in\Omega_a}\alpha_p )\otimes(\otimes_{q\in\Omega_b}\beta_q)$
of vector space $V_{\Omega}$. 
Here $\Omega_a$ (resp. $\Omega_b$) is the set of all interior (resp. boundary)
special points of $\Omega$, $\alpha_p\in E_A$, $\beta_q\in E_B$.

For example, let $\Omega$ be isomorphic to a disc with three boundary
special points $q_1,q_2,q_3$.
Choose  an admissible set $\mathcal O$ of local orientations and fix
an admissible ordering $\mathcal N$ of special points. 
In this case $\mathcal O$ and $\mathcal N$ are uniquely defined by fixing
an orientation of the disc and choosing a first special point on its
boundary. The order of special points 
defines the isomorphism $V_{\Omega}\eqsim B_1\otimes B_2\otimes B_3$,
where $B_i$ is a copy of vector space $B$ and $i$ means that the isomorphism 
brings a vector $\beta_i$ assigned to $i$-th special point $q_i$ to this 
multiplier  ($i=1,2,3$).
For basic element  $\otimes_{i=1}^3 \beta_i\in V_{\Omega}$ put 
$\Phi_{(\Omega,\mathcal O)}(\otimes_{i=1}^3 \beta_i)
=T_{\beta_1,\beta_2,\beta_3}$, where tensor 
$T_{\beta_1,\beta_2,\beta_3}$ belongs to structure constants of 
structure algebra $\mathcal H$ (see subsection 2.1).
For a nonadmissible set $\mathcal O$ of local orientations
define $\Phi_{(\Omega,\mathcal O)}$ by invariance of 
changing of local orientations (axiom $2^\circ$ from the definition
of Klein topological field theory).
Properties of structure constants of a structure algebra
(see lemma \ref{2.2}) provide the correctness of
the definition of $\Phi_{(\Omega,\mathcal O)}$,
its topological invariance and invariance with respect
to changing of local orientations.

Analogous reasoning leads to correct definitions of 
$\Phi_{(\Omega,\mathcal O)}$ for all 
trivial and basic surfaces. Key formulas for a surface $\Omega$, 
an admissible set $\mathcal O$ of local orientations and 
an admissible numbering $\mathcal N$ of  special points  are as follows.
Right-hand sides of these formulas are expressed via structure constants of
$\mathcal H$.

{\it Trivial surfaces}

\begin{itemize}
\item[(1)] $\Omega$ is a sphere $S^2$ without special points:
$V_\Omega=\mathbb C$, 
$\Phi_{(\Omega,\mathcal O)}(1)=J^{\alpha}J_{\alpha}$.

\item[(2)]  $\Omega$ is a projective plane $\mathbb{R}P^2$ without special points:
$V_\Omega=\mathbb C$, 
$\Phi_{(\Omega,\mathcal O)}(1)=J_{\alpha}D^{\alpha}$.

\item[(3)]  $\Omega$ is a disc  $D^2$ without special points:
$V_\Omega=\mathbb C$, 
$\Phi_{(\Omega,\mathcal O)}(1)=J^{\beta}J_{\beta}$.
 
\item[(4)]  $\Omega$ is a sphere $S^2$ with an interior special point  $p$:
$V_\Omega=A$,
$\Phi_{(\Omega,\mathcal O)}(\alpha)=J_\alpha$.

\item[(5)]  $\Omega$ is a disc $D^2$ with a boundary special point $q$:
$V_\Omega=B$,
$\Phi_{(\Omega,\mathcal O)}(\beta)=J_\beta$.
 
\item[(6)]  $\Omega$ is a sphere $S^2$ with two interior special
points $p_1$, $p_2$:
$V_\Omega\eqsim A\otimes A$,
$\Phi_{(\Omega,\mathcal O)}(\alpha_1\otimes\alpha_2)=F_{\alpha_1,\alpha_2}$.

\item[(7)] $\Omega$ is a projective plane $\mathbb{R}P^2$ with an interior
special point $p$:
$V_\Omega=A$, 
$\Phi_{(\Omega,\mathcal O)}(\alpha)=D_{\alpha}$.

\item[(8)]  $\Omega$ is a torus $T^2$ without special points: 
$V_\Omega=\mathbb C$, 
$\Phi_{(\Omega,\mathcal O)}(1)=
F_{\alpha',\alpha''}F^{\alpha',\alpha''}=\dim A$.
 
\item[(9)]  $\Omega$ is a Klein bottle $\Kl$ without special points: 
$V_\Omega=\mathbb C$, 
$\Phi_{(\Omega,\mathcal O)}(1)=
I_{\alpha',\alpha''}F^{\alpha',\alpha''}=\tr(x\mapsto x^*)$.

\item[(10)] $\Omega$ is a disc  $D^2$ with an interior special point $p$:
$V_\Omega=A$, 
$\Phi_{(\Omega,\mathcal O)}(\alpha)=J^{\beta}R_{\alpha,\beta}$.

\item[(11)] $\Omega$ is a disc $D^2$ with two boundary special
points $q_1$, $q_2$:
$V_\Omega\eqsim B\otimes B$, 
$\Phi_{(\Omega,\mathcal O)}(\beta_1\otimes\beta_2)=F_{\beta_1,\beta_2}$.

\item[(12)] $\Omega$ is a M\"obius band $\Mb$ without special points: 
$V_\Omega=\mathbb C$, 
$\Phi_{(\Omega,\mathcal O)}(1)=
J^{\beta}D^{\alpha}R_{\alpha,\beta}$.

\item[(13)] $\Omega$ is a cylinder $\Cyl$ without special points:
$V_\Omega=\mathbb C$, 
\\ $\Phi_{(\Omega,\mathcal O)}(1)=
J^{\beta'}R_{\alpha',\beta'}
F^{\alpha',\alpha''}J^{\beta''}R_{\alpha'',\beta''}$.
\end{itemize}

{\it Basic surfaces}

\begin{itemize}

\item[(1)] $\Omega$ is a sphere $S^2$ with three interior special
points $p_1,p_2,p_3$:
$V_\Omega\eqsim A\otimes A\otimes A$, 
$\Phi_{(\Omega,\mathcal O)}(\alpha_1\otimes\alpha_2\otimes\alpha_3)=
S_{\alpha_1,\alpha_2,\alpha_3}$.

\item[(2)] $\Omega$ is a disc $D^2$ with three boundary special
points $q_1,q_2,q_3$:
$V_\Omega\eqsim B\otimes B\otimes B$, 
$\Phi_{(\Omega,\mathcal O)}(\beta_1\otimes\beta_2\otimes\beta_3)=
T_{\beta_1,\beta_2,\beta_3}$.

\item[(3)] $\Omega$ is a disc $D^2$ with a boundary special
point $q$ and an interior special point $p$:
$V_\Omega= A\otimes B$, 
$\Phi_{(\Omega,\mathcal O)}(\alpha\otimes\beta)=
R_{\alpha,\beta}$.

\end{itemize}

Multiplicativity (axiom  $6^\circ$ from the definition of Klein
topological field theory) allows us to define
correctly linear forms $\Phi_{(\Omega,\mathcal O)}$ also for 
stratified surfaces $\Omega$ that are disjoint union of several trivial and 
basic surfaces. 

\begin{lemma} \label{5.6} The already defined linear forms
$\Phi_{(\Omega,\mathcal O)}$ for trivial and basic surfaces and
their disjoint unions satisfy topological invariance axiom, invariance
of a change of local orientations axiom, nondegeneracy axiom,
invariance of adding trivial field axiom (when applicable), cut
invariance axiom (for cuts of trivial and basic surfaces) and
multiplicativity axiom.
\end{lemma}

\begin{proof}
Topological invariance and  invariance of a change of local orientations
follows from the definitions of functions $\Phi_{(\Omega,\mathcal O)}$ as
it was shown above for a disc with three boundary special points. 
Nondegeneracy follows from nondegeneracy of bilinear forms 
$(x',x'')|_A$ and $(y',y'')|_B$ in algebra $\mathcal H$.
Invariance of adding a trivial field and 
cut invariance for  trivial and basic surfaces are checked for all such
surfaces. Any verification turns into a simple identity for 
structure constants of $\mathcal H$. We skip the details.
Multiplicativity follows from the definition of 
$\Phi_{(\Omega,\mathcal O)}$ for 
disjoint unions of  trivial and  basic surfaces. 
\end{proof}

Define $\Phi_{(\Omega,\mathcal O)}$ for a connected 
stable surface $\Omega$. 
Choose a complete cut system  $\Gamma$  of $\Omega$ and orientations
of all simple cuts $\gamma\subset\Gamma$. Denote by
$\eta:(\Omega,\mathcal O)\to(\Omega_\#,\mathcal O_\#)$
cut morphism corresponding to $\Gamma$ (see subsection 4.1).
By definition of a complete cut system, 
all connected components of contracted cut 
surface $\Omega_\#$ are basic surfaces. 
Therefore, function 
$\Phi_{(\Omega_\#,\mathcal O_\#)}$ is already defined.
Put $\Phi_{(\Omega,\mathcal O,\Gamma)}(x)=
\Phi_{(\Omega_\#,\mathcal O_\#)}(\eta_*(x))$. 
We have to prove that 
this function does not depend on the choice of $\Gamma$. 

\begin{lemma} \label{5.7} Let $\Gamma$ and $\Gamma'$ be equivalent complete cut systems
endowed with orientations of all simple cuts of a stable surface $\Omega$. 
Then
$\Phi_{(\Omega,\mathcal O,\Gamma)}(x)=\Phi_{(\Omega,\mathcal O,\Gamma')}(x)$.
Particularly, $\Phi_{(\Omega,\mathcal O,\Gamma)}(x)$ does not depend 
on a choice of orientations of simple cuts.
\end{lemma}

\begin{proof} Suppose, we change an orientation of a simple cut 
$\gamma\subset\Gamma$ of a complete cut system $\Gamma$.
Denote the same complete cut system with
changed orientation of the cut by $\Gamma'$. We have two 
cut morphisms $\eta_{\Gamma}:(\Omega,\mathcal O)\to (\Omega_\#,\mathcal O_\#)$ and 
$\eta_{\Gamma'}:(\Omega,\mathcal O) \to (\Omega_\#,\mathcal O'_\#)$ and 
changing local orientations morphism 
$\psi:(\Omega_\#,\mathcal O_\#)\to(\Omega_\#,\mathcal O_\#')$.
Vector space $V_{\Omega_\#}$ is isomorphic to 
$V_\Omega\otimes X$, where vector space $X$ depends on $\Gamma$. 
By definition, linear map $\eta_*:V_\Omega\to V_{\Omega_\#}$
brings element $x\in V_\Omega$ to $x\otimes z$, where $z$ is a 
tensor product of several elements $K_{A,*}\in A$, 
$K_{B,*}\in B$ and $U\in A$ (see subsection 2.1). It 
can be easily checked that these elements of structure 
algebra are invariants of involution $x\mapsto x^*$. Therefore,
$(\eta_{\Gamma'})_*(x)=(\eta_{\Gamma})_*(\psi_*(x))=(\eta_{\Gamma})_*(x)$
and 
$\Phi_{(\Omega,\mathcal O,\Gamma)}(x)=\Phi_{(\Omega,\mathcal O,\Gamma')}(x)$.

The equality  $\Phi_{(\Omega,\mathcal O,\Gamma)}=\Phi_{(\Omega,\mathcal O,\Gamma')}$
for different but equivalent complete cut systems $\Gamma$ and $\Gamma'$
follows from the topological invariance for trivial and basic surfaces
and their disjoint unions applied to 
linear forms $\Phi_{(\Omega/\Gamma,\mathcal O_\#)}$ and  
$\Phi_{(\Omega/\Gamma',\mathcal O_\#')}$.  
\end{proof}

\begin{lemma} \label{5.8} If $\Omega$ is a simple stratified surface then
$\Phi_{(\Omega,\mathcal O,\Gamma)}=\Phi_{(\Omega,\mathcal O,\Gamma')}$
for any pair $\Gamma$, $\Gamma'$  of complete cut systems of $\Omega$.
\end{lemma}

\begin{proof} 
Consider separately simple surfaces (1) - (8) (see lemma \ref{3.9}).
For example, let $\Omega$ be a disc with two interior special points $p_1$, $p_2$
and one boundary special point $q$.
By lemma \ref{5.7}, it is sufficient to  consider only nonequivalent
complete cut systems $\Gamma$ and $\Gamma'$.
There are two nonequivalent complete cut 
systems $\Gamma$ and $\Gamma'$ of $\Omega$. One of them, $\Gamma$, 
consists of one separating contour $\gamma$ and $\Omega/\Gamma$ is isomorphic 
to the disjoint union of a sphere with three interior special points
$p_1,p_2,p'$ and  a disc with one interior special point $p''$  and 
one boundary special point $q$. Another one, $\Gamma'$, consists of two 
separating segments cut $\gamma_1$, $\gamma_2$ and 
$\Omega/\Gamma'$ is isomorphic to the disjoint union of three components: 

$\circ$ a disc with
interior special point $p_1$ and boundary special point $q_1'$,

$\circ$ a disc with 
interior special point $p_2$ and boundary special point $q_2'$,

$\circ$ a disc with 
three boundary special points $q,q_1'', q_2''$.

From definitions of $\Phi_{(\Omega,\mathcal O,\Gamma)}$ and 
$\Phi_{(\Omega,\mathcal O,\Gamma')}$ we obtain that the equality
\\ $\Phi_{(\Omega,\mathcal O,\Gamma)}=\Phi_{(\Omega,\mathcal O,\Gamma')}$
is equivalent to the equality 
$$R_{\alpha_1,\beta_1'}F^{\beta_1',\beta_1''} 
R_{\alpha_2,\beta_2'}F^{\beta_2',\beta_2''} 
T_{\beta'\beta''\beta}=
S_{\alpha_1\alpha_2\alpha'}F^{\alpha',\alpha''}R_{\alpha'',\beta}$$
for structure constants.
The latter coincides with  relation (5) for 
structure constants (see lemma \ref{2.2}) of structure algebra $\mathcal H$.

Analogous reasoning are applicable in all cases.
Equalities for structure constants derived from equalities  
$\Phi_{(\Omega,\mathcal O,\Gamma)}=\Phi_{(\Omega,\mathcal O,\Gamma')}$
for complete cut systems of simple surfaces are just 
relations for structure constants of a structure algebra. 
We omit the details.
\end{proof}

It follows from the lemma that formula 
$\Phi_{(\Omega,\mathcal O)}=\Phi_{(\Omega,\mathcal O,\Gamma)}$ correctly
defines a function $\Phi_{(\Omega,\mathcal O)}$ for all simple surfaces.
Now, we can define by  multiplicativity  linear form 
$\Phi_{(\Omega,\mathcal O)}$ for 
any  surface that is the disjoint union of trivial, basic and simple 
surfaces.
Topological invariance for these forms 
and invariance of change of local orientations can be easily checked.

\begin{lemma} \label{5.9} Let $\Gamma$ and $\Gamma'$ be adjacent
complete cut systems of a stable stratified surface 
$\Omega$. Then  
$\Phi_{(\Omega,\mathcal O,\Gamma)}=
\Phi_{(\Omega,\mathcal O,\Gamma')}$
\end{lemma}

\begin{proof} Let  $\Gamma$ and $\Gamma'$  be two  nonequivalent
adjacent complete cut systems of $\Omega$.
By lemma \ref{5.7}, the replacement of $\Gamma$ by equivalent cut system leads to the same 
linear form. Hence, we may assume 
that the set $\Gamma^\circ$ of cuts common to 
$\Gamma$ and $\Gamma'$ has the following property.
Contracted cut surface $\Omega/\Gamma^\circ$ is the 
disjoint union of  a simple surface $\Omega^\circ$ and several
basic surfaces and cut systems 
$\overline{\Gamma}=\Gamma\setminus\Gamma^\circ$, 
$\overline{\Gamma'}=\Gamma'\setminus\Gamma^\circ$ belong 
to $\Omega^\circ$; therefore $\overline{\Gamma}$ and
$\overline{\Gamma'}$ are complete cut systems 
of simple surface $\Omega^\circ$. 

We have that linear forms $\Phi_{(\Omega/\Gamma^\circ,\mathcal O^\circ_\#)}$
are already correctly defined for any set $\mathcal O^\circ_\#$ of local 
orientations. From {\it transitivity property} of cut morphisms 
we obtain that 
$\Phi_{(\Omega,\mathcal O,\Gamma)}(x)=
\Phi_{(\Omega/\Gamma^\circ,\mathcal O^\circ_\#)}(\eta^\circ(x))$
and $\Phi_{(\Omega,\mathcal O,\Gamma')}(x)=
\Phi_{(\Omega/\Gamma^\circ,\mathcal O^\circ_\#)}(\eta^\circ(x))$,
where $\eta^\circ:V_\Omega\to V_{\Omega/\Gamma^\circ}$ is linear map
induced by cut morphism 
$\eta^\circ:(\Omega,\mathcal O)\to(\Omega/\Gamma^\circ,\mathcal O^\circ_\#)$.
Therefore, $\Phi_{(\Omega,\mathcal O,\Gamma)}=
\Phi_{(\Omega,\mathcal O,\Gamma')}$.
\end{proof}

By theorem \ref{th4.1} any pair of complete cut systems of $\Omega$ is
connected by a chain of adjacent complete cut systems. Thus, formula 
$\Phi_{(\Omega,\mathcal O)}=\Phi_{(\Omega,\mathcal O,\Gamma)}$ correctly
defines function $\Phi_{(\Omega,\mathcal O)}$
and we define linear forms $\Phi_{(\Omega,\mathcal O)}$ for all 
connected stratified surfaces. 
Define $\Phi_{(\Omega,\mathcal O})$ for disconnected surfaces by multiplicativity.

\begin{theorem} \label{th5.3} A tuple  of data
$\mathcal T(\mathcal H)=
(A,x\mapsto x^*,B,y\mapsto y^*,\{\Phi_{(\Omega,\mathcal O)}\})$
is a Klein topological field theory.
\end{theorem}

\begin{proof}

Let us verify axioms of a Klein topological field theory.

$1^\circ$ It is sufficient to check topological invariance for 
connected stable surfaces that are not trivial or basic surfaces. 
Let $\phi:(\Omega,\mathcal O)\to (\Omega',\mathcal O')$ be an isomorphism.
Choose a complete cut system $\Gamma$ of $\Omega$ and denote
by $\Gamma'$ complete cut system $\phi(\Gamma)$ of $\Omega'$. 
The isomorphism $\phi$ induces the isomorphism 
$\phi_\#:(\Omega_\#,\mathcal O_\#)\to (\Omega_\#',\mathcal O_\#')$
of contracted cut surfaces 
$\Omega_\#=\Omega/\Gamma$ and $\Omega_\#'=\Omega'/\Gamma'$.
The topological invariance with respect of $\phi$ follows from
the topological invariance of $\phi_\#$. The latter was
proven in lemma \ref{5.6}. 

$2^\circ$ Invariance of changing of local orientations
is verified analogously to  $1^\circ$.

$3^\circ$ Nondegeneracy follows from  nondegeneracy of bilinear forms
$(x',x'')|_A$ and $(y',y'')_B$.

$4^\circ$ Invariance of adding trivial field follows from the same
statement proved for basic and trivial surfaces (lemma \ref{5.6}).
Indeed, let $\Omega$ be a stable surface and $r$ be a nonspecial point.
Choose complete cut system $\Gamma$ that does not
meet $r$. Then the image of $r$ in contracted cut surface belongs to 
a connected component $\Omega_{\#i}$, which is basic surface.
Invariance of adding new special point $r$ follows from the same
statement for its image in $\Omega_{\#i}$.

$5^\circ$ It is sufficient to check cut invariance for 
a simple cut $\gamma$ 
of a connected stable surface $\Omega$ that is not trivial or basic surface. 

Suppose that $\gamma$ is a stable cut. Then it can be included in a complete 
cut system $\Gamma$ of $\Omega$. By properties of cut morphisms
$(\eta_{\Gamma})_*= (\eta_{\Gamma\setminus\gamma})_*\circ(\eta_{\gamma})_*$,
where $(\eta_{\Gamma})_*:V_\Omega\to V_{\Omega/\Gamma}$,
$(\eta_{\Gamma\setminus\gamma})_*:V_{\Omega/\gamma}\to V_{\Omega/\Gamma}$,
$(\eta_{\gamma})_*:V_\Omega\to V_{\Omega/\gamma}$ are linear maps
corresponding to cut morphisms. Therefore,
$\Phi_{(\Omega,\mathcal O)}(x)=\Phi_{(\Omega/\gamma,\mathcal O_\#)}
(\eta_{\gamma})_*(x))$.

Suppose that $\gamma$ is trivial (unstable) cut. Then there exists
a complete cut
system $\Gamma$ such that $\Gamma\cap\gamma=\emptyset$. Hence, 
the image of $\gamma$ in contracted cut surface is a simple cut of
a connected
component  $\Omega_{\#i}$, which is basic surface. The invariance with respect 
to cut ${\gamma}$ follows from the invariance with respect of its image in
the basic surface, which was proved in lemma \ref{5.6}.

$6^\circ$ Multiplicativity is clear because linear forms
$\Phi_{(\Omega,\mathcal O)}$
for disconnected surfaces are defined using this axiom.

\end{proof}

\begin{theorem} \label{MMain} {\bf (Main)} The correspondences
$\mathcal T\rightsquigarrow \mathcal C(\mathcal T)
\rightsquigarrow\mathcal H
\mathcal C(\mathcal T))$ and $\mathcal H\rightsquigarrow \mathcal T
(\mathcal H)$ are reciprocal bijections between
isomorphism classes of Klein topological field theories 
$\mathcal T=\{A, x\mapsto x^*, B, y\mapsto y^*, \Phi_{(\Omega, \mathcal O)}\}$ 
and isomorphism classes of structure algebras
$\mathcal H=\{H=A\oplus B, (.,.), x\mapsto x^*, U\}$. The correlators for
Klein topological theory $\mathcal T(\mathcal H)$ have the following 
expressions in terms of structure algebra $\mathcal H$:
$$\langle x_1,\dots,x_m,({y_1}^1,\dots,y_{m_1}^1),\dots,(y_1^s,\dots,
y_{m_s}^s)\rangle_{g,1}=$$
$$=(x_1\dots x_m ( y_1^1\dots y_{m_1}^1)V_{K_B}( y_1^2\dots y_{m_2}^2)\dots V_{K_B}(y_1^s\dots
y_{m_s}^s),K_A^{g})$$
and
$$\langle x_1,\dots,x_m,(y_1^1,\dots,y_{m_1}^1),\dots,(y_1^s,\dots,y_{m_s}^s)
\rangle_{g,0}=$$
$$=(x_1\dots x_m ( y_1^1\dots y_{m_1}^1)V_{K_B}( y_1^2\dots y_{m_2}^2)\dots V_{K_B}(y_1^s\dots
y_{m_s}^s),U^{2g})$$
\end{theorem}

\begin{proof} The first statement follows from the exact constructions
of the correspondences 
$\mathcal T\rightsquigarrow \mathcal C(\mathcal T)$,
$\mathcal C\rightsquigarrow \mathcal H(\mathcal C)$,
$\mathcal H\rightsquigarrow \mathcal T(\mathcal H)$ and
from theorems \ref{th5.1}, \ref{th5.2}, \ref{th5.3}. The second
statement follows from a representation
$\langle X \rangle_{G}=
\conv_{(\Omega,\Gamma)}(\langle X_\# \rangle_{G_\#})$, where
$\Gamma$ is a specially chosen complete cut system. Namely, $\Gamma$ contains
a system of cuts $\Gamma^\circ$ between holes 1st and 2nd, 2nd and 3d, etc.
and $\Gamma\setminus\Gamma^\circ$ is a normal cut system of 
$\Omega/\Gamma^\circ$.
\end{proof}

%%%%%%%%%%%%%%%%%%%%%%%

\subsection{Open-closed topological field theory}
\label{subsec45}

Let $\mathcal H$ be a Klein topological field theory. One can restrict the basic category of surfaces to the category of oriented surfaces with 
local orientations at all special points induced by a global orientation 
of the surface. Call it {\it open-closed topological field theory} because it 
is equivalent to the open-closed topological field theory as defined 
in  \cite{Laz}, \cite{Moo2}.

It is convenient to discuss the correspondence between versions of open-closed TFTs in
terms of algebras. Note that both the involutive antiautomorphism $x\mapsto x^*$,
responsible for changes of local orientations, and  the element $U$ responsible 
for gluing the hole by a M\"obius band, are not used in the oriented case.

Define   a Lazaroiu-Moore algebra as "the structure algebra without an involutive antiautomorphism $x\mapsto x^*$, 
an element $U$, and axioms $5^\circ$, $6^\circ$,$7^\circ$,$8^\circ$" (compare with  \cite{Laz}, \cite{Moo2}).
That is

\begin{definition} A Lazaroiu-Moore  algebra 
$\mathcal H=\{H=A\oplus B, (.,.)\}$ 
is a finite dimensional associative algebra  $H$  
endowed with

a decomposition $H=A\oplus B$ of $H$ into the direct sum of two vector spaces; 

a symmetric invariant scalar product  $(.,.): H\otimes H\to\mathbb C$,
\\ such that the following axioms hold:

$1^\circ$ $A$ is a subalgebra belonging to the center of algebra $H$; 
algebra $A$ has an unit $1_A\in A$, and $1_A$ is also the unit of algebra $H$;

$2^\circ$ $B$ is a two-sided ideal of $H$ (typically noncommutative); algebra 
$B$ has a unit $1_B\in B$;

$3^\circ$ restrictions $(.,.)|_A$ and $(.,.)|_B$ are 
nondegenerate scalar products on algebras $A$ and $B$ resp. 

$4^\circ$  $(V_{K_B}(b_1),b_2)=(\widehat{K}_A,b_1\otimes b_2)$
for arbitrary $b_1,b_2\in B$ 
(this axiom reflects Cardy  relation \cite{Laz});

\end{definition} 

Reduplicating the proof ot theorem \ref{MMain} we prove

\begin{theorem} \label{op-cl} The correspondences
$\mathcal T\rightsquigarrow \mathcal C(\mathcal T)
\rightsquigarrow\mathcal H
\mathcal C(\mathcal T))$ and $\mathcal H\rightsquigarrow \mathcal T
(\mathcal H)$ are reciprocal bijections between
isomorphism classes of open-closed topological field theories 
$\mathcal T=\{A, B,\Phi_{(\Omega, \mathcal O)}\}$ 
and isomorphism classes of  Lazaroiu-Moore algebras
$\mathcal H=\{H=A\oplus B, (.,.)\}$. The correlators for
open-closed topological theory $\mathcal T(\mathcal H)$ have the following 
expressions in terms of open-closed string algebra $\mathcal H$:
$$\langle x_1,\dots,x_m,({y_1}^1,\dots,y_{m_1}^1),\dots,(y_1^s,\dots,
y_{m_s}^s)\rangle_{g,1}=$$
$$=(x_1\dots x_m ( y_1^1\dots y_{m_1}^1)V_{K_B}( y_1^2\dots y_{m_2}^2)\dots V_{K_B}(y_1^s\dots
y_{m_s}^s),K_A^{g})$$

\end{theorem}

The following claim is clear from  Main theorem, theorem \ref{op-cl}, and results of section
\ref{sec1}.

\begin{corollary} \label{cor4.1} Any massive (corresponding to  a semisimple
algebra) open-closed topological field
theory can be extended to a Klein topological field theory, and the number of such extentions in finite. Equivalent classes 
of extensions are in a bijection with pairs consisting of 
an involutive antiautomorphism $x\mapsto x^*$ and an 
element $U$ satisfying axioms of a structure algebra. 
\end{corollary}

Number of extensions of an open-closed TFT to a KTFT can be easely computed
by results of subsection \ref{2s2}

%%%%%%%%%%%%%%%%%%%%%%%%%%%%%%

\section{Hurwitz topological field theory}
\label{sec5}

In this section the coverings over a stratified surface are considered. It is shown (subsection 5.1) that 
singularities over boundary special points are classified by 'dihedral Yang diagrams', which in turn correspond to conjugacy classes of pairs of involutions in a symmetric
group $S_n$. 

Classical Hurwits numbers are generalized to coverings over 
stratified surfaces. A Klien topological field theory is associated (subsections 5.2, 5.3)
with these generalized Hurwitz numbers. This KTFT is called a {\it Hurwits topological field 
theory}. 

It is proved that a Hurwits topological field theory corresponds to the structure
algebra associated with a symmetric group $S_n$. It allows us to obtain
the expressions (subsection 5.4) for generalized Hurwits numbers via structural constants of the structural
algebra.

\subsection{Stratified coverings}
\label{subsecSC}

Fix a stratified topological space $\Omega=\coprod_{\lambda\in\Lambda}
\Omega_\lambda$. Let $\pi:P\to\Omega$ be a continuous epimorphic map of a
topological space $P$ onto $\Omega$. 
Denote by $P_\lambda$ the preimage of a stratum $\Omega_\lambda$ and by 
$\pi_\lambda:P_\lambda\to\Omega_\lambda$ the restriction of $\pi$ to
$P_\lambda$. Analogously, for an arbitrary subset $X\in\Omega$  denote by 
$\pi_X:P_X\to X$ the restriction of $\pi$ to $P_X=\pi^{-1}(X)$.

\begin{definition} 
\label{defcov}
A continuous map $\pi:P\to\Omega$ 
is called a covering over the stratification (or a stratified covering)
if and only if
for any stratum $\lambda$ the restriction $\pi_\lambda:P_\lambda\to
\Omega_\lambda$ is a local homeomorphism, i.e., a covering.        

A stratified covering $\pi:P\to\Omega$  is called special coverings if and only if
\\ (1) $\Omega$ is a specially stratified manifold (i.e., any codimension 1
stratum belongs to the boundary of $\Omega$);
\\ (2) $P$ is a manifold, possibly with a boundary.
\end{definition}              

In this work we deal with specially stratified coverings   over 
surfaces and call them {\it coverings} for short.

A covering is called {\it $n$-sheeted covering} if and only if the full preimage of any point of
a stratum of maximal dimension consists of $n$ points. $n$-sheeted coverings 
over disconnected base can (and sometimes, will) be considered.
In general, it is not assumed that a cover space $P$ is connected.

An {\it isomorphism}  of two coverings
$\pi_1:P_1\to\Omega$ and $\pi_2:P_2\to\Omega$ over the same base is a 
homeomorphism  $\phi:P_1\to P_2$ such that $\pi_2\circ\phi=\pi_1$.
By ${\mathcal Cov_n}(\Omega)$ we denote the category of $n$-sheeted 
stratified coverings over $\Omega$ with morphisms being isomorphisms of coverings.
By ${\rm Cov_n}(\Omega)$ is denoted the set of 
isomorphism classes in category  ${\mathcal Cov_n}(\Omega)$.
Clearly,  ${\rm Cov_n}(\Omega)$ is a finite set.

\begin{lemma} An isomorphism $\phi:\Omega\to\Omega'$ of stratified
manifolds induces the bijection between sets
${\rm Cov_n}(\Omega)$ and ${\rm Cov_n}(\Omega')$.
\end{lemma}

Indeed, the bijection is induced by morphism 
$\pi\mapsto\phi\circ\pi$. 

\begin{definition} Let $\pi:P\to\Omega$ and $\pi':P'\to\Omega'$ be 
two coverings over stratified manifolds $\Omega$ and $\Omega'$. 
A pair of homeomorphisms $\hat\phi:P\to P'$ and $\phi:\Omega\to\Omega'$
is called a topological equivalence of 
coverings if and only if $\pi'\circ\hat\phi=\phi\circ\pi$.
A class of topologically equivalent coverings is called a topological type
of coverings.

\end{definition}

Obviously, two isomorphic coverings over a stratified manifold 
$\Omega$ are of the same topological type. The converse is not true:
there exist nonisomorphic coverings over 
a surface that are  of the same topological type.

Let $\Aut\Omega$ be the group of all homeomorphisms of $\Omega$, that
preserve the stratification of $\Omega$ (a permutation of strata is allowed).
The group $\Aut\Omega$  acts on the set 
${\rm Cov}(\Omega)$: if $\phi\in\Aut\Omega$ and
$\pi\in{\rm Cov}(\Omega)$ then
$\phi(\pi)=\phi\circ\pi$. Clearly, two coverings
$\pi,\pi'\in {\rm Cov}(\Omega)$ are of the same topological
type if and only if they lie
in the same orbit of the group $\Aut(\Omega)$.

Note, that one can correctly assign a topological type
of a covering $\pi:P\to\Omega$ to a topological type of $\Omega$. 
In contrast, the isomorphism class $[\pi]$
of covering $\pi$ can be associated with just manifold $\Omega$ 
but not with its topological type because one cannot choose an isomorphism
between two different (but isomorphic) bases $\Omega$ and $\Omega'$
canonically.

For an oriented base $\Omega$ we define an {\it oriented
topological type} of a covering with respect to homeomorphisms preserving the
orientation.

There exist only two, up to isomorphism, 
compact connected special stratified one-dimensional manifolds,
a circle $S^1$ and a closed segment $I$. It is convenient 
to consider also one noncompact manifold, namely,
a ray $R$ (with two obvious strata).

%%%%%%%%%%%%%%%%%%%%%%%%%%%%%%%%%%%%%%%%%%%%%%%%%%%%%%%%%%%%%%%%%%%

	\subsubsection{Coverings over a circle}
                                                       
A covering $\pi:P\to S^1$ over a circle has no special points, i.e., it is a
nonramified covering. Isomorphism classes of $n$-sheeted coverings over $S^1$
are in one-to-one correspondence with unordered partitions of the number 
$n$ into a sum $n=n_1+\dots+n_s$. If $\pi:P\to S^1$ is a covering then $s$ is
the number of connected components of $P$ and $n_i$ is the degree of the
restriction of $\pi$ to $i$-th connected component of $P$. 

Equivalently, the isomorphism class of a covering over $S^1$ can be presented
1) as a Young diagram of order n; 2) as a conjugacy class of a symmetric group
$S_n$. We denote by $\alpha$ a partition of $n$, as well, as corresponding
Young diagram of order $n$ and a conjugacy class of $S_n$. Group $\Aut S^1$
acts trivially on the set ${\rm Cov_n}(S^1)$ of isomorphism classes of coverings.
Therefore,  topological type of a covering can be also presented by the
partition $\alpha$.

%%%%%%%%%%%%%%%%%%%%%%%%%%%%%%%%%%%%%%%%%%%%%%%%%%%%%%%%%%%%%%%%%%%
	\subsubsection{Coverings over a ray}

Let $\pi:P\to R$ be an $n$-sheeted covering over a ray. Denote by $M$ the set
of connected components of the preimage of an open interval
$R^\circ=R-\partial R$.
Obviously, $|M|=n$ and the restriction of $\pi$ to any component is a
homeomorphism. Any point $p'$ of the preimage of the end-point $p\in R$ 
belong to the closure of either just one or just two components of $M$.
Therefore, a special covering defines an involutive element $s$ of 
group  $\Aut M$ of all bijections $M\to M$. Group $\Aut M$ is isomorphic
to symmetric group $S_n$; there is no canonical isomorphism, therefore, 
involution $s\in M$ correctly defines only conjugasy class of its image in $S_n$.

Evidently, isomorphism classes of coverings over the ray $R$  are in one-to-one
correspondence with conjugacy classes of involutions in $S_n$.
On other hand, conjugacy classes of involutions are in one-to-one
correspondence with decompositions of $n$ into the sum $n=k+2l$.

Group $\Aut R$ acts trivially on the set ${\rm Cov_n}(R)$
of isomorphism classes of coverings. 
Therefore, topological type of a covering can be also presented as
a conjugacy class of an involutive permutation.

%%%%%%%%%%%%%%%%%%%%%%%%%%%%%%%%%%%%%%%%%%%%%%%%%%%%%%%%%%%%%%%%%%%

	\subsubsection{Coverings over a segment} (Fig.5)

Formal description of isomorphism classes of $n$-sheeted  coverings
over an interval $I$ is as follows.
Denote by $p_1$ and $p_2$ end points of $I$.
Let $\pi:P\to I$ be an $n$-sheeted  covering.
Denote by $M$ the set of connected components of the preimage of $I^\circ$ 
in $P$. Obviously, $\#M=n$ and the restriction of 
$\pi$ to any component is a homeomorphism. 

A point $p_i$ ($i=1,2$) defines an involutive permutation of set $M$,
which we denote by $s_i$. Thus, 
to each special covering $\pi$  we assign the
ordered pair $(s_1,s_2)$ of involutions $s_1,s_2\in\Aut M\approx S_n$. Two pairs
$(s_1,s_2)$ and $(s_1',s_2')$ are called conjugated if there exists $g\in S_n$
such that $s_1'=gs_1g^{-1}$ and $s_2'=gs_2g^{-1}$.

\medskip
\epsfxsize=5truecm
\begin{figure}[tbph]
\centerline{\epsffile{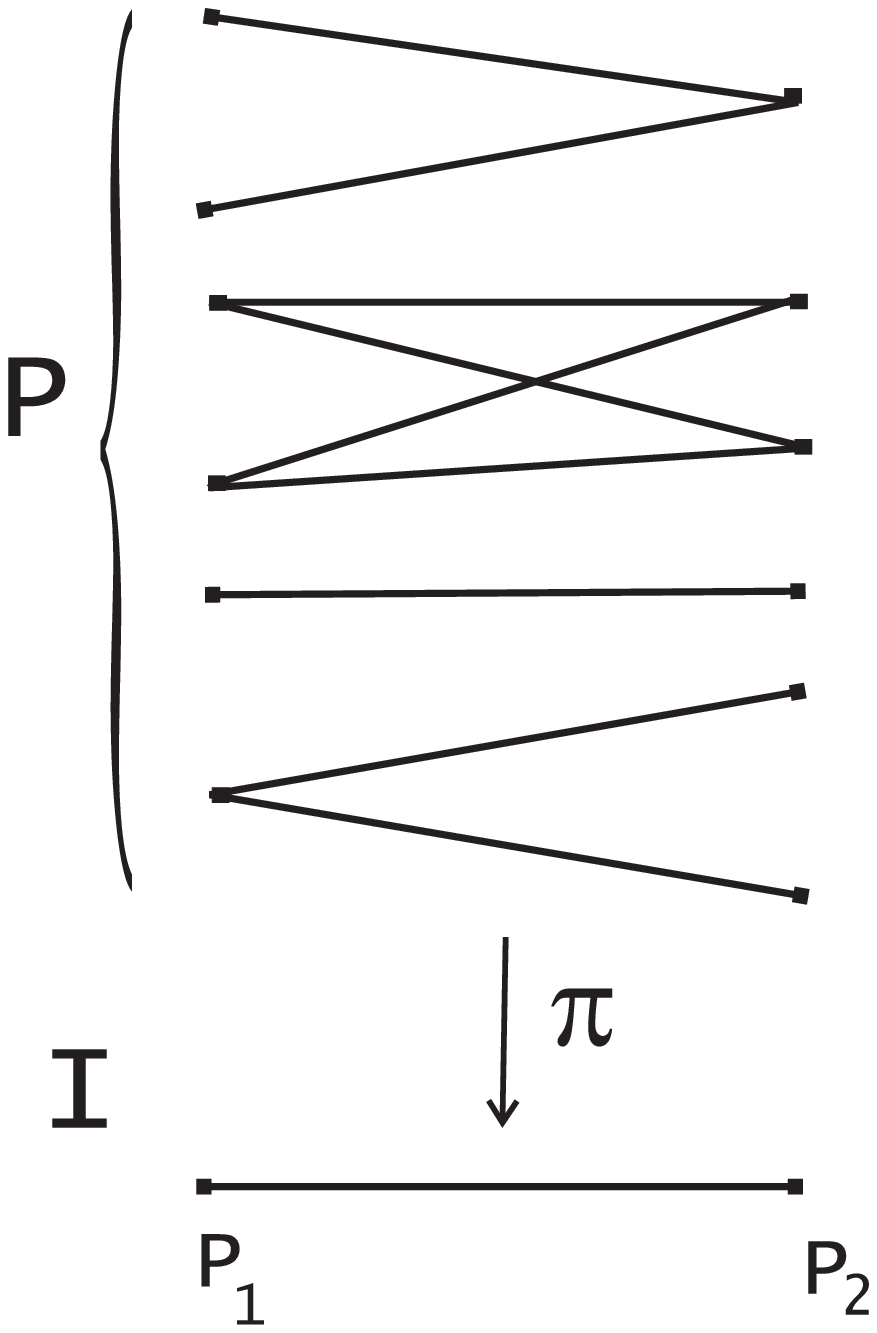}}
\caption{Stratified covering $\pi$ over segment $I$.}
\end{figure}
\medskip

\begin{lemma} There exists a natural one-to-one 
correspondence between set ${\rm Cov_n}(I)$ of isomorphism classes 
of $n$-sheeted coverings over an interval  $I$ and the set 
of conjugacy classes of pairs  of involutions in $S_n$.
\end{lemma}
Proof is skipped.

%%%%%%%%%%%%%%%%%%%%%%%%%%%%%%%%%%%%%%%%%%%%%%%%%%%%%%%%%%%%%%%%%%%
  
\subsubsection{Analogue of Young diagram for pairs of involutions}

The analogue of Young diagram can be defined for an ordered 
pair of involutions. Let $D_\infty=\langle \Sigma_1, \sigma_2 | \sigma_1^2=
\sigma_2^2=1\rangle$ be the infinite  dihedral group and $(s_1,s_2)$ be an
ordered pair of involutive permutations of the symmetric group $S_n$ acting
by permutations of the set $M=\{1,\dots,n\}$. Then the correspondence 
$\sigma_i\mapsto s_i$ can be uniquely extended to the representation $\rho$
of $D_\infty$ by permutations of the set  $M$.
Thus, $M$ can be decomposed into the orbits $M_1,\dots,M_l$ of $D_\infty$.
Any orbit can be considered as a transitive representation of $D_\infty$. 

Let $M_i$ be an orbit of $D_\infty$. Then it can be  of one of  four types:
\begin{itemize}
\item (type 1) $\sigma_1$ and $\sigma_2$ has no fixed points in $M_i$; 
\item (type 2) $\sigma_1$ has one fixed point and $\sigma_2$ has one fixed point;
\item (type 3) $\sigma_1$ has two fixed point and $\sigma_2$ has no fixed point;
\item (type 4) $\sigma_1$ has no fixed points and $\sigma_2$ has two fixed points. 
\end{itemize}

An orbit of type 1 consists of even number $2n_i$ of elements  and any two
orbits of type 1 of $2n_i$ elements are equivalent 
as representations of $D_\infty$.
An orbit of type 2 consists of odd number $2n_i-1$ of elements  and any two
orbits of type 2 of $2n_i-1$ elements are equivalent 
as representations of $D_\infty$.
An orbit of type 3 consists of even number $2n_i$ of elements  and any two
orbits of type 3 of $2n_i$ elements are equivalent as representations of
$D_\infty$. An orbit of type 4 consists of even number $2n_i$ of elements
and any two orbits of type 4 of $2n_i$ elements are equivalent as
representations of $D_\infty$.

Thus, to each pair $(s_1,s_2)$ we can assign a decomposition of $n$ into
the sum of natural numbers; the sum is separated into four blocks as follows:
$n=(2n_1^1+\dots+2n_{s_1}^1)+((2n_1^2-1)+\dots+(2n_{s_2}^2-1))+(2n_1^3+\dots+
2n_{s_3}^3)+(2n_1^4+\dots+2n_{s_4}^4)$ Blocks correspond to types of orbits,
summands in a block are unordered. We denote this decomposition by a single
letter $\beta$ and call it a {\it dihedral decomposition of $n$}. 
 
\begin{lemma} Let $(s_1,s_2)$ and $(s_1',s_2')$ be two ordered pairs of
involutive elements of $S_n$ and $\beta$, $\beta'$ be corresponding dihedral 
decompositions of $n$.
Then $(s_1,s_2)$ is conjugated to $(s_1',s_2')$ if and only if $\beta=\beta'$.
\end{lemma}         
Proof is skipped.

Thus, dihedral decompositions of $n$ are in bijection with conjugacy
classes of ordered pairs of involutions. Obviously, one can present $\beta$
also as the set of four Young diagrams with easily derived special properties.
We call this presentation a dihedral Young diagram of a pair $(s_1,s_2)$.
In sequel we consider terms 'a dihedral decomposition of $n$' and 'a dihedral Young
diagram of order $n$' as synonyms. Let $\beta$ be the dihedral Young diagram
of a pair $(s_1,s_2)$. Denote by $\beta^*$ a dihedral Young diagram of the
pair $(s_2,s_1)$.
The operation $\beta\mapsto\beta^*$ defines an involutive automorphism 
of the set ${\rm Cov_n}(I)$, which we denote by the same sign "$^*$".

\begin{lemma} Dihedral Young diagram $\beta^*$ can be obtained
from $\beta$ by replacing blocks of summands of types 3 and 4.
\end{lemma}

Proof is elementary. 

Assign to each conjugacy class $[(s_1,s_2)]$ conjugacy classes of
involutions $[s_1]$ and $[s_2]$ and denote them by 
$\iota_1(\beta)$ and $\iota_2(\beta)$ respectively.

A conjugacy classes of a pair $[(s,s)]$, $s\in S_n$, is called
{\it a trivial
class}; corresponding dihedral Young diagram is said to be trivial.
Evidently, they are in one-to-one correspondence with
conjugacy classes of involutive elements in $S_n$. The dihedral Young diagram
$\beta$ of the class $[(s,s)]$ is $n=(2+\dots+2)+(1+\dots+1)$ (blocks of
summands of types 3 and 4 are empty).

Fix an orientation of interval $I$. Then homeomorphisms of $I$
preserving the orientation fixes points $p_1$ and $p_2$. These homeomorphisms 
generate the index two subgroup $\Aut^+I$ of automorphism group $\Aut I$.
Evidently, $\Aut^+I$ acts trivially on the set ${\rm Cov}(I)$. Thus,
oriented topological  types of special coverings over a segment $I$ 
are in one-to-one correspondence with dihedral Young diagrams.

Any element from $\Aut I\setminus\Aut^+I$ acts on the set ${\rm Cov}(I)$
as $(s_1,s_2)\mapsto (s_2,s_1)$. Thus, coverings with dihedral Young
diagrams $\beta$ and $\beta^*$ are of the same topological  type. 
Therefore, topological type and isomorphism class of coverings do not
coincide for the segment. 

\subsubsection{Dimension 2} (Fig.6)

Let $\pi:P\to\Omega$ be a covering over a stratified surface
$\Omega=\coprod_{\lambda\in\Lambda}\Omega_\lambda$ and 
$p\in\Omega$ be any point of base. Topological invariants
of the restriction $\pi_U:P_U\to U$ of the covering to a sufficiently
small neighborhood $U$ of $p$ depends only on a stratum
containing $p$. More precisely, if $p$ belongs to
the generic stratum, then
$\pi_U$ is trivial covering, i.e., $\pi_U$ mappings $n$ copies 
of $U$ onto $U$. Thus, $n$ is a unique invariant at a generic point.

If $p$ is an interior special point, then one can choose
$U$ isomorphic to an open disc. In this case  $\pi_U$ is defined,
up to topological equivalence, by the topological type 
of its restriction to the boundary $\partial U$ (it is a circle) \cite{Ker}.
Thus, a Young diagram describes the covering in the neighborhood of $p$.

If $p$ belongs to a one-dimensional boundary stratum $E$, then it is
convenient to choose $U$ homeomorphic to the direct product  
of a neighborhood $U_E$ of $p$ in $E$ and a ray $R$ . Clearly, 
topological type of the covering $\pi_U$ is uniquely defined 
by the topological type of its restriction to the ray $R$. Therefore,
the covering in the neighborhood of $p\in E$ is described by
a conjugacy class of an involutive element of $S_n$.

If $q$ is a boundary special point, 
then it is convenient to choose a neighborhood $U$ homeomorphic to 
an open cone ${\mathcal C}^\circ(I)$ over a closed interval $I$.
Clearly, isomorphism classes of coverings over $U$ are in one-to-one
correspondence with isomorphism classes of coverings over $I$. 
As it was mentioned above, topological type of a covering over $I$
can include more than one isomorphism classes of coverings over  $I$. 
In order to control isomorphism classes of 
coverings over $U$, let us  fix an admissible set of local 
orientations at special points of the surface.
Hence, to each boundary special point $q$
we can assign an isomorphism class of coverings over the segment, 
i.e., dihedral Young diagram. This diagram describes the covering in the 
neighborhood of $q$.

Clearly, local invariants of a covering coincide for all points from the same 
stratum. Thus, to each  stratum of $\Omega$ we assigned a combinatorial 
invariant.

\medskip
\epsfxsize=9truecm
\begin{figure}[tbph]
\centerline{\epsffile{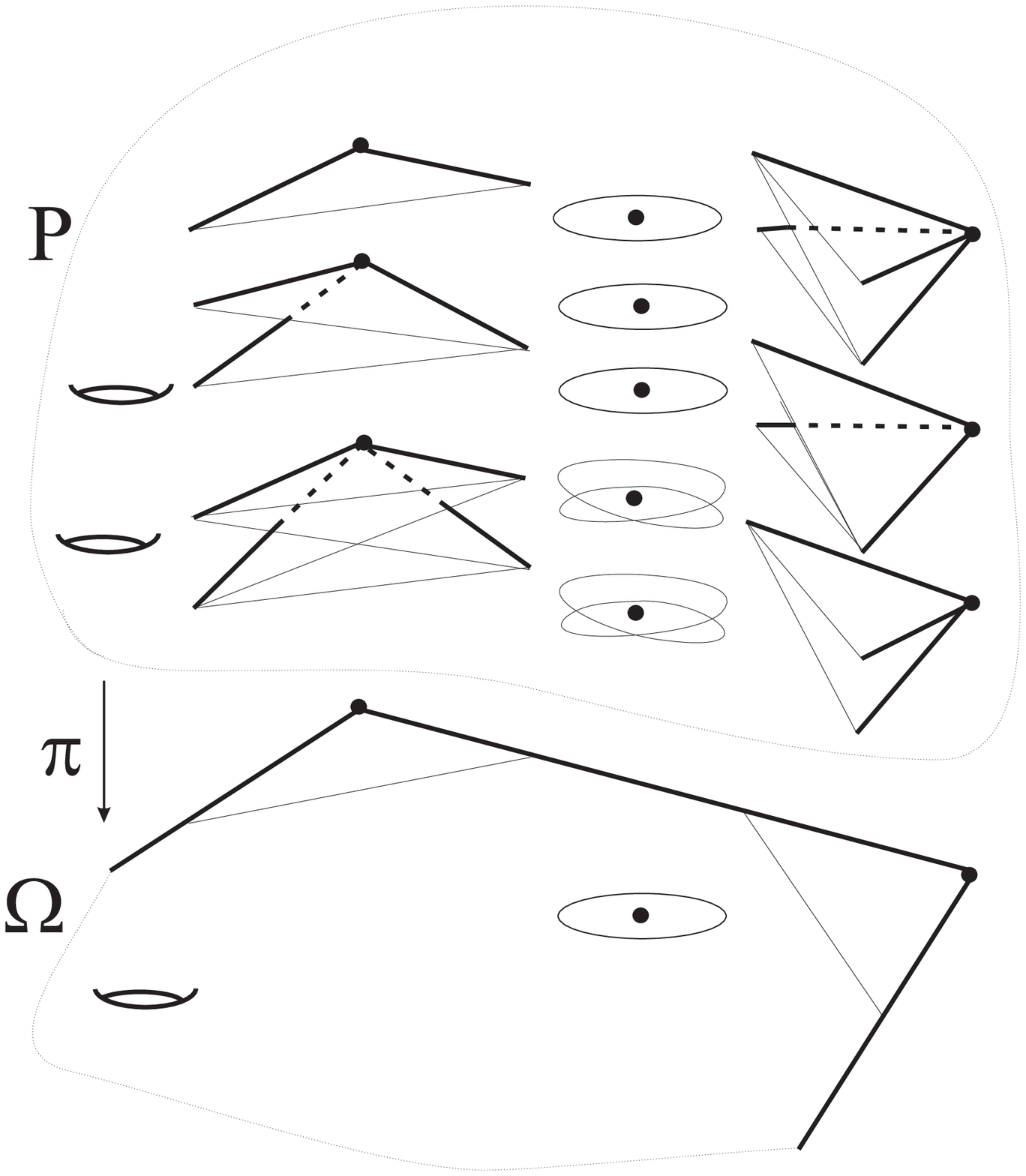}}
\caption{Stratified covering $\pi$ over two-dimensional stratified 
surface $\Omega$.}
\end{figure}
\medskip

\subsection{Cut coverings} 

Recall that a cut surface $\Omega_*=(\Omega_*,\Gamma_*,\tau)$ is a 
triple, consisting of 

a stratified surface $\Omega_*$, 

a set $\Gamma_*\subset\Omega_*$ such that $\Gamma_*$ coincides with the joint of 
several closed pairwise nonintersecting one-dimensional boundary strata,

an involutive homeomorphism $\tau:\Gamma_*\to\Gamma_*$ having no fixed points
(see section \ref{cutsys}).
  
Fix a cut surface $\Omega_*=(\Omega_*,\Gamma_*,\tau)$. 

Let $\pi_*:P_*\to \Omega_*$ be an $n$-sheeted covering. 
Denote by  $\widehat{\Gamma}_*$ the preimage of $\Gamma_*$ in $P_*$.
Suppose that the following condition holds:
$$
\widehat{\Gamma}_*\subset\partial P_*
\leqno (*)
$$
In sequel we deal with coverings $\pi_*$ satisfying this condition.

Denote by 
$\xi: \widehat{\Gamma}_*\to\Gamma_*$ the restriction
of $\pi_*$ to $\widehat{\Gamma}_*$. 
Due to $(*)$ the mapping $\xi:\widehat{\Gamma}_*\to\Gamma_*$ is an $n$-sheeted
covering over $\Gamma_*$. 

An involutive
homeomorphism
	$\widehat{\tau}:\widehat{\Gamma}_*\to\widehat{\Gamma}_*$ 
is called {\it an admissible involution} if 
	$\pi_*\circ\widehat{\tau}=\tau\circ\pi_*$.

\begin{definition} Let $\pi_*:P_*\to\Omega_*$ be a 
covering over a cut surface $\Omega_*$. Suppose condition $(*)$ is satisfied. 
Let $\widehat{\tau}:\widehat{\Gamma}_*\to\widehat{\Gamma}_*$ be 
an admissible involution. Then pair $(\pi_*,\widehat{\tau})$ is called 
a cut covering.
 
An isomorphism of two cut coverings
$(\pi_*,\widehat{\tau})$ and $(\pi_*',\widehat{\tau}')$ is an isomorphism 
$\phi:P_*\to P_*'$ of coverings such that $\widehat{\tau}'\circ\phi=\phi\circ
\widehat{\tau}$.
\end{definition}

Denote by ${\mathcal CutCov}(\Omega_*)$
the category of cut coverings over
a  cut surface $\Omega_*=(\Omega_*,\Gamma_*,\tau)$ with morphisms 
being isomorphisms of cut coverings;
denote by ${\rm CutCov}(\Omega_*)$ 
the set of isomorphism classes of cut coverings.
In sequel, we shall use  the same difference
in fonts in order to distinguish a category and the set of
isomorphism classes of objects in it and in other cases.

Clearly, ${\rm CutCov}(\Omega_*)$ is a finite set.

%%%%%%%%%%%%%%%%%%%%%

Let $\pi:P\to \Omega$ be a covering over stratified surface $\Omega$, 
$\Gamma\subset \Omega$ be a cut system, 
$(\Omega_*, \Gamma_*,\tau)$ be a cut surface
obtained by cutting $\Omega$ along $\Gamma$ and $\glue: \Omega_*\to\Omega$
be gluing map.

Clearly, the  preimage $\widehat{\Gamma}=\pi^{-1}(\Gamma)$
is a cut system of $P$. 
Denote by $(P_*, \widehat{\Gamma}_*,\widehat{\tau})$ 
the cut surface obtained by cutting  $P$ along $\widehat{\Gamma}$ and denote 
by ${\widehat{\glue}}$ the glueing map $P_*\to P$.
Obviously, there is a natural covering $\pi_*:P_*\to \Omega_*$ such that 
the following diagram is commutative.
$$
 \begin{CD}
P_*               @>{\widehat{\glue}}>>         P\\
@V{\pi_*}VV                    @V{\pi}VV\\
\Omega_*          @>>{\glue}>    {\Omega}\\ 
 \end{CD}
$$

Evidently, covering $\pi_*$ satisfies $(*)$ and the correspondence 
$\pi\mapsto \pi_*$ is a functor           
$\cut :{\mathcal Cov}(\Omega)\to{\mathcal Cov}(\Omega_*)$. This 
functor generates a map $\cut:{\rm Cov}(\Omega)\to {\rm Cov}(\Omega_*)$ of 
isomorphism classes of coverings. Functor $\cut$ is
a superposition of two functors 

(1) ${\mathcal Cov}(\Omega)\to{\mathcal CutCov}(\Omega_*)$ and 

(2) ${\mathcal CutCov}(\Omega_*)\to {\mathcal Cov}(\Omega_*)$,

where (1) is induced by cutting along $\Gamma$ and (2) is
just "forgetting" $\widehat \tau$.

\begin{lemma} 
\label{a6.5}           
Functor ${\mathcal Cov}(\Omega)\to {\mathcal CutCov}(\Omega_*)$ 
induces the equivalence of categories.
\label{equivalence}
\end{lemma}

\begin{proof} One can easily construct an inverse functor
${\mathcal CutCov}(\Omega_*)\to{\mathcal Cov}(\Omega)$.
It 'glues'  a cut covering $(\pi_*,\widehat\tau)$ into
a covering $\pi:P_*/\widehat{\tau}\to\Omega$.
\end{proof}

Let $(\Omega_*, \Gamma_*,\tau)$ be a cut surface.
Denote by ${\mathcal Cov}(\Omega_*,\Gamma_*)$
the subcategory of  ${\mathcal Cov}(\Omega_*)$, that consists of 
coverings $\pi_*:P_*\to \Omega_*$ over $\Omega_*$ 
satisfying condition $(*)$.

Fix a covering $\pi_*\in{\mathcal Cov}(\Omega_*,\Gamma_*)$;
denote by  $[\pi_*]$ its class.
Denote by $\xi$ the restriction $\pi_*\vert_{\widehat{\Gamma}_*}$. 

Clearly, an automorphism of covering $\pi_*$ induces 
the automorphism of covering $\xi$.  Hence, 
we obtain a homomorphism $\phi:\Aut\pi_*\to\Aut\xi$.
The preimage of a generic point of $\Gamma_*$ consists 
of $n$ elements where $n$ is the degree of $\pi_*$.
Therefore, $\phi$ is a monomorphism.
We shall consider group $\Aut\pi_*$ as a subgroup of $\Aut\xi$.

Denote by $T(\pi_*)$ the set of all admissible involutions
$\widehat{\tau}:\widehat{\Gamma}_*\to \widehat{\Gamma}_*$.

Clearly, the composition $\widehat{\tau}\circ\widehat{\tau}'$ of
two admissible involutions is an automorphism of covering
$\xi$ and vice versa, the composition of an admissible involution
and an automorphism of $\xi$ is an admissible involution. Hence,
the number of admissible involutions $|T(\pi_*)|$ coincides 
with the order of group $\Aut\xi$.

Group $\Aut\xi$ acts on $T(\pi_*)$ as follows:
$g(\widehat{\tau})(x)=g\circ\widehat{\tau}\circ g^{-1}(x)$, where
$g\in\Aut\xi$, $\widehat{\tau}\in T(\pi_*)$, $x\in\widehat{\Gamma}_*$.

Denote by ${\rm CutCov}(\Omega_*)_{[\pi_*]}$ the 
preimage of $[\pi_*]$ under the map
${\rm CutCov}(\Omega_*)\to{\rm Cov}(\Omega_*)$  (coming from 
functor (2) above). 
\begin{lemma} 
\label{6.6}
There exists a bijection between the set
${\rm CutCov}(\Omega)_{[\pi_*]}$ and the set of orbits
$T(\pi_*)/\Aut(\pi_*)$.
\end{lemma}

\begin{proof} Choose a representative 
$\pi_*\in[\pi_*]$. Any class from 
${\rm CutCov}(\Omega_*)_{[\pi_*]}$ contains
a representative $(\pi_*,\widehat{\tau})$ that 
includes just the covering $\pi_*$.
Clearly, cut coverings $(\pi_*,\widehat{\tau})$ and $(\pi_*,\widehat{\tau}')$ 
with the same $\pi_*$ are
isomorphic, if and only if $\widehat{\tau}$ and $\widehat{\tau}'$ belong
to the same orbit of $\Aut\pi_*$.
\end{proof}

Evidently, isomorphic coverings have isomorphic automorphism groups. 
Therefore, we can assign the isomorphism class of groups $\Aut [\pi]$
to an isomorphism class of coverings $[\pi]$.
Analogously, to $[\pi_*]$ we assign a class of isomorpic sets of admissible 
involutions $T([\pi])$.
 
 By lemma \ref{a6.5} there is a bijection 
between sets ${\rm CutCov}(\Omega_*)$ and  ${\rm Cov}(\Omega)$. 
Denote by
${\rm Cov}(\Omega)_{[\pi_*]}$ the image of  
${\rm CutCov}(\Omega_*)_{[\pi_*]}$ in  ${\rm Cov}(\Omega)$.

\begin{theorem}\label{Main} Let $\Gamma$ be a cut system of 
a stratified surface $\Omega$. Denote by  
$(\Omega_*, \Gamma_*,\tau)$  a cut surface obtained by 
cutting $\Omega$ along $\Gamma$.
Let $[\pi_*]$  be an arbitrary class of isomorphic coverings 
$\pi_*:P_*\to\Omega_*$ satisfying $(*)$. Then 
$$\sum_{[\pi]\in{\rm Cov}(\Omega)_{[\pi_*]}}\frac{1}{\vert\Aut[\pi]\vert}=
\vert T([\pi_*])\vert\cdot\frac{1}{\vert\Aut[\pi_*]\vert}.$$
\end{theorem}

\begin{proof}
Clearly, the stabilizer of $\widehat{\tau}\in T(\pi_*)$ in group
$\Aut\pi_*$ coincides with the group of all automorphisms of 
cut covering $(\pi_*,\widehat{\tau})$. By lemma \ref{a6.5}, there is 
the bijection between cut coverings over cut surface $\Omega_*$ and 
coverings over $\Omega$. Let $\pi\in{\mathcal Cov}(\Omega)$ corresponds to
cut covering $(\pi_*,\widehat{\tau})$. 
Therefore, the stabilizer of $\widehat{\tau}$ in $\Aut\pi_*$
is isomorphic to $\Aut\pi$. Denote by $\Delta([\pi])$ the 
orbit in $T([\pi_*])$ corresponding to $[\pi]$ (see lemma \ref{6.6}). 
Hence
$|\Delta([\pi])|=\frac{\vert\Aut[\pi_*]\vert}{\vert\Aut[\pi]\vert}$.

In order to complete the proof substitute the latter equality 
into the equality 
$\vert T(\pi_*)\vert=
\sum_{[\pi]\in{\rm Cov}(\Omega)_{[\pi_*]}} | \Delta(\pi)|$ which
reflects the decomposition of $T(\pi_*)$ into orbits of $\Aut\pi_*$.
The claim of theorem coincides with obtained equality up 
to algebraic transformations.
\end{proof}

\note Theorem  \ref{Main} can be  easily generalized to 
stratified coverings in any dimension.
\medskip

Let $\Gamma$ be a cut system of a stratified surface $\Omega$.
Denote by $(\Omega_*, \Gamma_*,\tau)$ a cut surface
obtained by cutting $\Omega$ along $\Gamma$ and denote by
$\Omega_\#$ the contracted cut surface.

Let $\pi:P\to \Omega$ be a stratified covering. Denote by 
$(\pi_*,\widehat{\tau})$ the cut covering corresponding to $\pi$.
By definition, $\pi_*\in{\mathcal Cov}(\Omega_*,\Gamma_*)$.
Clearly, covering $\pi_*:P_*\to\Omega_*$ induces the covering 
$\pi_\#:P_\#\to\Omega_\#$ where $P_\#$ is the surface obtained 
by contracting to a point each  connected component 
of $\Gamma_*$. 

\begin{lemma} \label{a6.7} The correspondence 
$\pi_*\to\pi_\#$ induces the 
bijection between sets  ${\rm Cov}(\Omega_*,\Gamma_*)$
and ${\rm Cov}(\Omega_\#)$ of classes of equivalent coverings.
\end{lemma}

Proof is skipped.
\medskip

Let $\Gamma$ be a cut system of a stratified surface $\Omega$.
The composition of maps 
${\rm Cov}(\Omega)\to{\rm CutCov}(\Omega_*)
\to{\rm Cov}(\Omega_*,\Gamma_*)\to{\rm Cov}(\Omega_\#)$
defines the map ${\rm Cov}(\Omega)\to{\rm Cov}(\Omega_\#)$.
By theorem  \ref{Main} and lemma \ref{a6.7}, there is a relation
between number $\frac{1}{\vert\Aut[\pi_\#]\vert}$ where 
$[\pi_\#]\in {\rm Cov}(\Omega_\#)$
and numbers $\frac{1}{\vert\Aut[\pi]\vert}$ for all preimages 
$[\pi]\in{\rm Cov}(\Omega)$ of $[\pi_\#]$.
We shall rewrite them below in slightly 
different form for simple cuts of all classes.

Let $\gamma$ be a simple cut of $\Omega$. Thus, in our 
notations $\Gamma=\{\gamma\}$. We will write $\gamma$, $\gamma_*$, 
etc. instead of $\Gamma$, $\Gamma_*$ etc.

Clearly, there are three possibilities for $\gamma_*$: 
(1) $\gamma_*$ consists of two ovals;
(2) $\gamma_*$ consists of two segments;
(3) $\gamma_*$ consists of one oval.

Denote by $r'$, $r''$ two special points of $\Omega_\#$ coming
from components of $\gamma_*$ in cases (1) and (2) and by $r$
the analogous point in case (3).

We have to fix local orientations 
at points $r',r''$ or at $r$. It is required that
these local orientations can be completed  to an 
admissible set of local orientations at all special 
points of $\Omega_\#$. In any case there are 
local orientations at $r',r''$ or $r$ (case (3))
induced by an orientation of simple cut $\gamma$. Note,
that these local orientations 
may not satisfy the requirement. We correct them 
by the following rule. If $\gamma$ is a simple cut 
of class 1, 2, 5, 6, 7 than change the local orientation 
at point $r''$ by opposite one. Fix obtained local orientations
and denote them by $\mathcal o',\mathcal o''$ or $\mathcal o$ (case (3)).
It can be checked that these orientations satisfy 
the requirement. 

Let $\pi$ be a covering over $\Omega$ and $\gamma_*'$ be 
a connected component of $\gamma_*$. There is an orientation
of $\gamma_*'$ that is compatible with fixed local orientation
at point $r'\in\Omega_\#$ coming from the contraction of $\gamma_*'$.
Denote by $\xi'$ the restriction of 
covering $\pi_*$ to the preimage $\pi_*^{-1}
(\gamma_*')$. Thus,
we have that $\xi'$  is defined, up to isomorphism, either by 
a Young diagram or by a dihedral Young diagram. 

Therefore, we obtain a local invariant of covering 
$\pi_\#:P_\#\to\Omega_\#$ at this point as either 
a Young diagram or a dihedral Young diagram.

Let $\pi_\#$ be an arbitrary covering over $\Omega_\#$. Denote by 
$\alpha'$, $\alpha''$ the Young diagram describing local invariants
of $\pi_\#$ at points $r'$, $r''$ resp. in case (1). Analogously,
denote by $\beta'$, $\beta''$ dihedral Young diagrams in case (2) and by 
$\alpha$ Young diagram in case (3).

%It follows from the construction, that local invariant of $\pi_\#$
%at point $\gamma_\#'$ coincides with (oriented) topological type of
%$\xi'$.
Denote by ${\rm Cov}(\Omega)_{[\pi_\#]}$ the preimage of 
$[\pi_\#]\in{\rm\Omega}_\#$ in ${\rm Cov}(\Omega)$.

Put $\alpha^*=\alpha$ for a Young diagramm $\alpha$. We use 
$\alpha^*$ below in order to obtain forumas that will be correct also
for coverings with structure group $G$ different from $S_n$.

\begin{corollary}
\label{cor5.1} Let $\gamma$ be a simple cut of a stratified surface 
$\Omega$ and $\pi_\#$ be a covering over contracted cut 
surface $\Omega_\#$. Then the following identities hold:
\\ (1) If $\gamma$ belongs to one of classes 1, 2 then
$$\sum_{[\pi]\in{\rm Cov}(\Omega)_{[\pi_\#]}}\frac{1}{\vert\Aut[\pi]\vert}=
\delta_{\alpha',{\alpha''}^*}
\vert \Aut\alpha'\vert\cdot\frac{1}{\vert\Aut[\pi_\#]\vert}$$
\\ (2) If $\gamma$ belongs to class 3 then
$$\sum_{[\pi]\in{\rm Cov}(\Omega)_{[\pi_\#]}}\frac{1}{\vert\Aut[\pi]\vert}=
\delta_{\alpha',\alpha''}
\vert \Aut\alpha'\vert\cdot\frac{1}{\vert\Aut[\pi_\#]\vert}$$
\\ (3) If $\gamma$ belongs to one of classes 5,6,7 then
$$\sum_{[\pi]\in{\rm Cov}(\Omega)_{[\pi_\#]}}\frac{1}{\vert\Aut[\pi]\vert}=
\delta_{\beta',{\beta''}^*}
\vert \Aut\beta'\vert\cdot\frac{1}{\vert\Aut[\pi_\#]\vert}$$
\\ (4) If $\gamma$ belongs to one of classes 8,9 then
$$\sum_{[\pi]\in{\rm Cov}(\Omega)_{[\pi_\#]}}\frac{1}{\vert\Aut[\pi]\vert}=
\delta_{\beta',\beta''}
\vert \Aut\beta'\vert\cdot\frac{1}{\vert\Aut[\pi_\#]\vert}$$
\\ (5) If $\gamma$ belongs to class 4 then
$$\sum_{[\pi]\in{\rm Cov}(\Omega)_{[\pi_\#]}}\frac{1}{\vert\Aut[\pi]\vert}=
d_{\alpha}\cdot\frac{1}{\vert\Aut[\pi_\#]\vert}$$
where $d_\alpha$ is a number of involutions in symmetric group 
$S_n$ having no fixed points and commuting with an element $g\in S_n$ 
that belongs to conjugated class corresponding to Young diagram $\alpha$.
\end{corollary}

Proof follows from theorem  \ref{Main} and lemma \ref{a6.7}.

\subsection{Hurwitz topological field theory}

Denote by $\mathcal A=\mathcal A_n$ the set of Young diagrams of order $n$
and denote by $A$ a vector space of  formal linear combinations of Young
diagrams. Denote by $\mathcal B=\mathcal B_n$ the set of dihedral 
Young diagrams of order $n$ and denote by  $B$ a vector space of 
formal linear combinations of dihedral Young diagrams.
Define an involutive linear transformation $^*:A\to A$ as
identical map.  Define an involutive linear 
transformation $^*:B\to B$ as linear continuation
of the map $\beta\mapsto\beta^*$ for dihedral Young diagram $\beta$. 

Let $(\Omega,\mathcal O)$ be a pair consisting of a 
stratified surface and a set of local orientations $\mathcal O$  
at special points of $\Omega$. Denote by 
$\Omega_a$ the set of interior special points of $\Omega$ and 
by $\Omega_b$ the set of boundary special points.

Let  $\alpha:\Omega_a \to \mathcal A$ and 
$\beta :\Omega_b\to\mathcal B$ be two maps.
Denote by $\alpha_p$ (resp., $\beta_q$) the image 
of an interior special point $p$ (resp., boundary special point $q$)
in set $\mathcal A$ (resp.,$\mathcal B$).
Denote by ${\rm Cov}(\Omega,\{\alpha_p\},\{\beta_q\})$
the set of isomorphism classes of coverings 
having a local invariant $\alpha_p$ at each point $p\in\Omega_a$ and 
a local invariant $\beta_q$ at each point  $q\in\Omega_b$. 

We shall define linear function
$H_{(\Omega,\mathcal O)}: V_{(\Omega,\mathcal O)}\to \mathbb C$,
where $V_{(\Omega,\mathcal O)}=
(\otimes_{p\in\Omega_a}A_p)\otimes(\otimes_{q\in\Omega_b}B_q)$.

Assume first that $\mathcal O$ is an admissible set of local orientations.
Then for element 
$(\otimes_{p\in\Omega_a}\alpha_p)\otimes(\otimes_{q\in\Omega_b}\beta_q)\in
V_{(\Omega,\mathcal O)}$ put                                                     
$$
H_{(\Omega,\mathcal O)}(\otimes_{p\in\Omega_a}\alpha_p)
\otimes(\otimes_{q\in\Omega_b}\beta_q))=
\sum_{[\pi]\in {\rm Cov}(\Omega,\{\alpha_p\},\{\beta_q\})}\frac{1}{\vert\Aut[\pi]
\vert}.
$$

Elements $(\otimes_{p\in\Omega_a}\alpha_p)\otimes
(\otimes_{q\in\Omega_b}\beta_q)$
form a basis of $V_{(\Omega,\mathcal O)}$, hence we can expand 
$H_{(\Omega,\mathcal O)}$ by linearity.
For nonadmissible $\mathcal O$
define linear form $H_{(\Omega,\mathcal O)}$ 
in such a way that  axiom $2^\circ$ of Klein field theory is
satisfied. This rule correctly defines $H_{(\Omega,\mathcal O)}$ 
in all cases because $\sum_{[\pi]\in {\rm Cov}(\Omega,\{\alpha_p\},\{\beta_q\})}
\frac{1}{\vert\Aut\pi\vert}$ does not depend on the choice of 
an admissible set of local orientations.

\begin{theorem} The set of data 
$\mathcal H=\{A,^*:A\to A,B,^*:B\to B, H_{(\Omega,\mathcal O)}\}$
is a Klein topological field theory.
We call it Hurwitz topological field theory of degree $n$.
\end{theorem}

\begin{proof} Evidently, axioms $1^\circ$, $2^\circ$, $3^\circ$ and $4^\circ$
are satisfied. Bilinear forms can be easely computed. 
We obtain equalities: 
$(\alpha',\alpha'')_A=\delta_{\alpha',{\alpha''}^*}\frac{1}{|\Aut\alpha'|}$ and
$(\beta',\beta'')_B=\delta_{\beta',{\beta''}^*}\frac{1}{|\Aut\beta'|}$.
Element $U\in A$ is equal to 
$\sum_{\alpha}d_{\alpha}(\alpha,\alpha)_A$ where numbers 
$d_\alpha$ are defined in corollary \ref{cor5.1}. Cleary, there 
are the following identities for tensors corresponding to 
bilinear forms and $U$:
$F^{\alpha',\alpha''}=\delta_{\alpha',{\alpha''}^*}|\Aut\alpha'|$,
$F^{\beta',\beta''}=\delta_{\beta',{\beta''}^*}|\Aut\beta'|$,
$D^{\alpha}=d_\alpha$.

Axiom $5^\circ$  essentially follows from theorem  \ref{Main} and 
corollary \ref{cor5.1}. 

For example, let us verify the axiom for $\gamma$ of class 5.

Let $\Omega$ be a stratified surface, $\mathcal O$ be a set
of local orientations at special points and  
and  $\gamma$ be a simple cut of $\Omega$ of class 5, i.e.,
a cut between two holes. Denote by $\Omega_\#$ 
the contracted cut surface. 

The set $(\Omega_\#)_0$ of all special points of $\Omega_\#$ consists of 
the image of set $\Omega_0$ and two additional boundary 
points $q'$, $q''$ coming from connected 
components of  $\gamma_*$. 
Choose a set of local orientations $\mathcal O$ such that the following
conditions are satisfied. First, $\mathcal O$ 
is an admissible set of local orientations. Second, induced 
set of local orientations $\mathcal O_\#$ at special points of
$\Omega_\#$ coincides with an admissible set of local
orientations $\mathcal O_{\#adm}$ in all points except $q''$.
Clearly, for any $\Omega$ and $\gamma$ of class 5 it is possible to
choose $\mathcal O$ satisfying these conditions. (Note, that  
there is no admissible set of local orientations  $\mathcal O$ such that
$\mathcal O_\#$ is also admissible set of local orientations.)

Fix an interior primary field $\alpha_p$ at each interior special 
point $p$ of $\Omega$ and a dihedral primary field $\beta_q$ 
at each boundary special point $q$ of $\Omega$. Denote by 
$x=(\otimes_{p\in\Omega_a}\alpha_p)\otimes 
(\otimes_{q\in\Omega_b}\beta_q))$ the element of 
vector space $V_{(\Omega,\mathcal O)}$.
Then, by definition, element $\eta_*(x)\in V_{(\Omega,\mathcal O)}$ 
is equal to
$y=\sum_{\beta',\beta''}F^{\beta'\beta''}
(\otimes_{p\in\Omega_a}\alpha_p)
\otimes(\otimes_{q\in\Omega_b}\beta_q)\otimes \beta'\otimes {\beta''}^*$.

Sustitute the values $F^{\beta'\beta''}$:
$y=\sum_{\beta}|\Aut\beta|
(\otimes_{p\in\Omega_a}\alpha_p)
\otimes(\otimes_{q\in\Omega_b}\beta_q))\otimes \beta\otimes\beta$.

%$y(\alpha')=(\otimes_{p\in\Omega_a}\alpha_p)\otimes \alpha'\otimes \alpha'
%\otimes(\otimes_{q\in\Omega_b}\beta_q))\in V_{(\Omega_\#,\mathcal O_\#)}$,

%$x=(\otimes_{p\in\Omega_a}\alpha_p)\otimes 
%(\otimes_{q\in\Omega_b}\beta_q))\in V_{(\Omega,\mathcal O)}$.

Introduce notation: 

$R=H_{(\Omega_\#,\mathcal O_{\#})}(\eta_*(x))$

Substitute $y=\eta_*(x)$:

$R=\sum_{\beta}|\Aut\beta|H_{(\Omega_\#,\mathcal O_{\#})}
((\otimes_{p\in\Omega_a}\alpha_p)
\otimes(\otimes_{q\in\Omega_b}\beta_q)\otimes \beta\otimes\beta)$

By (already established) 
invariance of change of local orientations, we obtain:

$R=\sum_{\beta}|\Aut\beta|H_{(\Omega_\#,\mathcal O_{\#adm})}
((\otimes_{p\in\Omega_a}\alpha_p)
\otimes(\otimes_{q\in\Omega_b}\beta_q)\otimes \beta\otimes\beta^*)$

By definition of $H_{(\Omega_\#,\mathcal O_{\#adm})}$,

$R=\sum_{\beta}{\vert\Aut \beta\vert}
\sum_{[\pi_\#]\in {\rm Cov}(\Omega_\#,\{\alpha_p\},
\{\{\beta_q\},\beta,\beta^*\})}
                \frac{1}{\vert\Aut \pi_\#\vert}$

By theorem  \ref{Main} and corollary \ref{cor5.1} we have 

$R=\sum_{\beta}{\vert\Aut \beta\vert}
        \sum_{[\pi_\#]\in 
	{\rm Cov}(\Omega_\#),\{\alpha_p\},\{\{\beta_q\},\beta,\beta^*\})
	 \frac{1}{|\Aut\beta|}}
\sum_{[\pi]\in {\rm Cov}(\Omega)_{[\pi_*]}}
                \frac{1}{\vert\Aut \pi\vert}=$

$=\sum_{\beta}
	\sum_{[\pi_\#]\in
	{\rm Cov}(\Omega_\#),\{\alpha_p\},\{\{\beta_q\},\beta,\beta^*\})}
	\sum_{[\pi]\in{\rm Cov}(\Omega)_{[\pi_*]}}
                \frac{1}{\vert\Aut \pi\vert}=$

$=\sum_{[\pi]\in {\rm Cov}(\Omega)}
                \frac{1}{\vert\Aut \pi\vert}=
H_{(\Omega,\mathcal O)}(x)$.
\medskip

For other class of simple cuts axiom $5^\circ$ can be verified 
quit similar.

Axiom $6^\circ$ is satisfied because the sum 
$\sum_{[\pi]\in {\rm Cov}(\Omega,\{\alpha_p\},\{\beta_q\})}\frac{1}{\vert\Aut[\pi]
\vert}$ 
is multiplicative with respect to disjoint union of surfaces.

\end{proof}

%%%%%%%%%%%%%%%%%%%%%%%%%%

Structure algebra of Hurwitz topological field theory is a
vector space $H=A\oplus B$ with multiplication constructed 
by tensors from $\mathcal C(\mathcal H)$.

Let $\bar H=\bar A\oplus \bar B$ be a structure algebra of symmetric group.
Note that Young diagram $\alpha$ corresponds to a
class of conjugated elements of  $S_n$. This class we  denote
by the same sign $\alpha$.
Define a linear map $\phi:A\to\bar A$ as 
$\phi(\alpha)=E_\alpha$, where $E_\alpha=\sum_{g\in\alpha}g$ is an element
of the center of group algebra $\mathbb C[S_n]$.  
Analogously, dihedral Young diagram $\beta$ 
corresponds to a  conjugacy class of ordered pairs of
involution in $S_n$. Put $\phi(\beta)=E_\beta$, where 
$E_\beta=\sum_{(s_1,s_2)\in\beta}E_{s_1,s_2}$ and 
$E_{s_1,s_2}=((\delta_{s_1,s_2}))\in M(S_n)$.

\begin{theorem} \label{th53} The linear map $\phi$ is an isomorphism between
the structure algebra of Hurwitz topological field theory
of degree $n$ and the structure algebra of symmetric group $S_n$.
\end{theorem}

\begin{proof}

Obviously, $\phi$ defines an isomorphism of linear spaces
$H$ and $\bar H$. Moreover, $\phi$ preserves scalar products
on $A$ and $B$ and commutes with '$^*$' involutions.

In order to prove that $\phi$ is an isomorphism 
of structure algebras it is sufficient to check 
the following equalities:                                                  
\\ (1) $S_{\alpha_1,\alpha_2,\alpha_3}=S_{E_{\alpha_1},E_{\alpha_2},E_{\alpha_3}}$,
\\ (2) $T_{\beta_1,\beta_2,\beta_3}=T_{E_{\beta_1},E_{\beta_2},E_{\beta_3}}$,      
\\ (3) $R_{\alpha,\beta}=R_{E_\alpha,E_\beta}$.

Proof of (1). After Hurwitz, it is known that
computation of weighted numbers of isomorphism classes
of coverings over two-dimensional  sphere is equivalent to the solution 
of factorization problem in $S_n$ (see \cite{Go-Jac}). 

Therefore,
$S_{\alpha_1,\alpha_2,\alpha_3}=\vert(\{ (a_1,a_2,a_3) | a_i
\in \alpha_i, a_1a_2a_3=1\}/S_n)\vert$,
where $g\in S_n$ acts as $(a_1,a_2,a_3)\mapsto
(ga_1g^{-1},ga_2g^{-1},ga_3g^{-1})$.
By direct calculations, we obtain that
$E_{\alpha_1}E_{\alpha_2}=S_{\alpha_1,\alpha_2}^{\alpha_3}
 E_{\alpha_3}$, where
$S_{\alpha_1,\alpha_2}^{\alpha_3}=S_{\alpha_1,\alpha_2,\alpha_3}
\vert\Aut\alpha_3\vert$.
Hence, $S_{\alpha_1,\alpha_2,\alpha_3}=S_{E_{\alpha_1},E_{\alpha_2},E_{\alpha_3}}$.

Proof of (2).
By definition $T_{\beta_1,\beta_2,\beta_3}$ is a number of equivalent classes 
of coverings over triangle $\Delta$  with local invariants at vertices $q_i$
equal to $\beta_i$ ($i=1,2,3$). 

Let $\pi:P\to\Delta$ be a stratified covering. Evidently, preimage of 
$\Delta^\circ$ 
consists of $n$ connected components and each one if homeomorphic to 
$\Delta^\circ$.
An edge $e_{i,j}$ of $\Delta$ defines a involutive permutation $s_{i,j}$ of these 
preimages.
Therefore, the covering is uniquely
defined by triple of involutions  $(s_{0,1},s_{1,2},s_{2,0})$ 
of symmetric group $S_n$ up to equivalence 
$(s_{0,1},s_{1,2},s_{2,0})\sim (gs_{0,1}g^{-1},gs_{1,2}g^{-1},gs_{2,0}g^{-1})$. 
Conversely  an equivalent  class $[(s_{0,1},s_{1,2},s_{2,0})]$ of triple of 
involutions defines an equivalent class of coverings over $\Delta$. Local 
invariant at vertex $q_i$ is a dihedral Young diagram $\beta_i$ corresponding
to a conjugacy class $[(s_{i-1,i},s_{i,i+1})]$ (where $i-1$ and $i+1$ are 
taken modulo $3$) of the pair of involutions.

Hence, $T_{\beta_1,\beta_2,\beta_3}=
\vert(\{(s_{0,1},s_{1,2},s_{2,0}) | s_{i,i+1}\in S_n, s_{i,i+1}^2=1, 
[(s_{i-1,i},s_{i,i+1})]=\beta_i\}/\sim)\vert$.
By direct calculations we obtain that 
$E_{\beta_1}E_{\beta_2}=T_{\beta_1,\beta_2,\beta_3}
\vert\Aut\beta_3\vert E_{\beta_3^*}$.
Therefore,
$T_{\beta_1,\beta_2,\beta_3}=T_{E_{\beta_1},E_{\beta_2},E_{\beta_3}}$

Proof of (3).
By definition $R_{\alpha,\beta}$ is the number of equivalent classes 
of coverings over disc $(D,p,q)$ with local invariant $\alpha$  at 
interior special point $p$ and local invariant $\beta$  at 
boundary special point $q$. 

Let $\pi:P\to D$ be a stratified covering. Let us connect points $p$ and 
$q$ by a segment $\delta$ and denote by $D^\circ$ the set of all interior 
point except $\delta$. Evidently, the preimage of $D^\circ$ 
consists of $n$ connected components and each one is homeomorphic to $D^\circ$.

Segment $\delta$  defines a permutation $a$ of these preimages and 
the boundary circle defines an involutive permutation $s$ of these preimages.
(If $X$ is a connected component of the preimage of $D^\circ$ then 
$Y=s(X)$ is another  component such that 
$\overline{X}\cap\overline{Y}$ is equal to a connected component 
of the preimage of a boundary without point $q$; if there no such component 
then $s(X)=X$). Therefore, the covering is uniquely
defined by a pair $(a,s)$ 
up to equivalence $(a,s)\sim (gag^{-1},gsg^{-1})$. 

Conversely  an equivalent  class $[(a,s)]$
defines an equivalent class of coverings over $D$. Local 
invariant at vertex $p$ is  Young diagram $\alpha=[a]$ 
and local invariant at vertex $q$ is dihedral  Young diagram $\beta$ 
corresponding to a pair of involutions $(s,asa^{-1})$.
Hence, $R_{\alpha,\beta}=
\vert(\{(a,s) | a,s\in S_n, a\in\alpha, (s,asa^{-1})\in \beta \}/\sim)\vert$.
By direct calculations we obtain that 
$(E_\alpha,E_{\beta})=R_{\alpha,\beta}$ (see subsection 2.3)

\end{proof}

%In sequel we identify Young diagrams and dihedral Young diagrams
%with their images in $S_n$ by $\phi$.

\subsection{Hurwitz numbers}

Classical Hurwitz numbers are weighted numbers of algebraic maps 
$f:P\to Q$ to Riemann sphere $Q=S^2$ having  prescribed branchings 
at finite number of fixed points $z_1,\dots,z_m\in Q$.
Young diagramms $\alpha_1,\dots,\alpha_m$ are fixed. 
It is required that $z_1,\dots,z_m$ are critical values 
of $f$ and 
the branching of $f$ at $z_i$ is described by Young diagramm $\alpha_i$.

We use term 'generalized Hurwitz numbers' or 'Hurwitz numbers' for short
for morphisms of both complex and  real algebraic curves, an arbitrary base $Q$
and an arbitrary local types of branchings. Definitions are given below.

A real algebraic curve is a pair $(P, \tau)$, where $P$ is 
a complex algebraic curve, i.e., a compact Riemann surface, and $\tau:P\to P$ 
is an antiholomorphic involution \cite{Al-Gr, Nat3}. Involution $\tau$
is called complex conjugation.
Fixed points of $\tau$ are called  real points of the real curve
$(P, \tau)$. An algebraic map $f: (P, \tau )\to(Q, \omega)$ of real algebraic curves is
defined as a holomorphic map $f:P\to Q$ such that $f\tau=\omega f$.

Algebraic maps $f: (P, \tau )\to(Q, \omega)$ of a fixed degree $n$ 
to a fixed real algebraic curve 
form a category ${\mathcal Map}(Q,\omega)$; a  morphism $\phi:f\to f'$, where
$f: (P, \tau )\to(Q, \omega)$ and $f': (P', \tau' )\to(Q, \omega)$ are two 
algebraic maps to $(Q, \omega)$, is defined as 
a holomorphic map
$\phi:P'\to P$ such that $\phi\tau'=\tau\phi$ and $f'=f\phi$.

Any complex algebraic curve $P_c$ generates real algebraic curve 
$(P,\tau)$. Namely, $P=P_c\coprod\bar P_c$ where $\bar P_c$
is $P_c$ endowed with  conjugated complex structure and $\tau$
is just permutation of components $P_c$ and $\bar P_c$.
Thus, complex algebraic curves can be considered as a particular case 
of real algebraic curves. 
 
Denote by $\Omega$ the factor-space $Q/\tau$. It is known that 
$\Omega$ carries the structure of Klein surface \cite{Al-Gr, Nat1}.
This implies that
$\Omega$ is nooriented, possibly nonorientable, possibly with boundary, surface
endowed with an atlas such that any transaction function is dianalitic, i.e.,
either holomorphic or antiholomorphic function .
Morphisms 
of Klein surfaces are defined as dianalitic functions.
According to \cite{Al-Gr} the category of real algebraic 
curves is isomorphic to the category of Klein surfaces. 

Fix a real curve $(Q,\omega)$ and put $\Omega=Q/\tau$. 
Any morphism $f:(P,\tau )\to(Q,\omega)$ of
real algebraic curves generates the dianalytic map 
$\pi(f):P/\tau \to \Omega$. Denote by $p_1,\dots,p_m\in\Omega$ critical values of 
$\pi(f)$ that are interior points, and denote by $q_1^j,\dots,q_{m_j}^j$ 
critical values of $\pi(f)$ that belong to $j$-th boundary contour of 
$\Omega$. Clearly, points $p_i, q_i^j$  come from critical values 
$z_0,\dots,z_t\in Q$ of $f$: an interior point $p_i$ is the image of two points
$z_{k},z_{l}$, a boundary point $q_i^j$ is the image of one point
$z_{r}$ and $z_{r}$ is a fixed point of $\tau$.
Points $p_i, q_i^j$ induce the stratification of $\Omega$.
Forgetting dianalitic structures on $P/\tau$ and $\Omega$
we obtain that $\pi(f)$ is a stratified covering, i.e., an element 
of  ${\mathcal Cov}(\Omega)$. 

Fix the stratification of $\Omega$ with special points $p_i, q_i^j$. 
According to \cite{Nat2}
forgetting dianalitic structures induces the bijection 
between the set of isomorphism classes 
of dianalitic maps to $\Omega$, such that all critical values are 
special points of the 
stratification, and the set ${\rm Cov}(\Omega)$
of isomorphism classes of stratified covering over $\Omega$.
Thus, up to isomorphisms, algebraic maps to 
$(Q,\tau)$ can be identified with stratified coverings over $\Omega$.

The latter fact allows us to define Hurwitz numbers 
for real algebraic curves in topological terms.
Namely, 
define {\it Hurwitz numbers} for $\Omega$ as

$$
{\rm Hurw}_n(\Omega,\{\alpha_{p}\},\{\beta_{q}\})=
\sum_{[\pi]\in {\rm Cov}(\Omega,\{\alpha_{p}\},\{\beta_{q}\})}\frac{1}
{\vert\Aut\pi\vert}$$ 

Here $\alpha_p$ (resp., $\beta_q$) is a Young diagram 
(resp., a dihedral Young diagram) of degree $n$
assigned to an iterior special point $p$ (resp., to a boundary special point $q$) of 
$\Omega$.

By previous considerations, these Hurwitz numbers are weighted numbers 
of classes of alrebraic maps with prescribed branchings
to real algebraic curve $(Q,\tau)$.

For $G=(0,1)$ our definition gives weighted numbers of meromorphic
functions with fixed critical values of fixed topological types.
They includes numbers that were first introduced by Hurwitz \cite{Hur}.

By the results of previous subsection, 
Hurwitz numbers are equal to correlators for Hurwitz
topological field theory. 
 Thus, from theorem \ref{th53} and theorem
\ref{MMain} we get the following theorem.

\begin{theorem}
\label{th5.4} Let $\Omega$ be a stratified surface
of type $G$. Fix an admissible set of local orientations at 
special points of $\Omega$.Then  
\\ (1) $${\rm Hurw}_n(\Omega,\{\alpha_p\},\{\beta_q\})=
\langle \alpha_1,\dots,\alpha_m,({\beta_1}^1,\dots,\beta_{m_1}^1),
\dots,(\beta_1^s,\dots, \beta_{m_s}^s)\rangle_G$$
\\ (2) If $G=(g,1,m,m_1,\dots,m_s)$ then 
$${\rm Hurw}_n(\Omega,\{\alpha_p\},\{\beta_q\})=$$
$$=(\alpha_1\dots \alpha_m(\beta_1^1\dots \beta_{m_1}^1)V_{K_B}
(\beta_1^2\dots \beta_{m_1}^2)
\dots V_{K_B}(\beta_1^s\dots \beta_{m_s}^s),K_A^{g})$$
\\ (3) If $G=(g,0,m,m_1,\dots,m_s)$ then 
$${\rm Hurw}_n(\Omega,\{\alpha_p\},\{\beta_q\})=$$
$$=(\alpha_1\dots \alpha_m (\beta_1^1\dots \beta_{m_1}^1)V_{K_B}
(\beta_1^2\dots \beta_{m_1}^2)
\dots V_{K_B}(\beta_1^s\dots \beta_{m_s}^s),U^{2g})$$

In (2) and (3) right hand sides are expressions in the structure 
algebra associated with symmetric group $S_n$ (see subsection 2.3). 
\end{theorem}

\note 

Let $G$ be a finite group and $\Omega$ be a stratified surface. 
Denote by ${\mathcal Cov}_G(\Omega)$ the caregory of stratified 
coverings 
with structure group $G$ ($G$-coverings). A definition of a 
$G$-covering in the case of stratified surfaces possibly with 
boundary is a generalization of a standard one and should be 
explained.

Fix a generic point $x\in\Omega$. Let $\delta\subset\Omega$ 
be a path with ends at $x$. Path $\delta$ is called generic
if it does not cross any special point and cross the boundary 
of $\Omega$ in finitely many points.

Let $\pi:P\to\Omega$ be a stratified covering and $\hat x\in P$ be 
one of preimages of $x$. Lift $\delta$ to $\hat\delta\subset P$ 
using the following 
rule of lifting through a point $\hat y$ such that 
$\pi(y)\in\delta\cap\partial\Omega$:
if the covering $\pi_U:\hat U\to U$ of a 
neighborhood $\hat U$ of point $\hat y$ over 
$U=\pi(\hat U)$ is two-sheeted covering then $\hat\delta$ must 
go from one sheet of $\hat U$ to another sheet through point $\hat y$.

This rule allows to define a fundamental group, a monodromy group etc.
for stratified surfaces and stratified covering. All theory of
fundamental group is consistent for this generalization. Note that these 
definitions are valid for all dimensions in the case of specially
stratified manifolds and specially stratified coverings (see 
definition \ref{defcov}).
(Recall that in this work we deal with specially stratified surfaces
and specially stratified coverings only). 
  
Fix an action of group $G$ on the set of preimages of $x$. A covering
$\pi$ is said to be $G$-covering if (generalized) monodromy group
at $x$ is contained in $G$. All definitions and considerations
of this section, 
particularly, the definition of Hurwitz numbers, can be applied to
$G$-coverings without changes. Thus, we can construct '$G$-Hurwitz
topological field theory'. Structure algebra of it is isomorphic to
the structure algebra of group $G$. $G$-Hurwitz numbers have just 
the same representation as in theorem \ref{th5.4}; these formulas 
coincide with correlators for $G$-Hurwitz topological field theory.

$G$-Hurwits numbers for oriented surfaces without a boundary 
were considered in \cite{Dij}.

\bigskip
A.Alexeevski

Belozersky Institute,
Moscow State University

Moscow 119899 

aba@belozersky.msu.ru

\bigskip

S.Natanzon

Belozersky Institute,
Moscow State University

Moscow 119899

Independent University of Moscow

Institute Theoretical and Experimental Physics 
 
natanzon@mccme.ru


\bigskip

\begin{thebibliography}{99}

\bibitem{Ab} Abrams L, Two-dimensional topological quantum field theory   and frobenius algebras
 / Jornal of Knot and Its Ramifications v.5 No.5 (1995), 569-587 

\bibitem{Al-Gr} Alling N.L., Greenleaf N., Foundations of the theory of
Klein surfaces / Lecture Notes in Math. N 219 // Berlin --
Heideb. -- N.Y.: Springer -- Verlag. -- 1971. p. 117.

\bibitem{Arn} Arnold V., Topological classification of complex trigonometric
polynomials and the combinatorics of graphs with the same number
of edges and vertices, Funct. Anal. and Appl. 30 (1996), 1, 1-17.

\bibitem{At} Atiyah M., Topological Quantum Field Theories, Inst. Hautes
Etudes Sci. Publ. Math., 68 (1988), 175-186.

\bibitem{Bar} Barannikov S.A.,  About space of real polynoms
without multiple critical values, Funktsional. Anal.i Prilozhen.,
26 (2) (1992), 10-17.

\bibitem{Car} Cardy J. L., Boundary conditions, fusion rules and the
Verlinde formula, Nucl. Phys. B 324 (1989), 581-596.

\bibitem{Co-Mo-Ra} Cordes S., Moore G., Ramgoolam S.,
Large $N$ 2D Yang--Mills theory and topological string theory,
Commun. Math. Phys., 185 (1997), 543-619.

\bibitem{Dev} Devadoss S.L., Tessellations of moduli spaces and
the mosaic operad, Contemp. Math. 239 (1999), 91-114.

\bibitem{Dij} Dijkgraaf  R., Mirror symmetry and elliptic curves,
The moduli spaces of curves, Progress in Math., 129 (1995), 149-163,
Brikh\"auser.

\bibitem{Dij1} Dijkgraaf  R., Geometrical Approach to Two-Dimesional Conforal Field Theory, Ph.D.Thesis(Utrecht,1989)
\bibitem{Dub} Dubrovin B., Geometry of $2D$ topological field theories,
in: Lecture Notes in Math., 1620, Springer, Berlin (1996), 120-348.

\bibitem{Eke-Lan} Ekedahl T., Lando S,K., Shapiro M., Vainshtein A.,
On Hurwitz numbers and Hodge integrals, C.R.Acad. Sci. Paris
S/'er. I Math. 328 (1999), 1175-1180.

\bibitem{Gaw} Gawedzki K., Boundary WZW, G/H, G/G and CS theories.
arXiv: hep-th/0108044.

\bibitem{Go-Jac} Goulden I.P., Jackson D.M., Transitive factorizations
into transpositions and holomorphic mappings on the sphere,
Proc. Amer. Math. Soc. 125 (1997), 51-60.

\bibitem{Go-Jac2} Goulden I.P., Jackson D.M., The combinatorial
relationship between trees, cacti and certain connection
coefficients for the symmetric group, Europ. J. Combinatorics, 13
(1992), 357-365.

\bibitem{Gr-Tay} Gross D.J., Taylor IV W., Twists and Wilson loops in the
string theory of two dimensional QCD, Nuclear Phys., B 403
(1993), N 1-2, 395-449.

\bibitem{Hug-Wei}  Hughes B., Weinberger S.,
Surgery and stratified spaces.
Cappell, Sylvain (ed.) et al., Surveys on surgery theory. Vol. 2:  
Princeton, NJ: Princeton University Press. Ann. Math. Stud. 149
(2001), 319-352.

\bibitem{Hur} Hurwitz A., \"Uber Riemann'sche Fl\"achen mit gegeben
Verzweigungspunkten, Math. Ann., (1891), Bn.39, 1-61.

\bibitem{Kar-Mos} Karimipour V., Mostafazadeh A., Lattice
topological field theory on nonorientable surfaces, J. Math. Phys.
38 (1) (1997), 49-66.

\bibitem{Ker} Kerekjarto B., Vorlesungen \"uber Topologie. I.
Fl\"achentopologie. -- Berlin, Springer. -- 1923.

\bibitem{Kee} Keel S., Intersection theory of moduli space of stable
$N$-pointed curves of genus zero, Transactions of the American
Mathematical Society, 330:2 (1992), 545-574.

\bibitem{Kos-Sta-Wyn} Kostov I.K., Staudacher M., Wynter T.,
Complex matrix models and statistics of branched coverings of 2D
surfaces, Commun. Math. Phys. 191 (1998), 283-298.

\bibitem{Laz} Lazaroiu C.I., On the structure of open-closed topological
field theory in two-dimensions, Nucl. Phys. B 603 (2001),
497-530.

\bibitem{Med} Mednykh A.D., Pozdnjakova G.G., On the number of
nonequivalent coverings over compact nonorientable surface, Sib.
Mat. Zb. 27 (1) (1986), 123-131.

\bibitem{Moo1} Moore G., Some comments on branes, G-flux and K-theory,
Int. J. Mod. Phys. A16 (2001), 939-944.

\bibitem{Moo2}  Moore G., D-Branes, RR-Fields and K-Theory,
\\ http://online.itp.ucsb.edu/online/mp01/moore1/

\bibitem{Nat1} Natanzon S.M., Klein surfaces, Russian Math. Surveys,
45:6 (1990), 53-108.

\bibitem{Nat2} Natanzon S.M., Topology of 2-dimensional coverings and
meromorphic functions on real and complex algebraic curves,
Selecta Mathematica, v.12, N 3 (1993), 251-291.

\bibitem{Nat3} Natanzon S.M., Moduli of real algebraic surfaces,
and their superanalogues. Differentials, spinors and Jacobians of
real curves, Russian Math. Surveys, 54:6 (1999), 1091-1147.

\bibitem{Nat-Tur} Natanzon S., Turaev V., Systems of correlators and
solutions of the WDVV equations, Commun. Math. Phys. 196 (1998),
399-410.

\bibitem{Ok-Pan} Okounkov A., Pandharipande R., Gromov--Witten theory,
Hurwitz numbers and matrix models, I. math.AG/0101147.

\bibitem{Seg} Segal G.B., Two--dimensional conformal field theory and
modular functions, Swansea Proceedings, Mathematical Physics,
1988, 22-37.

\bibitem{Tu} Turaev V., Homotopy field theory in dimension 2 and
group-algebras, arXiv:math.QA/9910010

\bibitem{Wit} Witten E., Some geometrical applications of quantum field theory, In IXth International Congress on Mathematical Physics(1988) 77-116. 

\end{thebibliography}
\end{document}